\input amstex
\documentstyle{amsppt}
\magnification=\magstep1

\NoBlackBoxes

\pageheight{9.0truein}
\pagewidth{6.5truein}

\long\def\ignore#1\endignore{#1}

\ignore
\input xy \xyoption{matrix} \xyoption{arrow}
          \xyoption{curve}  \xyoption{frame}
\def\edge{\ar@{-}}
\def\dttdar{\ar@{.>}}
\def\drbl{\save+<0ex,-2ex> \drop{\bullet} \restore}
\def\levelpool#1{\save [0,0]+(-3,3);[0,#1]+(3,-3)
                  **\frm<10pt>{.}\restore}
\def\dashedge{\ar@{--}}

\endignore

\def\la{{\Lambda}}
\def\lamod{\Lambda\text{-}\roman{mod}}
\def\Lamod{\Lambda\text{-}\roman{Mod}}
\def \len{\operatorname{length}} 
\def\AA{{\Bbb A}}
\def\CC{{\Bbb C}}
\def\PP{{\Bbb P}}
\def\SS{{\Bbb S}}

\def\NN{{\Bbb N}}
\def\QQ{{\Bbb Q}}

\def\aut{\operatorname{Aut}}
\def\Aut{\operatorname{Aut}}

\def\pdim{\operatorname{proj\,dim}}
\def\idim{\operatorname{inj\,dim}}

\def\GL{\operatorname{GL}}

\def\Ext{\operatorname{Ext}}
\def\Tor{\operatorname{Tor}}
\def\autlap{\operatorname{Aut}_\la(P)}
\def\Gal{\operatorname{Gal}(K/\Kcirc)}
\def\End{\operatorname{End}}

\def\h{\operatorname{hmg}}

\def\pr{\operatorname{pr}}
\def\trdeg{\operatorname{trdeg}}

\def\A{{\Cal A}} 
\def\B{{\Cal B}} 
\def\C{{\Cal C}}
\def\D{{\Cal D}}
\def\E{{\Cal E}}

\def\h{\text{hmg}}
\def\I{{\Cal I}}

\def\T{{\Cal T}}

\def\Q{{\Cal Q}}

\def\S{{\sigma}}

\def\T{{\Cal T}}
\def\t{{\frak t}}

\def\V{{\frak V}}
\def\Q{{\frak Q}}

\def\phat{\widehat{p}}
\def\qhat{\widehat{q}}
\def\Phat{\widehat{P}}
\def\zhat{\widehat{z}}
\def\Kcirc{{K_\circ}}

\def\modlad{\operatorname{\bold{Mod}}_d(\Lambda)}
\def\mod{\operatorname{\bold{Mod}}}

\def\toptd{\operatorname{\bold{Mod}}^T_d}
\def\top#1{\operatorname{\bold{Mod}}^{#1}_d}
\def\Top{\operatorname{\bold{Mod}}}
\def\primetoptd{\operatorname{\bold{Mod}}^{T'}_d}
\def\laySS{\operatorname{\bold{Mod}}(\SS)}
\def\lay#1{\operatorname{\bold{Mod}}(#1)}
\def\modS{\operatorname{\bold{Mod}}(\S)}
\def\Hom{\operatorname{Hom}}
\def\SchuS{\operatorname{\Schu}(\sigma)}
\def\grasstd{\operatorname{\frak{Grass}}^T_d}
\def\Grfrak{\operatorname{\frak{Gr}}}

\def\autlap{\aut_\Lambda(P)}
\def\gautlap{\operatorname{Gr-Aut}_\Lambda(P)}
\def\unirad{\bigl(\aut_\la(P)\bigr)_u}
\def\grassS{\operatorname{\frak{Grass}}(\S)}
\def\ggrasstd{\operatorname{Gr-\frak{Grass}}^T_d}
\def\grassSS{\operatorname{\frak{Grass}}(\SS)}
\def\ggrassSS{\operatorname{Gr-\frak{Grass}}(\SS)}
\def\ggrassS{\operatorname{Gr-\frak{Grass}}(\S)}
\def\ggrass#1{\operatorname{Gr-\frak{Grass}}(#1)}
\def\grass{\operatorname{\frak{Grass}}}
 
\def\Soc{\operatorname{Soc}}

\def\Schu{\operatorname{Schu}}
\def\GrG{\operatorname{Gr-}G}
\def\GrC{\operatorname{Gr-}C}

\def\grassgraded{{\bf 1}}
\def\codes{{\bf 2}}
\def\GeomII{{\bf 3}} 
\def\GeomIV{{\bf 4}}
\def\CB{{\bf 5}}
\def\CBII{{\bf 6}}
\def\CBS{{\bf 7}}
\def\DW{{\bf 8}}
\def\DSW{{\bf 9}}
\def\DHL{{\bf 10}}
\def\Hil{{\bf 11}}
\def\domino{{\bf 12}}
\def\GeomI{{\bf 13}}
\def\menace{{\bf 14}}
\def\grassI{{\bf 15}}
\def\grassII{{\bf 16}}
\def\grassIII{{\bf 17}}
\def\Kac{{\bf 18}}
\def\Scho{{\bf 19}}
\def\Sma{{\bf 20}}

\topmatter

\title Generic representation theory of quivers with relations 
\endtitle

\rightheadtext{Generic representation theory}

\author E. Babson, B. Huisgen-Zimmermann, and R. Thomas
\endauthor

\address Department of Mathematics, University of California, Davis, CA
95616
\endaddress

\email babson\@math.ucdavis.edu \endemail

\address Department of Mathematics, University of California, Santa
Barbara, CA 93106-3080 \endaddress

\email birge\@math.ucsb.edu \endemail

\address Department of Mathematics, University of Washington, Box 354350,
Seattle, WA 98195-4350 \endaddress

\email thomas\@math.washington.edu \endemail

\thanks The authors were partly supported by grants from the National
Science Foundation. 
\endthanks

\abstract  The irreducible components of varieties parametrizing the
finite dimensional representations of a finite dimensional algebra $\la$ are
explored, with regard to both their geometry and the structure of the
modules they encode.  Provided that the base field is algebraically closed
and has infinite transcendence degree over its prime field, we establish the
existence and uniqueness (not up to isomorphism, but in a strong sense to be
specified) of a {\it generic module\/} for any irreducible component
$\C$, that is, of a module which displays
\underbar{all} categorically defined generic properties of the modules
parametrized by $\C$; the crucial point of the existence statement  --  a
priori almost obvious  --  lies in the description of such a module in a
format accessible to representation-theoretic techniques.  Our approach to
generic modules over path algebras modulo relations, by way of minimal
projective resolutions, is largely constructive.  It is explicit for large
classes of algebras of wild type.  We follow with an investigation of the
properties of such generic modules in terms of quiver and
relations.  The sharpest specific results on all fronts are obtained for
truncated path algebras, that is, for path algebras of quivers modulo ideals
generated by all paths of a fixed length; this class of algebras
extends the comparatively thoroughly studied hereditary
case, displaying many novel features. 
\endabstract

\endtopmatter

\document

\head 1.  Introduction \endhead

Let $\la$ be a basic finite dimensional algebra over an algebraically
closed field $K$.  Tameness of the representation type of $\la$  --  the
only situation in which one can, at least in principle,  meaningfully
classify all finite dimensional representations of $\la$  --  is a
borderline phenomenon.  However, for wild algebras, it is often
possible to obtain a good grasp of the ``bulk" of
$d$-dimensional representations, for any dimension $d$, by understanding
finitely many individual candidates of dimension $\le d$.  The underlying
approach was initiated by Kac in 1982 for the hereditary case, refined by
Schofield in 1992, and extended to arbitrary finitely generated
$K$-algebras by Crawley-Boevey-Schr\"oer in 2002 (\cite{\Kac},
\cite{\Scho}, \cite{\CBS}).  The idea is to explore the ``generic behavior"
of the modules represented by the irreducible components of the affine
variety, $\modlad$, which parametrizes the $d$-dimensional left
$\la$-modules.  

More precisely:  Suppose that $\C$ is a (locally closed) irreducible
subvariety of $\modlad$; for instance, take $\C$ to be an irreducible
component of $\modlad$.  Then a property $(*)$ of modules is said to be
$\C$-{\it generic\/} in case there exists a dense open subset $U$ of $\C$
such that all $\la$-modules parametrized by points in $U$ have property
$(*)$.  As is common, we will, more briefly, refer to the modules {\it in\/}
$U$. 

For instance, due to \cite{\Kac} and \cite{\CBS}, the number of
indecomposable summands of a module is $\C$-generic for any irreducible
component $\C$ of
$\modlad$, as is the family of corresponding dimension vectors of the
indecomposable summands.  To date, these numerical data, next to $\Ext$- and
$\Hom$-space dimensions and numbers of subrepresentations, have been the
main objects of study along this line (see also \cite{\DW} for the canonical
decomposition, and \cite{\CBII}, \cite{\DSW} for subrepresentations).  

The primary purpose of this paper is to more
broadly study generic properties of the modules in the irreducible 
components of $\modlad$  --  these properties include the generic
behavior of their syzygies  --  first in general (Sections 2-4), then in a
more specialized setting (Section 5).  In intuitive terms, the goal of such
an investigation is to obtain structural information on a substantial part
of the
$d$-dimensional representations, irrespective of the representation type
of the underlying algebra.  Here ``substantial part" means ``a
Zariski-dense open set's worth".

The foundation consists of an existence and
uniqueness result (Theorem 4.3).  Roughly, the
existence part says that, provided the base field $K$ has infinite
transcendence degree over its prime field, {\it all
\/} categorically defined generic properties of an irreducible component
$\C$ of $\modlad$ are on display in a single module $G = G(\C)$ in $\C$;
we will explain in a moment what we mean by ``categorically defined". 
Existence of $G$ is essentially obvious, unless one insists (as we do) on a
specific format permitting structural and homological evaluation.  Our
construction provides a minimal projective presentation of
$G$, which explicitly depends on a comparatively small
subvariety of $\C$: namely, on the subvariety $\E$ whose points
correspond to those modules in
$\C$ that share a type of ``path basis" (a basis of this ilk will be called
a {\it skeleton\/}  --  see below).  Two advantages of these varieties
$\E$ are: (a) they are readily obtained from a presentation of $\la$ by
quiver and relations, and (b), they allow for an effective back and forth
between their points and the first syzygies of the modules they encode.
Along this line, we obtain, for any irreducible component $\C$ of
$\modlad$, a module $G$ in the intersection
of $\C$ with a suitable $\E$, such that $G$ has all
$\C$-generic properties satisfying a mild invariance condition; namely
we ask that the pertinent module properties be invariant under all Morita
self-equivalences of the following restricted type.  Suppose that
$\Kcirc$ is an algebraically closed subfield of $K$ with
$\trdeg(K:\Kcirc) \ge \aleph_0$ such that $\la$ is defined over $K_\circ$,
i.e.,
$\la = KQ/I$, where $Q$ is a quiver and $I$ an ideal defined over
$\Kcirc$.  The Morita self-equivalences of $\lamod$ under consideration are
those which are induced by the $\Kcirc$-automorphisms of
$K$; we dub them $\Gal$-{\it equivalences\/}, and refer to a module property
which is invariant under such equivalences as $\Gal$-{\it invariant\/}.  
Calling a module $G$ in $\C$ {\it generic for $\C$\/} if it has all
$\Gal$-invariant generic properties of $\C$, we moreover prove the
generic modules for
$\C$ to be unique up to $\Gal$-equivalence.  More generally, we establish
a uniqueness result for generic orbits of
$K$-varieties which are defined over $\Kcirc$ and which carry, next to an
algebraic group action,  a suitably compatible action of $\Gal$. 
(We point out that the above concept
``generic for $\C$", or ``$\C$-generic" more briefly, must not be confused
with ``generic" in the sense of Crawley-Boevey, as defined in \cite{\CB}; in
particular, the generic objects considered here are finite dimensional over
$K$.)
Our general description of $\C$-generic modules is essentially
constructive, provided that $\la$ is given by quiver and relations.  It is
representation-theoretically manageable in case $\C$ is rational. 

With regard to a specific class of algebras, the general
background suggests a program:  Namely, to gain a better
understanding of the geometry of the irreducible components
$\C$ of the varieties $\modlad$  -- or, often far easier, the geometry of
the varieties $\E$ -- in order to more effectively explore the
representation-theoretic properties of their generic modules.  A sampling
of the generic properties to be addressed can be found in Corollary 4.7; in
particular, we re-encounter results from \cite{\CBS} mentioned above,
addressing direct sum decompositions of generic modules.    For our
approach to generic modules, it is crucial that the list of $\Gal$-stable
generic properties includes {\it skeleta\/} of modules (see Definitions
3.1).  These are preferred bases reflecting the radical layering and the
$KQ$-structure of a module $M$.  In light of the key role they play towards
useful projective presentations of generic modules, we precede their formal
introduction with a preliminary one.  Suppose that
$\la = KQ/I$, and let $M  = P/C$ be a
$\la$-module with projective cover $P$; moreover, fix elements $z_1,
\dots, z_t$ of $P$ which induce a basis for
$P/JP$.  A skeleton of $M$ is a set $\sigma$ of
elements of the form $y_{p,r} = (p + I) z_r \in P$, where $p$ runs through
certain paths in $KQ$.  Keeping track of the lengths of these paths, we
require that the following two conditions be satisfied by the $y_{p,r}$
in $\S$:
$\bullet$ the residue classes $y_{p,r} + C$ in $M$, corresponding to
paths $p$ of any fixed length $l$, induce a basis for $J^l M / J^{l+1} M$;
moreover, 
$\bullet$ the set $\S$ is closed under initial subpaths, that is, if
$y_{p,r}$ belongs to $\S$ and $p'$ is an initial subpath of $p$, then
also $y_{p', r}$ belongs to $\S$.  The purpose of the second requirement 
may not be immediately apparent.  It will turn out to be pivotal towards the
usefulness of the varieties
$\E$ from which our construction is launched (see Section 3.C).  We note
that every 
$\la$-module has a nonempty set of skeleta, and that the set of all skeleta
of the modules sharing a given radical layering is finite (as long as
$P$ and the $z_r$ are fixed).  In Reduction Step 3 below, we will sketch how
skeleta are tied into the search for explicit presentations of generic
modules. 
\medskip    

Our second major goal (Section 5) is to carry out the suggested program for
{\it truncated path algebras}, that is, for the algebras of the form
$KQ/ I$, where $Q$ is a quiver and
$I$ the ideal generated by all paths of a fixed length.  They include
the algebras with vanishing radical square and the hereditary algebras. 
Clearly, any basic finite dimensional algebra is isomorphic to a factor
algebra of a truncated one.  In fact, if $\la =  KQ/I$ is any algebra of
Loewy length
$L+1$ and $\la' = KQ/\langle \text{all paths of length}\ L+1 \rangle$,
then the parametrizing varieties for
$\lamod$ are closed subvarieties of the analogous varieties for
$\la'$-mod.  The case where $\la$ is truncated demonstrates, in particular,
the efficiency of the format in which our key Theorem 4.3 displays the first
syzygy $\Omega^1(G)$ of a generic module $G$:  For an arbitrary finite
dimensional factor algebra $\la$ of a path algebra, $\Omega^1(G)$ is given
in terms of a finite generating set $\bigl(g_i \bigr)_{ i \in \I}$ depending
on quiver, relations, and the considered skeleton; when $\la$ is a
truncated path algebra, this choice is irredundant.  Indeed, it then
provides us with a direct sum decomposition $\Omega^1(G) = \bigoplus_{i \in
\I} \la g_i$ in which all summands are nontrivial. 
\medskip

To provide better orientation, we outline
our initial three reduction steps; two of them are presented in Section 2,
and the final one in Section 3.  These reductions underlie our
construction of generic objects in Section 4.  Through Section 4,
we do not impose any condition on
$\la$, beyond the above assumption on the transcendence degree of $K$.  We
then preview some results from Section 5.  This last section fleshes out
the general theory in the special case of truncated path algebras.   

\subhead  The general case \endsubhead
\smallskip

Fix an irreducible component $\C$ of some variety
$\modlad$.  Our strategy for accessing $\C$-generic modules involves
reductions of the problem to successively smaller subvarieties of
$\C$.  That the subvarieties $\C' \supseteq \C'' \supseteq C'''$
we choose are better adapted to our purpose is only to a minor degree due
to the reduction in size.  Their main benefit lies in the availability of
alternate, more helpful varieties parametrizing the classes of modules
corresponding to the points of these subvarieties.  The largest, $\D'$, of
the alternates, 
$\D' \supseteq \D'' \supseteq \D'''$, is a projective variety.  All of the
reductions are in representation-theoretic terms, that is, they
are based on successively finer isomorphism invariants of the modules $M$
they encode:  First the {\it top\/}
$T(M) = M/JM$, then the {\it radical layering\/}
$\SS(M) = \bigl(J^l M/ J^{l+1} M \bigr)_{0 \le l \le L}$, where
$L$ is maximal with $J^L \ne 0$ (we identify the
semisimple modules in this sequence with their isomorphism classes), and
finally the {\it skeleta\/} $\S$ of $M$.    
\smallskip   

{\bf Reduction Step 1} (Subsection 2.A).   Given a semisimple
module $T$, we let $\toptd$ be the locally closed subvariety of
$\modlad$ consisting of the points which correspond to the modules with top
(isomorphic to) $T$.  Clearly,
$\modlad$ is the disjoint union of the subvarieties $\toptd$, where $T$
runs through finitely many choices. Since, generically, the modules in $\C$
have fixed top, there exist precisely one semisimple $T$ and
precisely one irreducible component
$\C'$ of $\toptd$ such that $\C$ is the closure of $\C'$ in
$\modlad$.  In particular, any 
$\C'$-generic module is also $\C$-generic.

The projective counterpart of the variety $\toptd$ (presented at the end of
2.A) is denoted by $\grasstd$.  It is the obvious closed subvariety of the
classical Grassmannian of all $(\dim P - d)$-dimensional $K$-subspaces of
$JP$, where $P$ is a fixed projective cover of $T$; this variety was
introduced and geometrically related to $\toptd$ by Bongartz and the author
in \cite{\GeomII} and \cite{\GeomIV}.  Since the irreducible components of
$\grasstd$ are in natural bijection with those of
$\toptd$  --  see Proposition 2.2 below  -- this means that, in studying the
$\C$-generic modules, we may restrict our focus to an irreducible 
component $\D'$ of $\grasstd$. 
\smallskip 

{\bf Reduction Step 2}  (Subsection 2.B).  Given a sequence $\SS = (\SS_0,
\SS_1, \dots, \SS_L)$ of semisimple modules that has total dimension $d$,
we let
$\laySS$ be the subvariety of $\modlad$ consisting of the points
corresponding to the modules $M$ with radical layering $\SS$.  Clearly,
each of the varieties $\toptd$ is the disjoint union of the subvarieties
$\laySS$, where $\SS$ traces those semisimple sequences for
which
$\SS_0 = T$.  Again, it is readily seen that, generically, the modules
in $\C$ have fixed radical layering, say $\SS$.  As a consequence, there
exists precisely one irreducible component $\C''$ of $\laySS$ with the
property that the closure of
$\C''$ in $\toptd$ coincides with $\C'$.  Consequently, the closure of
$\C''$ in $\modlad$ equals $\C$, whence any $\C''$-generic module
is also $\C$-generic. 

The counterpart of the subvariety $\laySS$ of $\toptd$ is the
subvariety
$\grassSS$ of $\grasstd$ consisting of the projective points which
parametrize the modules with radical layering $\SS$.  Since the irreducible
components of $\laySS$ are in natural bijection with those of $\grassSS$,
our task is thus reduced to the study of $\D''$-generic
modules, where $\D''$ is the irreducible component of $\grassSS$
corresponding to $\C''$.
\smallskip

{\bf Reduction Step 3} (Section 3).  This step relies on a presentation
$\la = KQ/I$, where $Q$ is a quiver and $I \subset KQ$ an admissible
ideal.  For any sequence $\SS$ of semisimple modules, the variety $\grassSS$
is a finite union of {\it open\/} affine
varieties
$\grassS$, where $\S$ runs through the skeleta ``compatible with $\SS$" in
an obvious sense, and $\grassS$ consists of those points in $\grassS$ which
correspond to the modules with skeleton $\S$.  Thus the closures of the
irreducible components of the $\grassS$ in $\grassSS$ are precisely the
irreducible components of $\grassSS$.  In other words, the generic modules
for the irreducible components of the varieties $\grassSS$ coincide with the
generic modules for the irreducible components of the $\grassS$.  So we are
left with the task of constructing generic modules for the irreducible
components of the  $\grassS$.  (The subvarieties of the
classical $\modlad$ corresponding to the latter are the varieties $\E$ to
which we referred in the preceding outline.)  On this level, useful
descriptions of generic objects are within reach.
\smallskip

In particular, the generic modules for the irreducible components of the
$\laySS$ include the generics for the irreducible components of
$\modlad$.  The ``redundant" generic modules on the resulting list  -- 
those that are generic on the $\laySS$-level, but not
generic for any irreducible component of the ambient variety
$\modlad$  --  are of interest in their own right:  They yield a more
complete generic picture of the representation theory of
$\la$ than restriction to the components of $\modlad$ would. 
On the other hand, the sifting required to reduce this larger list to the
modules which are generic on the $\modlad$-level is
not constructive in general, as far as we can tell.  In the present paper,
we only carry it out in examples, but will address it more systematically in
a sequel.  If one is solely interested in the generic modules for the
irreducible components of $\modlad$, our approach amounts to a tradeoff: 
We consider ``too many" generic objects, but obtain them in a useful
format.

\subhead Truncated path algebras \endsubhead
\smallskip

In the truncated situation, we may take $\Kcirc$ to be the algebraic closure
of the prime field of $K$.  The following geometric information paves the
way to explicit construction and analysis of the generic modules for the
irreducible components of the varieties
$\laySS$.       

\proclaim{Theorem A} Suppose $\la = KQ/I$ is a truncated path algebra
of Loewy length $L+1$.  For each sequence $\SS$ of semisimple
$\la$-modules, $\grassSS$ is covered by
dense open subsets {\rm {(}}the $\grassS${\rm {)}} which are isomorphic
copies of affine $N$-space $\AA^N$, where $N$ depends only on $\SS$.  In
particular,
$\grassSS$ is irreducible, rational and smooth.   

In more detail,
$\grassSS$ is an affine bundle with fibre $\AA^{N_1}$ over $\ggrassSS$,
where $\ggrassSS$ is the subvariety of $\grassSS$ consisting of the points
that correspond to graded modules generated in degree zero; both $N$ and
$N_1$ can be determined from the quiver $Q$ and the integer $L$ by a simple
count. Moreover, the variety
$\ggrassSS$ of graded objects is projective and, in turn, smooth and
rational: It is an iterated Grassmann bundle over a finite
direct product of classical Grassmannians. {\rm(For terminology, consult
Theorems 5.3, 5.9, and the paragraph preceding the latter.)}  
\endproclaim    

We follow with a slice of the representation-theoretic side of the
picture, stated somewhat informally.  For more precision, we refer to
Theorem 5.12, which addresses the {\it generic graded modules\/} in
tandem with the generic modules as originally defined.  

\proclaim{Theorem B} Keep the hypotheses of Theorem A.

The generic modules $G = G(\SS)$ for
$\laySS$, alias $\grassSS$, can be read off the quiver $Q$, as can be
several of their algebraic invariants. 

In particular, the {\rm {(}}generic{\rm {)}} skeleta of $G$
are available at a glance, and from those the generic syzygies.  
The syzygies $\Omega^k(G)$, for $k \ge 1$, are direct sums
of cyclic modules, which are determined by $\SS$ up to isomorphism {\rm
{(}}not only up to $\Gal$-equivalence{\rm {)}}.
\endproclaim

By means of this theorem, the generic homological dimensions of the
modules with fixed radical layering have been determined, in terms of
the quiver $Q$ and the cutoff length $L$ (see \cite{\DHL}).
A suitable simultaneous choice of the modules $G(\SS)$, where $\SS$ runs
through the eligible sequences of semisimple modules, additionally
yields the following ``relative" generic behavior:  Given any two distinct
semisimple sequences $\SS$ and $\SS'$, the pair $(G(\SS), G(\SS'))$
possesses any
$\Gal$-stable generic property of pairs of modules in $\grassSS
\times \grass{\SS'}$.   More concretely, this entails for instance that
$$\dim \Ext^i_\la \bigl(G(\SS), G(\SS') \bigr) = \min\{\dim
\Ext^i_\la(M,N) \mid \SS(M) = \SS,\ \SS(N) = \SS' \} \ \ \ \text{for} \ i
\ge 0$$(see Corollary 4.7). 
The two examples following Theorem 5.12 demonstrate the strength of this
result.  For illustrations of the more general Theorem 4.3 providing the
backdrop for Theorem B, see Example  4.8.
\smallskip  

The crucial concepts can be found in Definitions 3.1 (skeleta), 3.7
(critical paths and affine coordinates of the $\grassS$), 3.9 (hypergraphs
of modules), 4.1 and 4.2 (generic modules for the irreducible components of
the parametrizing varieties).  Moreover, Theorem 3.8 summarizes background
from
\cite{\grassIII}.

\head 2. Playing irreducible components back and forth among subvarieties
of $\modlad$ and their projective counterparts
\endhead  

We briefly review some of the constructions and results from
\cite{\grassI}, \cite{\grassII}, and \cite{\grassIII} needed in the
sequel.  Along the way, we line up observations aimed at reducing the
problem of finding the irreducible components of
$\modlad$ to finding the irreducible components of successively smaller
varieties, first $\toptd$, next $\laySS$.  The full collections of
irreducible components of $\toptd$ and
$\laySS$ are of interest in themselves towards refined pictures of the
generic representation theory of $\la$.

As already stated, we let $\la$ be a basic finite dimensional algebra over
an algebraically closed field $K$.  Hence, we may assume, without loss of
generality, that $\la = KQ/I$, where $Q$ is a quiver and
$I$ an admissible ideal in the path algebra $KQ$.  The vertices $e_1,
\dots, e_n$ of $Q$ will be identified with the primitive idempotents of
$\la$ corresponding to  the paths of length zero.  As is well-known, the
left ideals $\la e_i$ then represent all indecomposable projective (left)
$\la$-modules, up to isomorphism, and the factors $S_i = \la e_i /J e_i$,
where $J$ is the Jacobson radical of
$\la$, form a set of representatives for the simple (left)
$\la$-modules.  By
$L+1$ we will denote the {\it Loewy length\/} of $J$, that is, $L$ is
maximal with $J^L \ne 0$.   Moreover, we will observe the following
conventions:  The product $pq$ of two paths $p$ and $q$ in $KQ$  stands for
``first
$q$, then $p$"; in particular, $pq$ is zero unless the end point of $q$
coincides with the starting point of $p$.  In accordance with this
convention, we call a path $p_1$ an {\it initial subpath\/} of $p$ if $p =
p_2 p_1$ for some path
$p_2$.  A {\it path in
$\la$\/} is a residue class of the form $p + I$, where
$p$ is a path in $K Q \setminus I$; we will suppress the residue notation,
provided there is no risk of ambiguity.  Further, we will gloss over the
distinction between the left $\la$-structure of a module
$M \in \Lamod$ and its induced $KQ$-module structure when there is no
danger of confusion.  An element
$x$ of $M$ will be called a {\it top element\/} of
$M$ if $x \notin JM$ and $x$ is {\it normed\/} by some $e_i$, meaning that
$x = e_i x$.  Any collection $x_1, \dots, x_m$ of top elements of $M$
generating
$M$ and linearly independent modulo $JM$ will be called a {\it full
sequence of top elements of $M$\/}.

The two isomorphism invariants of a $\la$-module $M$ which will be pivotal
here are the {\it top\/} and the {\it radical layering\/} of $M$, the
latter being a refinement of the former.  The {\it top\/} of $M$ is defined
as
$M/JM$,  and the {\it radical layering\/} as the sequence
$\SS(M) = \bigl(J^lM / J^{l+1}M \bigr)_{1\le l \le L}$ of semisimple
modules.  We will identify  isomorphic semisimple modules and, in
particular, call any module $M$ with
$M/JM  \cong T$ a {\it module with top $T$\/}.    
\medskip

\noindent {\it The following choices and notation will be observed
throughout\/}: We fix a semisimple module $T$, say 
$$T = \bigoplus_{1 \le i \le n} S_i^{t_i},$$   set $\t = \sum_i t_i= \dim
T$, and denote by
$$P = \bigoplus_{1 \le i \le n} (\la e_i)^{t_i}$$ ``the" projective cover
of $T$.  Clearly, $P$ is also a projective cover of any module with top $T$;
in other words, the modules with top $T$ are precisely the quotients $P/C$
with $C \subseteq JP$, up to isomorphism.  Next, we fix a full  sequence
$z_1, \dots, z_\t$ of top elements of $P$.  This means $P = \bigoplus_{1
\le r \le \t} \la z_r$ with $\la z_r \cong \la e(r)$, where $e(r)$ is the
idempotent in $\{e_1, \dots, e_n\}$ norming $z_r$.   A natural choice of
such top elements of $P$ is to take the $z_r$ to be the  primitive
idempotents $e_i$, each with multiplicity
$t_i$, distinguished by their ``slots'' in the above decomposition of
$P$.  Finally, we fix a positive integer $d \ge \t$.

In the following, we will refer to $P$ as {\it the distinguished projective
cover of $T$\/} with {\it distinguished sequence $z_1, \dots, z_\t$ of top
elements\/}.

The hierarchies on sets of irreducible components which we introduce in
Observation 2.1(2) and Corollary 2.7 are subsidiary to our theoretical
development and will resurface only in examples.

\subhead 2.A. From the affine variety $\modlad$ to the quasi-affine variety
$\toptd$ and the projective variety $\grasstd$
\endsubhead
 
Let $a_1, \dots, a_r$ be a set of algebra generators for
$\la$ over $K$.  A convenient set of such generators consists of the
primitive idempotents (= vertices)
$e_1, \dots, e_n$ together with the (residue classes in $\la$ of the)
arrows in $Q$.  Recall that, for $d \in \NN$, the classical affine variety
of $d$-dimensional representations of
$\la$ can be described in the form
$$\modlad = \{(x_i) \in \prod_{1 \le i \le r} \End_K(K^d)
\mid \text{\ the\ } x_i \text{\ satisfy all relations satisfied by the
\ } a_i\}.$$
  As is well-known, the isomorphism classes of $d$-dimensional (left)
$\la$-modules are in one-to-one correspondence with the orbits of
$\modlad$ under the $\GL_d$-conjugation action.  Moreover, the connected
components of $\modlad$ are in natural bijection with the
dimension vectors ${\bold d} = (d_1, \dots, d_n)$ such that
$\sum_i d_i = d$; by $\mod_{\bold d}(\la)$ we denote the connected
component corresponding to ${\bold d}$.  If
$I = 0$, that is, if
$\la$ is hereditary, the connected components coincide with the irreducible
components, but this fails already in small non-hereditary examples, e.g.
for $d = 2$ and 
$\la = KQ/I$, where $Q$ is the quiver 
\ignore
$\xymatrixrowsep{1pc}\xymatrixcolsep{2pc}
\xymatrix{   1 \ar@/^0.2pc/[r] &2 \ar@/^0.2pc/[l] }$\endignore and $I$ is
generated by the paths of length $2$.   

The tops lead to a first rough subdivision of $\modlad$:  By
$\toptd$ we will denote the locally closed subvariety of $\modlad$ which
consists of the points representing the modules with top
$T$.  As is easily seen,
$\toptd$ is nonempty if and only if
$\dim T \le d$ and the projective cover $P$ of $T$ has dimension at least
$d$. In light of our identification of isomorphic semisimple modules, the
set of all finite dimensional semisimple modules is partially ordered
under inclusion. Specifically, we write $T \le T'$ to denote that a
semisimple module $T$ is (isomorphic to) a submodule of a semisimple module
$T'$. 
   
The first part of the following observation is due to the fact that
$\modlad$ is the finite disjoint union of the locally closed subvarieties
$\toptd$, where $T$ runs through the semisimple modules that arise as tops
of $d$-dimensional $\la$-modules.  The second part is an immediate
consequence of upper semicontinuity of the functions $\modlad
\rightarrow \NN$ given by 
$X \mapsto \dim \Hom_\la(X, S_i)$ for $1 \le i \le n$.

\proclaim{Observation 2.1 and Terminology} 

{\rm (1)} For every irreducible component $\C$ of $\modlad$, there exist
precisely one semisimple module $T$ and precisely one irreducible component
$\D$ of $\toptd$ such that
$\C$ is the closure of $\D$ in
$\modlad$.  In particular: Generically, the modules in $\C$ have top
$T$.  

{\rm (2)}  If $\toptd \ne \varnothing$ and $T$ is minimal among the
semisimple modules
$T'$ that give rise to nonempty sets $\primetoptd$, then $\toptd$ is open
in $\modlad$.  Consequently: If the minimal elements in the poset of $T'$ 
with $\primetoptd \ne \varnothing$ are
$T^{(0,1)}, \dots, T^{(0,m_0)}$, then the closures of the varieties
$\top{T^{(0,i)}}$, $1 \le i \le m_0$, are unions of irreducible components
of $\modlad$.  In fact, the irreducible components of the closures
$\overline{\top{T^{(0,i)}}}$  --  they are called the {\rm {irreducible
components of class $0$ of $\modlad$}}  --  are in bijective correspondence
with the irreducible components of $\top{T^{(0,i)}}$.    
 
To continue recursively, suppose $T^{(h,1)}, \dots, T^{(h,m_h)}$ are the
distinct semisimple modules which are minimal in the poset of those $T'$
for which 
$$\top{T'}\  \not\subseteq \ \bigcup_{k < h \text{\ and\ } i
\le m_k}  \overline{\top{T^{(k,i)}}}.$$   Then each $\top{T^{(h,i)}}$
intersects the closed subvariety $\modlad
\setminus \bigcup_{k < h, i \le m_k} \top{T^{(k,i)}}$ in a nonempty open
set.  In particular, the irreducible components of the various subvarieties
$\top{T^{(h,j)}}$ which are not contained in the closure of
$\, \bigcup_{k < h,\, i \le m_k} \top{T^{(k,i)}}$ yield distinct
irreducible components of $\modlad$, via passage to closures.  These are
called the {\rm {irreducible components of class
$h$ of $\modlad$}}. \qed 
\endproclaim

Observation 2.1 says that the problem of identifying the irreducible
components of
$\modlad$ and exploring their general modules can be played back to the
components of the subvarieties $\toptd$.  So, in studying the components of
$\modlad$, one of the first questions addresses the semisimple modules that
actually arise as ``generic tops".  (To give an easy example: if
$\la = KQ$ where $Q$ is the quiver $1 \rightarrow 2$, then $S_1$ is a
generic top for the $2$-dimensional modules  -- it gives rise to an
irreducible component of class $0$  --  as is
$S_1 \oplus S_1$  --  the latter gives rise to a component of class $1$ in
our partial order  --  whereas $S_1 \oplus S_2$ fails to be a generic top
for $\mod_2(\la)$.) 

In the next subsection, the partial order on tops will be extended to one
on sequences of semisimple modules which, again, are generic invariants of
the irreducible components of $\modlad$.   
\medskip

Our principal tool will be an alternate variety parametrizing the same
class of representations as
$\toptd$.  Namely, we consider the following closed subvariety of the
classical Grassmannian $\Grfrak(\tilde{d},JP)$ of  $\tilde{d}$-dimensional
subspaces of the $K$-space $JP$, where $\tilde{d} = \dim_K P - d$:
$$\grasstd = \{C \in \Grfrak(\tilde{d},JP) \mid C \text{\ is a\ }
\la\text{-submodule of\ } JP \}.$$  This variety comes with an obvious
surjection
$$\grasstd \longrightarrow \{\text{isomorphism classes of\ } d\text
{-dimensional modules with top\ } T\},$$   sending $C$ to the class of
$P/C$.  Clearly, the fibres of this map  coincide with the orbits of the
natural $\autlap$-action on
$\grasstd$.   While the global geometry of the projective variety
$\grasstd$ cannot be reasonably compared with that of the quasi-affine
variety $\toptd$, the ``relative geometry'' of the
$\autlap$-stable subsets of $\grasstd$ is tightly linked to that of the
$\GL_d$-stable subsets of $\toptd$ in the following sense:  

\proclaim{Proposition 2.2} {\rm{(}}See \cite{\GeomIV, Proposition
C}.{\rm{)}} The assignment $\autlap.C \mapsto \GL_d.x$, which pairs orbits
$\autlap.C \subseteq \grasstd$ and
$\GL_d.x \subseteq \toptd$  representing the same
$\la$-module up to isomorphism, induces an inclusion-preserving bijection
$$\Psi: \{ \autlap\text{-stable subsets of\ } \grasstd \} \rightarrow
\{\GL_d\text{-stable subsets of\ } \toptd \}$$  which preserves openness,
closures, connectedness, irreducibility, and types of singularities.
\qed
\endproclaim
  
This correspondence permits transfer of information concerning the
irreducible components of any locally closed
$\GL_d$-stable subvariety of $\toptd$ to the irreducible components of the
corresponding $\autlap$-stable subvariety of $\grasstd$, and vice versa. 
Indeed, since the acting groups, $\GL_d$ and
$\autlap$, are connected, all of their orbits are irreducible.  Therefore
all such irreducible components are again stable under the respective
actions.  In particular, we obtain:

\proclaim{Observation 2.3}  The irreducible components of
$\toptd$ are in natural one-to-one correspondence with the irreducible
components of $\grasstd$.  Moreover, if ${\bold d} = (d_1, \dots, d_n)$
is a partition of $d$ and
$\grass^T_{{\bold d}}$\ the closed subvariety of
$\grasstd$ consisting of the points that represent modules with dimension
vector ${\bold d}$, then $\grass^T_{{\bold d}}$ is a union of
irreducible components of $\grasstd$.  \qed
\endproclaim

Finding the irreducible components of $\toptd$ by way of the alternate
projective setting offers significant advantages, as we will see.

\subhead 2.B.  From $\toptd$ and $\grasstd$ to $\laySS$ and
$\grassSS$ \endsubhead      

The radical layering of modules provides us with a further partition of
$\grasstd$ into pairwise disjoint locally closed subvarieties.  (In
general, this partition fails to be a stratification in the technical
sense, however, even when $\la$ is hereditary.)  

A $d$-{\it dimensional semisimple sequence $\SS$ with top\/}
$T$ is any sequence  
$(\SS_0, \dots, \SS_L)$ of semisimple modules such that
$\SS_0 = T$ and
$\sum_{0 \le l \le L} \dim \SS_l = d$.   Since we identify semisimple
modules with their isomorphism classes, a semisimple sequence amounts to a
matrix of discrete invariants keeping count of the multiplicities of the
simple modules in the semisimples occurring in the slots of
$\SS$.  

Accordingly, we consider the following action-stable locally closed
subvarieties of
$\grasstd$ and $\toptd$, respectively:  
$$\grassSS  = \{ C \in \grasstd \mid \SS(P/C) = \SS\},$$ while $\laySS
\subseteq \toptd$ consists of those points in 
$\toptd$ that parametrize the modules with radical layering $\SS$. Clearly,
the above one-to-one correspondence between the $\autlap$-stable subsets of
$\grasstd$ and the $\GL_d$-stable subsets of $\toptd$ restricts to a
correspondence between the $\autlap$-stable subsets of $\grassSS$ and the 
$\GL_d$-stable subsets of $\laySS$.  Thus we obtain the following
counterpart to Observation 2.1(1):

\proclaim{Observation 2.4}  For each irreducible component $\C$ of
$\grasstd$ {\rm {(}}resp., $\toptd${\rm {)}}, there exist precisely one
$d$-dimensional semisimple sequence $\SS$ and precisely one irreducible
component $\D$ of $\grassSS$ {\rm {(}}resp., $\laySS${\rm {)}} such that
$\C$ equals the closure of $\D$ in $\grasstd$ {\rm {(}}resp.,
$\toptd${\rm {)}}.  In particular: Generically, the modules in $\C$ have
radical layering $\SS$.  \qed
\endproclaim

In parallel to Section 2.A, we next introduce a partial order on the
(finite) set of $d$-dimensional semisimple sequences.  It refines the
partial order for semisimple modules in that
$\SS \le \SS'$ implies that $\SS_0
\le \SS'_0$.

\definition{Definition 2.5} Let $\SS$ and $\SS'$ be two semisimple
sequences with the same dimension.  We say that
$\SS'$ {\it dominates\/} $\SS$ and write $\SS \le \SS'$ if and only if
$\bigoplus_{l \le r} \SS_l \le \bigoplus_{l \le r}
\SS'_l$ for all $r \ge 0$. 
\enddefinition 

Roughly speaking, $\SS'$ dominates $\SS$ if and only if $\SS'$ results from
$\SS$ through a finite sequence of leftward shifts of simple summands of
$\bigoplus_{0 \le l \le L} \SS_l$ in the layering provided by the
$\SS_l$. 

In intuitive terms, the next observation says that the simple summands in
the radical layers of the modules represented by $\grassSS$ are only
``upwardly mobile" as one passes to modules in the boundary of the closure
$\overline{\grassSS}$ of
$\grassSS$ in $\grasstd$.  All of the semisimple sequences are tacitly
assumed to be
$d$-dimensional.   

\proclaim{Observation 2.6}  {\rm{(}}For a proof, see
\cite{\grassII, Observation 3.4}.{\rm{)}} Suppose that
$\SS$ is a semisimple sequence with top $T$.  Then the union
$\bigcup_{\SS' \ge \SS, \,\SS'_0 = T} \grass{(\SS')}$ is closed in
$\grasstd$. {\rm {(}}As before, $\SS'_0$ is the semisimple module in the
$0$-th slot of the sequence $\SS'$.{\rm {)}}\qed
\endproclaim

In view of Proposition 2.2, the following corollary translates into an
observation about the classical variety 
$\toptd$, yielding a counterpart to Observation 2.1(2).  

\proclaim{Corollary 2.7 and Terminology}  All semisimple sequences are
assumed to be $d$-dimen\-sio\-nal with top $T$. 

If $\SS$ is minimal among the semisimple sequences $\SS'$ with top $T$ that
give rise to nonempty sets
$\grass{(\SS')}$, then $\grassSS$ is open in $\grasstd$.  Thus: If the
minimal elements in the poset of such sequences $\SS'$ are $\SS^{(0,1)},
\dots, \SS^{(0,n_0)}$, then the closures of the varieties
$\grass{\SS^{(0,i)}}$ in $\grasstd$, $1 \le i \le n_0$, are unions of
irreducible components of $\grasstd$.  In fact, the irreducible components
of $\overline{\grass{\SS^{(0,i)}}}$   --  they are  called {\rm {the
irreducible components of class $0$ of $\grasstd$}}  --  are in bijective
correspondence with those of
$\grass{\SS^{(0,i)}}$.  
 
To continue recursively, suppose $\SS^{(h,1)}, \dots, \SS^{(h,n_h)}$ are
the distinct semisimple sequences  which are minimal in the poset of
semisimple sequences $\SS'$ with 
$$\grass{(\SS')} \ \not\subseteq \ \bigcup_{k < h,\, i \le n_k}
\overline{\grass{\SS^{(k,i)}}}.$$   Then each $\grass{\SS^{(h,i)}}$
intersects the closed subvariety
$\grasstd \setminus \bigcup_{k < h,\, i \le n_k} \grass{\SS^{(k,i)}}$ in a
nonempty open set.  In particular, the irreducible components of the
various subvarieties $\grass{\SS^{(h,j)}}$ which are not contained in the
closure of $\bigcup_{k < h,\, i \le n_k} \grass{\SS^{(k,i)}}$ yield
distinct irreducible components of $\grasstd$, via passage to closures. 
They are called the {\rm {irreducible components of class $h$ of
$\grasstd$}}. 
\endproclaim

\demo{Proof}  For the first assertion, suppose that $\SS$ is minimal among
the semisimple sequences $\SS'$ with $\grass{(\SS')} \ne
\varnothing$.  Let
$\A$ be the collection of semisimple sequences that lie minimally above
$\SS$, and $\B$ the collection of those sequences which are not comparable
with $\SS$.  Moreover, let $\widetilde{\A}$ and
$\widetilde{\B}$ be the union of all $\grass(\SS')$ for those semisimple
sequences $\SS'$ that are larger than or equal to some member of $\A$ and
$\B$, respectively.  Then $\widetilde{\A}$ and $\widetilde{\B}$ are closed,
since the set of semisimple sequences (with top $T$) is finite.  The claim
thus follows from the equality
$\grasstd \setminus \grassSS = \widetilde{\A} \cup \widetilde{\B}$.

The other assertions are proved similarly. \qed
\enddemo

Provided we understand the irreducible components of the
$\grassSS$ and their closures in $\grasstd$, we can therefore find the
irreducible components of the latter variety.  Hence, the problem of
finding the generic semisimple sequences for $\grasstd$ imposes itself,
namely of finding those $\SS$ which generically arise as radical layerings
of the modules in the various irreducible components of $\grasstd$.

Finally, we point out that our findings in this section readily translate
into the classical affine setting with the aid of Proposition 2.2.

\head 3. Skeleta and irreducible components of $\grassSS$ in affine
coordinates
\endhead

In a final transition, we pass to still smaller parametrizing varieties,
which retain {\it stability under the action of the unipotent radical of
$\autlap$\/} and are endowed with particularly convenient affine
coordinates.  For that purpose, we will reduce the problem of exploring
generic properties of (the modules represented by) a fixed irreducible
component of $\grassSS$ --  or $\grasstd$, or
$\modlad$  --  to finding those of varieties parametrizing the modules
with fixed ``skeleton".  The latter varieties are readily accessible,
combinatorially, from quiver and relations of $\la$ (see the remarks
following Theorem 3.8).  As they constitute an {\it open\/} cover of
$\grassSS$, they will give us a good handle on the generic modules
representing the components of $\grassSS$.

Let $\SS$ be a $d$-dimensional semisimple sequence with top $\SS_0 = T$. 
Again, $P$ denotes a projective cover of $T$ with distinguished sequence of
top elements, $z_1, \dots, z_\t$, such that $z_r = e(r)z_r$ (cf\. beginning
of Section 2), our non-canonical input for specifying coordinates.  As
announced, the affine subvarieties to be described will turn out to be 
stable under the action of the unipotent radical
$\unirad$ of $\autlap$.  We can do a little better in fact, as we will
see.  Recall that
$\autlap$ is isomorphic to a semidirect product $\unirad \rtimes
\aut_\la(T)$, and let $\T$ be the following incarnation of a maximal torus
in $\autlap$: namely $\T \cong (K^*)^\t$, where $(a_1,
\dots, a_\t)$ represents the automorphism of $P$ determined by $z_i
\mapsto a_i z_i$.  If   
$G$ is the subgroup $\unirad \rtimes \T$ of $\autlap$, the patches of our
affine cover will in fact be $G$-stable.  If $T$ is squarefree (that is, if
$T$ does not contain any simple summands with multiplicity
$\ge 2$),
$G$ equals $\autlap$, whence our affine cover is $\autlap$-stable in that
case. 
 
Here we make no assumptions on $T$, and 
$\grasstd$ will in general not possess any finite cover consisting of
$\autlap$-stable affine charts, since the $\autlap$-orbits need not be
quasi-affine in general. (In fact, examples attest to the fact that
projective
$\autlap$-orbits of positive dimension may arise if $T$ has repeated simple
summands.)  

\subhead 3.A.  Skeleta \endsubhead

As before, $\la = KQ/I$, and $T \in \lamod$ is semisimple.  Roughly
speaking, skeleta allow us to carry over some of the benefits of the
path-length grading of projective $KQ$-modules to arbitrary
$\la$-modules.   Next to
$P$  --  the distinguished $\la$-projective cover of $T$ with top elements
$z_1, \dots, z_\t$  --   we therefore consider the projective $KQ$-module
$\Phat = \bigoplus_{1 \le r\le \t} (KQ) \zhat_r$, with corresponding
sequence of top elements $\zhat_1, \dots, \zhat_\t$ so that the class of
$\zhat_r$ modulo $I \Phat$ coincides with $z_r$; in particular, $\zhat_r =
e(r) \zhat_r$.  By a {\it path in $\Phat$ starting in $\zhat_r$\/} we mean
any element $\phat = p \zhat_r \in \Phat$, where $p$ is a path in $KQ$
starting in the vertex
$e(r)$.  The {\it length\/} of $\phat$ is defined to be that of $p$; ditto
for the {\it endpoint\/} of $\phat$.  If $p_1$ is an initial subpath of
$p$, meaning that $p = p_2 p_1$ for paths $p_1$, $p_2$, we call
$p_1\zhat_r$ an {\it initial subpath of\/} $p \zhat_r$.  So, in particular,
$\zhat_r = e(r)\zhat_r$ is an initial subpath of length $0$ of any path $p
\zhat_r$ in
$\Phat$.  The reason why, a priori, we do not identify $p \zhat_r \in
\Phat$ with
$p z_r \in P = \Phat/I \Phat$ lies in the fact that we require an
unambiguous notion of path length, which is not guaranteed for paths $pz_r$
in $P$.  However, in the sequel, we will often not make a notational
distinction between $p \zhat_r$ and $p z_r$, unless there is need to
emphasize well-definedness of path lengths.  Moreover, we emphasize that
the label
$r$ is a crucial attribute of a path in
$\Phat$:  given a path $p
\in KQ$ starting in a vertex $e(r) = e(s)$, the elements $p \zhat_r = p
\zhat_s$ of $P$ are distinct unless $r = s$. 

\definition{Definition 3.1 and Further Conventions}
  
\noindent {\bf {(1)}} An ({\it abstract\/}) {\it $d$-dimensional skeleton
with top $T$} is a set $\S$ of paths of lengths at most $L$ in
$\Phat$,  which has cardinality $d$, contains the paths $\zhat_1, \dots,
\zhat_\t$ of length zero (i.e., contains the full sequence of distinguished
top elements of $\Phat$), and is closed under initial subpaths; that is,
whenever
$p_2 p_1
\zhat_r \in \S$, then $p_1 \zhat_r \in \S$.  
\smallskip

\noindent Usually, we view a skeleton $\S$ as a forest of $\t$ tree graphs
each consisting of the paths in $\S$ that start in a fixed  top element
$\zhat_r$ of $\Phat$; see Example 3.3 below and the remarks preceding it. 
Further conventions:  
 
\noindent $\bullet$  By $\S_l$ we denote the set of paths of length $l$ in
$\S$.

\noindent $\bullet$  When we pass from the $KQ$-module $\Phat$ to the
$\la$-module $P$ by modding out $I \Phat$, we often identify $\zhat_r$ with
$z_r$ and view $\S$ as the set 
$$\{p z_r \in P \mid  p \zhat_r \in \S,\, r \le \t \} \subseteq P.$$  
\smallskip    

\noindent {\bf {(2)}}  Let $\S$ be an abstract $d$-dimensional skeleton
with top $T$, and $M$ a $\la$-module.  

We call $\S$ a {\it skeleton of\/} $M$ if there exist top elements $m_1,
\dots, m_\t$ of $M$ and a
$KQ$-epimorphism $f: \Phat \rightarrow M$ satisfying $f(\zhat_r) = m_r$ for
all $r$ such that, for each $l \in \{0, \dots, L\}$, the subset
$$\{f(p\zhat_r) \mid p\zhat_r \in \S_l\} = 
\{p m_r
\mid r \le \t, \, p \zhat_r \in \S_l\}$$  of $M$ induces a $K$-basis for
the radical layer $J^lM / J^{l+1} M$ (here we identify the
$\la$-mul\-ti\-pli\-ca\-tion on $M$ with the induced $KQ$-multiplication). 
In this situation, we also say that
$\S$ is a {\it skeleton of $M$ relative to the sequence $m_1, \dots, m_\t$
of top elements\/}, and observe that the union over $l$ of the above sets,
namely
$\{p m_r \mid  r \le \t, \, p \zhat_r \in \S \}$, is a $K$-basis for
$M$.  As a consequence, we recognize skeleta of $M$ as special $K$-bases
which respect the radical layering of $M$ and are closely tied to its
$KQ$-structure.  

If $M = P/C$, we call $\S$ a {\it distinguished\/} skeleton of $M$ provided
that
$\S$ is a skeleton of $M$ relative to the distinguished sequence $z_1 + C,
\dots, z_\t + C$ of top elements. (Note: If $P/C \cong P/C'$, any
distinguished skeleton of $P/C$ is a skeleton of $P/C'$, but not vice versa
in general, as the top elements are shuffled under the
$\Aut_{\la}(T)$-action; see Example 3.3 below.)
\smallskip    

\noindent {\bf {(3)}} Finally, we define 
$$\grassS = \{ C \in \grasstd \mid \S \text{\ is a distinguished skeleton
of\ } P/C \}.$$
\enddefinition

The set of {\it all\/} skeleta of $M$ is clearly an isomorphism invariant
of $M$.  It contains at least one distinguished candidate when $M = P/C$. 
Our definition of a skeleton coincides in essence with that given in
\cite{\grassI} for the situation of a squarefree top
$T = P/JP$.  However, in that special case, it is unnecessary to hook up
the elements of an abstract skeleton $\S$ with specific top elements of the
$KQ$-module $\Phat$, since the dependence on specific sequences of top
elements disappears.  In particular, every skeleton of $P/C$ is
distinguished in the case of squarefree
$T$.       

We return to the general case.  Again suppose that $M = P/C$, and let
$\S$ be an abstract skeleton.  Then $M \cong P/D$ for {\it some\/} point
$D$ in $\grassS$ if and only if $\S$ is a skeleton of $M$.  Yet, $\autlap.C
\cap \grassS \ne \varnothing$ need not imply $C \in \grassS$.  On the other
hand, the
$\grassS$ do enjoy the following partial stability under the
$\autlap$-action:

\proclaim{Observation 3.2} The sets $\grassS$ are locally closed
subvarieties of $\grasstd$, which are stable under the action of $G =
\unirad \rtimes \T$.  Moreover, $\grasstd$ is the union of the $\grassS$. 
\endproclaim

\demo{Proof}  Local closedness will follow from Observation 3.5 below,
since the $\grassSS$ are known to be locally closed in $\grasstd$.

As for stability: Suppose $C \in \grassS$, meaning that $\S$ is a
distinguished skeleton of $P/C$.  Then $\S$ remains a distinguished
skeleton of
$P/C$ after passage from our given distinguished sequence $z_1, \dots,
z_\t$ of top elements of $P$ to a new distinguished sequence of the form
$g.z_1, \dots, g.z_\t$ for any
$g \in G$.  But this is tantamount to saying that $\S$ is a distinguished
skeleton of $P/(g^{-1}.C)$ relative to the original sequence.  In other
words, $g^{-1}.C \in \grassS$.

The straightforward fact that each $\la$-module
$P/C$ with $C \in \grasstd$ has at least one distinguished skeleton yields
the remaining assertion.
\qed \enddemo  

Any abstract skeleton $\S$ can be communicated by means of an undirected
graph which is a {\it forest\/}, that is, a finite disjoint union of tree
graphs.  There are $\t$ trees if $\S$ is a skeleton with top $T$, one for
each $r \le \t$;  here $\zhat_r$, identified with $e(r)$, represents the
root of the $r$-th tree, recorded in the top row of the graph. The paths $p
\zhat_r$ of positive length in $\S$ are represented by edge paths of
positive length.  Instead of formalizing this convention, we will
illustrate it with an example to which we will return in Sections 3.C and
4.  

\definition{Example 3.3}  Let $\la = KQ$, where $Q$ is the quiver

\ignore
$$\xymatrixrowsep{1pc}\xymatrixcolsep{3pc}
\xymatrix{  &&&2 \ar[dl]_{\gamma} \\ 1 \ar[r]^{\alpha} &4
\ar[r]<0.6ex>^{\beta_1}
\ar[r]<-0.6ex>_{\beta_2} &6 &5 \ar[l]^{\epsilon} &3 \ar[l]^{\delta} }$$
\endignore

\noindent Moreover, let $T = S_1^2 \oplus S_2 \oplus S_3$ and $P =
\bigoplus_{1
\le r \le 4} \la z_r$ with distinguished top elements $z_1 = (e_1,0,0,0)$,
$z_2 = (0,e_1,0,0)$, $z_3 = (0,0,e_2,0)$, $z_4 = (0,0,0,e_3)$.  Choose $d =
9$.  We consider the module $M = P/C$, where $C$ is the submodule of $P$
generated by
$\beta_2 \alpha z_1$, $\beta_1 \alpha z_2$, $\gamma z_3 - \epsilon
\delta z_4$, and $\beta_1 \alpha z_1 + \beta_2 \alpha z_2 + \gamma z_3$. 
Then $M$ has precisely three distinguished skeleta as follows:
\roster
\item $\quad \{e_1 \zhat_1, \alpha \zhat_1, \beta_1 \alpha \zhat_1 \} 
\sqcup \{e_1 \zhat_2, \alpha \zhat_2, \beta_2 \alpha \zhat_2 \}
\sqcup \{e_2 \zhat_3\} \sqcup \{e_3 \zhat_4, \delta \zhat_4\}$;
\item $\quad \{e_1 \zhat_1, \alpha \zhat_1, \beta_1 \alpha \zhat_1\} 
\sqcup \{e_1 \zhat_2, \alpha \zhat_2 \} \sqcup \{e_2 \zhat_3\} \sqcup
\{e_3 \zhat_4, \delta \zhat_4, \epsilon \delta \zhat_4\}$; 
\item $\quad \{e_1 \zhat_1, \alpha \zhat_1\} 
\sqcup \{e_1 \zhat_2, \alpha \zhat_2, \beta_2 \alpha \zhat_2\} \sqcup
\{e_2 \zhat_3\} \sqcup
\{e_3 \zhat_4, \delta \zhat_4, \epsilon \delta \zhat_4\}$.
\endroster

\noindent There are five additional, non-distinguished, skeleta of $M$. 
The first two of the three graphs below show skeleta of $M$ relative to the
permuted sequence of top elements $(z_2 + C, z_1 + C, z_3 + C, z_4 + C)$. 
The rightmost graph displays a skeleton of $M$ relative to the sequence
$(z_1 + z_2 + C, z_2 + C, z_3 + C, z_4 + C)$ of top elements.

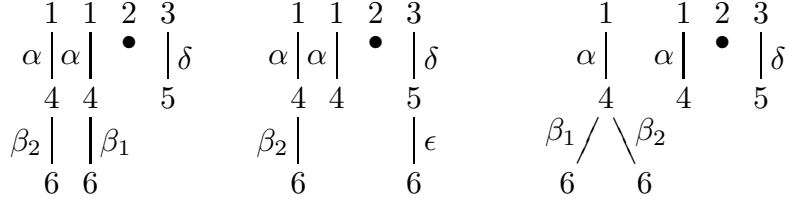
\ignore
$$\xymatrixrowsep{1.2pc}\xymatrixcolsep{0.1pc}
\xymatrix{   1 \edge[d]_{\alpha} &1 \edge[d]_{\alpha} &2 \drbl &3
\edge[d]^{\delta} &&&&&1 \edge[d]_{\alpha} &1 \edge[d]_{\alpha} &2 \drbl &3
\edge[d]^{\delta} &&&&&&&1 \edge[d]_{\alpha} &&1 \edge[d]_{\alpha} &2
\drbl &3 \edge[d]^{\delta} \\  4 \edge[d]_{\beta_2} &4 \edge[d]^{\beta_1}
&&5 &&&&&4 \edge[d]_{\beta_2} &4 &&5 \edge[d]^{\epsilon} &&&&&&&4
\edge[dl]_{\beta_1} \edge[dr]^{\beta_2} &&4 &&5 \\   6 &6 && &&&&&6 &&&6
&&&&&&6 &&6 }$$
\endignore

Note that no skeleton of $M$ contains the path $\gamma \zhat_3$, since
$\gamma M \subseteq J^2 M$.
\qed

\enddefinition

\subhead 3.B.  $\grassSS$ as a union of $\grassS$'s \endsubhead
\medskip

For unproven statements in this subsection and the next, we refer to
\cite{\grassIII}.

Suppose that $\SS$ is a
$d$-dimensional semisimple sequence with top $T$.  Then $\grassSS$ is the
union of those subvarieties $\grassS$ which have non-empty intersection
with $\grassSS$; indeed, $\grassSS \cap \grassS \ne
\varnothing$ implies  that $\grassS$ is contained in
$\grassSS$.  These are precisely those nonempty candidates among the
$\grassS$ which are based on skeleta $\S$ compatible with
$\SS$ in the following sense:  

\definition{Definition 3.4} Given a semisimple sequence $\SS$, we call a
skeleton $\S$   {\it compatible with $\SS$ = $=$ $(\SS_0,
\dots, \SS_L)$\/} if, for each
$l \le L$ and $i \le n$, the number of paths in  $\S_l$ ending in the
vertex $e_i$ coincides with the multiplicity of the simple module $S_i$ in
$\SS_l$. 
\enddefinition

Evidently, each abstract skeleton compatible with $\SS$ shares the 
dimension and top with $\SS$. In fact, given any skeleton $\S$ of a
$\la$-module
$M$, the radical layering $\SS(M)$ of $M$ is the only semisimple sequence
with which $\S$ is compatible.  

Suppose that $\S$ is compatible with
$\SS$. Provided that $\grassS \ne \varnothing$, the sum of the subspaces
$Kp z_r \subseteq JP$, with $p \zhat_r \in \S$ of positive length, is
direct  --  this follows from the definition of a skeleton  --  and the
variety $\grassS$ is the intersection of $\grassSS$ with the following big
Schubert cell $\SchuS$ of $\grasstd$: Namely,
$$\SchuS = \{C \in \grasstd \mid  JP = C \oplus \bigoplus_{p \zhat_r \in
\S,\,
\len(p) > 0} Kp z_r \}.$$ Thus we obtain:

\proclaim{Observation 3.5 and Terminology}  The $\grassS$, where $\S$ runs
through the skeleta compatible with $\SS$, form an open cover of
$\grassSS$; in general, the $\grassS$ fail to be open in $\grasstd$,
however.

Suppose $\C$ is an irreducible component of $\grassSS$ and
$\frak S$ the {\rm(}finite{\rm)} set of skeleta $\sigma$ with the property
that
$\grassS \cap \C \ne \varnothing$.  Then $\C \cap \bigcap_{\S \in {\frak
S}} \grassS$ is a dense open subset of $\C$.   We call $\frak S$ the {\rm
generic set of skeleta\/} of the modules in
$\C$.
\qed
\endproclaim

Thus, finding the irreducible components of $\grassSS$ can be played back
to finding those of the $\grassS$ (as mentioned, this task is
computationally mastered).  Indeed, we have the following correspondence;
its elementary proof is left to the reader. 

\proclaim{Observation 3.6}  Let $\SS$ be a semisimple sequence such that
$\grassSS \ne \varnothing$.           

{\rm{(1)}}  Whenever
$\C$ is an irreducible component of $\grassSS$ and 
$\S$ a skeleton such that $\grassS$ intersects $\C$ nontrivially, there
exists a unique irreducible component $\D$ of
$\grassS$ whose closure in $\grassSS$ equals $\C$.  In fact, $\D = \C
\cap \grassS$.

{\rm{(2)}}  Conversely, given any irreducible component $\D$ of a nonempty
$\grassS$, where $\S$ is a skeleton compatible with $\SS$, the closure of
$\D$ in $\grassSS$ is an irreducible component of $\grassSS$.
\qed \endproclaim

\subhead 3.C.  Affine coordinates and irreducible components of the
$\grassS$ \endsubhead
\medskip

As we saw in the preceding subsection, the map assigning to each point
$x$ in an irreducible component $\C$ of $\modlad$ the set of skeleta of the
module corresponding to $x$ is constant on a dense open subset of
$\C$.  This generic behavior of skeleta singles them out as relevant
for the construction of generic modules representing the irreducible
components of $\modlad$.

Let $\S$ be a $d$-dimensional skeleton with top $T$.  Recall from
3.A that this makes $\S$ a set of ``paths" with well-defined lengths in the
projective $KQ$-module
$\Phat = \bigoplus_{1 \le r \le \t} KQ \zhat_r$.  We supplement Definition
3.1 as follows:

\definition{Definition 3.7} A {\it $\S$-critical path\/} is a path
$\alpha p \zhat_r$ of length at most $L$ in $\Phat \setminus \S$, where
$\alpha$ is an arrow and $p \zhat_r \in \S$.  Moreover, for any such
$\S$-critical path, the $\S$-{\it set of $\alpha p \zhat_r$\/}  is the set
$\S(\alpha p \zhat_r)$ consisting of all paths $q \zhat_s
\in
\S$  which are at least as long as $\alpha p \zhat_r$ and end in the same
vertex as $\alpha p \zhat_r$.  
\enddefinition

In other words, a path in $\Phat$ is $\S$-critical if and only if it fails
to belong to $\S$, whereas every proper initial subpath does. If
$\S$ is Skeleton (1) in Example 3.3, the
$\S$-critical paths are $\beta_2 \alpha \zhat_1$, $\beta_1 \alpha
\zhat_2$, $\gamma \zhat_3$ and $\epsilon \delta \zhat_4$, with
$\S(\beta_2 \alpha \zhat_1)$ $=$ $\{\beta_1 \alpha \zhat_1, \beta_2
\alpha
\zhat_2\}$ $=$ $\S(\beta_1 \alpha \zhat_2)$ and 
$\S(\gamma \zhat_3)$ $=$ $\{\beta_1 \alpha \zhat_1, \beta_2 \alpha
\zhat_2\}$
$=$ $\S(\epsilon \delta \zhat_4)$.  If we add the arrows $\alpha$ giving
rise to $\S$-critical paths $\alpha p \zhat_r$ to the graph of this
skeleton, marking them by broken edges, we obtain the following picture:

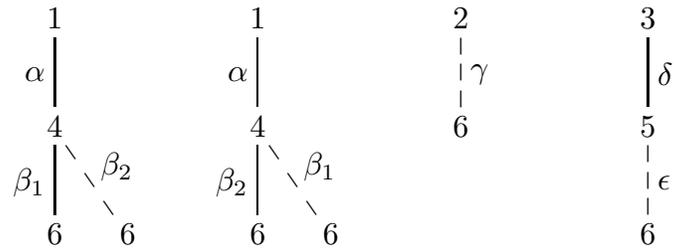
\ignore
$$\xymatrixrowsep{1.8pc}\xymatrixcolsep{1pc}
\xymatrix{  1 \edge[d]_{\alpha} &&&1 \edge[d]_{\alpha} &&&2
\dashedge[d]^{\gamma} &&&3
\edge[d]^{\delta}  \\ 4 \edge[d]_{\beta_1} \dashedge[dr]^{\beta_2} &&&4
\edge[d]_{\beta_2}
\dashedge[dr]^{\beta_1} &&&6 &&&5 \dashedge[d]^{\epsilon}  \\ 6 &6 &&6 &6
&&&&&6 }$$
\endignore

That $C$ be a point in $\grassS$ obviously entails the existence of unique
scalars
$c(\alpha p\zhat_r,q\zhat_s)$ with the property that
$$\alpha p z_r + C =  \sum_{q\zhat_s \in \S(\alpha p\zhat_r)} c(\alpha
p\zhat_r,q\zhat_s)\, q z_s + C,$$  in $P/C$, whenever $\alpha p \zhat_r$ is
a $\S$-critical path.  Clearly, the isomorphism type of $M = P/C$ is
completely determined by the family of these scalars.  Thus we obtain a map 
$$\grassS \rightarrow \AA^N,\ \ \ \ C \mapsto c = \bigl(c(\alpha
p\zhat_r,q\zhat_s)\bigr)_{\alpha p\zhat_r\ \S\text{-critical},\ q\zhat_s \in
\S(\alpha p\zhat_r)}$$     where $N$ is the (disjoint) union of the sets
$\{\alpha p\zhat_r\}\times \S(\alpha p\zhat_r)$.   Note: Since, a priori,
the sets
$\S(\alpha p\zhat_r)$, where $\alpha p\zhat_r$ traces the $\S$-critical
paths, need not be disjoint, we have paired their elements with the
pertinent
$\sigma$-critical paths to force disjointness.    

By \cite{\grassIII}, the map $\grassS \rightarrow \AA^N$,
$C \mapsto c$, described above is an isomorphism from $\grassS$ onto a
closed subvariety of $\AA^N$.  The point $C \in \grassS$ corresponding to a
point $c$ in the image of this map is the submodule of $JP$ generated over
$\la$ by the elements $\alpha p z_r -  \sum_{q \zhat_s \in
\S(\alpha p \zhat_r)} c(\alpha p\zhat_r,q\zhat_s)\, q z_s$, where $\alpha
p\zhat_r$ runs through the
$\S$-critical paths.  We will identify $C$ with $c$ whenever convenient. 

The following result from \cite{\grassIII} summarizes those properties of
the cover
$\bigl(\grassS \bigr)_\S$ of $\grasstd$ which will be relevant here; the
first statement partially overlaps with Sections 3.A/B.  

\proclaim{Theorem 3.8}  For every $d$-dimensional skeleton $\S$ with top
$T$, the set $\grassS$ is a locally closed affine subvariety of
$\grasstd$ which is stable under the action of the group $G = \unirad
\rtimes
\T$.  Moreover, given any
$d$-dimensional semisimple sequence $\SS$ with top
$T$, the varieties $\grassS$, where $\S$ traces the skeleta compatible with
$\SS$, form an open affine cover of $\grassSS$. 

Polynomial equations determining the $\grassS$ in $\AA^N$ can be
algorithmically obtained from $Q$ and generators for the admissible ideal
$I \subseteq KQ$, where $I$ is viewed as a left ideal.  If $\Kcirc$ is a
subfield of
$K$ over which such generators for $I$ are defined, the resulting
polynomials determining $\grassS$ are defined over $\Kcirc$ as well.
\qed \endproclaim

The authors have implemented the mentioned algorithm for $\grassS$ at
\cite{\codes}, and combined it with a software package that computes the
irreducible components.

The usefulness of the graphing technique of \cite{\domino} and
\cite{\menace} is limited, since it calls for display of {\it all\/}
relations of a module $P/C$, including those that are consequences of
others.  The following more sparing graphs will better serve our purpose of
graphically representing ``generic modules" for the irreducible components
of various parametrizing varieties.  Recall that the term {\it
hypergraph\/} is commonly used to refer to an undirected graph which not
only allows for conventional edges coupling pairs of vertices, but for {\it
hyperedges\/} connecting more than two vertices.  The hypergraphs of a
module $P/C$ considered here presuppose a fixed choice of top elements of
$P$.  Roughly speaking, they consist of a distinguished skeleton of $P/C$
(which may be viewed as a forest of traditional tree graphs  --  see
Subsection 3.A), combined with suitable subsets of the $\S$-sets of the
$\S$-critical paths which in turn depend on $C$;  these latter subsets
provide additional (hyper)edges.   

\definition{Definition 3.9. Hypergraphs of Modules} Suppose $M$ is a
$d$-dimensional module with top $T$ and skeleton $\S$.  This means that
$M \cong P/C$ for some $C \in \grassS$.  Denote by 
$$\bigl(c(\alpha p\zhat_r,q\zhat_s) \bigr)_{\alpha p\zhat_r\ 
\S\text{-critical},\, q\zhat_s \in
\S(\alpha p\zhat_r)}$$  the local affine coordinates of $C$ relative to
$\S$.  Then the following will be called a {\it hypergraph of $M$\/}:  The
skeleton $\S$, paired with the family $\bigl( M(\alpha p\zhat_r)
\bigr)$ of subsets $M(\alpha p\zhat_r) \subseteq \S(\alpha p\zhat_r)$,
respectively, where $\alpha p\zhat_r$ ranges through the $\S$-critical paths
and
$M(\alpha p\zhat_r)$ consists of those paths $q\zhat_s \in \S(\alpha
p\zhat_r)$ for which
$c(\alpha p\zhat_r,q\zhat_s) \ne 0$; here it is important that the sets
$M(\alpha p\zhat_r)$ be recorded in conjunction with the $\S$-critical paths
by which they are labeled.
\enddefinition 

In depicting such a hypergraph, we start with the graph of $\S$, add on
edges representing the terminating arrows $\alpha$ of the
$\S$-critical paths $\alpha p\zhat_r$  --  the latter distinguished from
the edges of $\S$ via broken lines  --  and finally mark the {\it
hyperedges\/}, one for each
$M(\alpha p\zhat_r) \subseteq \S(\alpha p\zhat_r)$, by means of closed
curves including the sets of endpoints of the paths in $\{\alpha p\zhat_r\}
\cup M(\alpha p\zhat_r)$.  In case $M(\alpha p\zhat_r)$ is a singleton, say
$\{q\zhat_s\}$, we simply join the endpoints of $\alpha p\zhat_r$ and
$q\zhat_s$.  If, for the algebra
$\la$ of Example 3.3, we take $T = \S_1^2$, $P = (\la e_1)^2$  --  meaning
that
$z_1$, $z_2$ are both normed by $e_1$  --  and $C = \la(\beta_1 \alpha -
\beta_2
\alpha) z_1 +
\la[\beta_1 \alpha z_1 + (\beta_1 \alpha - \beta_2 \alpha) z_2]$, then the
hypergraph of $P/C$ with respect to the skeleton $\S = \{\zhat_1 = e_1
\zhat_1,\, \alpha \zhat_1,\, \beta_1
\alpha \zhat_1,\, \zhat_2 = e_1 \zhat_2,\, \alpha \zhat_2,\, \beta_1
\alpha \zhat_2\}$ consists of $\S$, together with the sets $M(\beta_2
\alpha \zhat_1) = \{\beta_1 \alpha \zhat_1\}$ and
$M(\beta_2 \alpha \zhat_2) = \{ \beta_1 \alpha \zhat_1,\, \beta_1
\alpha \zhat_2\}$.  It looks as follows:

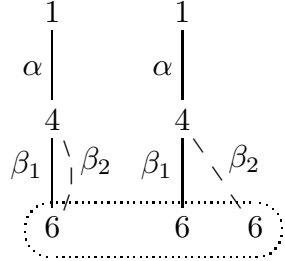
\ignore
$$\xymatrixrowsep{1.8pc}\xymatrixcolsep{1pc}
\xymatrix{  1 \edge[d]_{\alpha} &&1 \edge[d]_{\alpha}  \\ 4
\edge[d]_(0.4){\beta_1} \dashedge@/^/[d]^(0.4){\beta_2} &&4
\edge[d]_(0.4){\beta_1} \dashedge[dr]^{\beta_2}  \\ 6 \levelpool3 &&6 &6 }$$
\endignore

\noindent Note that this hypergraph of $P/C$ is connected. More strongly,
$P/C$ is indecomposable.
 
For further illustrations, see Examples 4.7, 5.8, and 5.10.

The final observation of this section is an immediate consequence of the
definitions.  We leave the easy proof to the reader. 

\proclaim{Observation 3.10 and Notational Convention}  Let $\SS$ be a
$d$-dimensional semisimple sequence with top $T$ and $\S$ a skeleton
compatible with $\SS$. Then the number
$$N(\SS) = \sum_{\alpha p\zhat_r \ \sigma\text{-critical}} |\sigma(\alpha
p\zhat_r)|$$ depends only on $\SS$, not on $\S$.

Moreover, $N(\SS) = N_0(\SS) + N_1(\SS)$, where both of the numbers
$$N_0(\SS) =  \sum_{\alpha p\zhat_r \ \sigma\text{-critical}} |\{ q\zhat_s
\in \sigma(\alpha p\zhat_r) \mid \len(q) = \len(\alpha p)\}|$$  and
$$N_1(\SS) = \sum_{\alpha p\zhat_r \ \sigma\text{-critical}} |\{ q\zhat_s
\in \sigma(\alpha p\zhat_r) \mid \len(q) > \len(\alpha p)\}|$$  depend
solely on $\SS$. \qed
\endproclaim

\head 4. $\SS$-generic, $(T,d)$-generic, and $d$-generic modules \endhead

{\it Throughout this section, we assume the base field $K$ to have infinite
transcendence degree over its prime field. Moreover, we fix an
algebraically closed subfield
$\Kcirc$ of $K$ such that $K$ has infinite transcendence degree over
$\Kcirc$ and $I \subseteq KQ$ is defined over $\Kcirc$.\/}  The latter means
that, as a $KQ$-ideal,
$I$ is generated by $I_0 = I\cap \Kcirc Q$.  

The assumption on the transcendence degree of our base field will enable us
to realize and study ``$\SS$-generic", ``$(T,d)$-generic", and
``$d$-generic" modules inside the category
$\lamod$.  Since we will obtain one generic object for every irreducible
component of $\grassSS$ (or $\grasstd$, or $\modlad$), unique up to
suitable equivalence, there will be only finitely many generic modules to
be studied in each dimension, which permits assembling essential
information on the representation theory of $\la$ in a finite frame.  (Here
the attribute ``generic" is unrelated to the generic modules introduced by
Crawley-Boevey in \cite{\CB}.  However, it is unlikely to  conflict with
existing usage, as we attach specifying prefixes;  moreover,
Crawley-Boevey's generic modules are infinite dimensional by definition,
whereas we exclusively consider finite dimensional modules.)

In the present setting, $\la$ is clearly isomorphic to the algebra  $\la_0
\otimes_{\Kcirc} K$ via $\lambda_0 \otimes k \mapsto  k \lambda_0$, where
$\la_0 = \Kcirc Q / I_0$, and all of the varieties $\grasstd$,
$\grassSS$ and
$\grassS$ are defined over $\Kcirc$.  For the $\grassS$, this follows from
Theorem 3.8; both $\grasstd$ and $\grassSS$ are finite unions of
$\grassS$'s.  We denote by $\bigl(\grasstd\bigr)_0$,
$\bigl(\grassSS\bigr)_0$ and $\bigl(\grassS\bigr)_0$ the restrictions of
the mentioned varieties to $\Kcirc$.  Hence, given an irreducible component
$\C$ of $\grasstd$, $\grassSS$, or $\grassS$, it makes sense to refer to
the corresponding irreducible component $\C_0$ of $\bigl(\grasstd\bigr)_0$,
$\bigl(\grassSS\bigr)_0$, or
$\bigl(\grassS\bigr)_0$.    

Next to the described algebraic group actions on the parametrizing
varieties, we will consider the following action of the Galois group
$\Gal$.  It results from the obvious fact
that any projective $K$-space carries such an action and the closed
subvariety
$\grasstd$ of $\PP(\bigwedge^{\tilde{d}} JP)$ is defined over $\Kcirc$ (again,
$\tilde{d} =\dim P - d$).  We briefly discuss the forms this $\Gal$-action
takes on in terms of the coordinate systems we are using.  Given any
skeleton $\S$ with $\grassS \ne \varnothing$, we view $\S$ as a (necessarily
$K$-linearly independent) subset of
$P$, and supplement $\S \cap JP$ to a path basis for $JP$.  This provides a
homogeneous coordinatization of
$\PP(\bigwedge^{\tilde{d}} JP)$, and we let a map $\phi \in \Gal$ act on a
point in projective space by applying $\phi$ to the homogeneous
coordinates.  The action does not depend on the choice of path basis (and
hence does not depend on $\S$), since the transition matrices from one such
basis to another have coefficients in $\Kcirc$.  Moreover, we note: The
subvarieties $\grasstd$,
$\grassSS$, and $\grassS$ are stable under this $\Gal$-action, and all maps
$\phi \in \Gal$ leave the restrictions
$(\grasstd)_0$, etc., pointwise fixed.  In particular, the $\Gal$-action
stabilizes the irreducible components of the considered varieties because
they are all defined over $\Kcirc$.  Given
$C \in \grasstd$, we will write $C_\phi$ for $\phi.C$ to set it off from
the $\autlap$-action.   For a point $C \in \grassS$, the affine coordinates
$(c_\nu)_{\nu \in N}$ of Section 3.C are readily seen to be $\Kcirc$-linear
combinations of certain Pl\"ucker coordinates of $\PP(\bigwedge^{\tilde{d}}
JP)$ relative to the described path basis for $JP$; in fact, it is routine
to check that $C$ is determined by a subcollection of the $c_\nu$ which
coincide with selected Pl\"ucker coordinates (that we do not need to keep
track of all of them in order to pin down $C$, is due to the fact that $C$
is not only a
$K$-space, but a $\la$-submodule of $JP$).  Thus $C_\phi =
\bigl(\phi(c_\nu)\bigr)$ in the distinguished affine coordinates of
$\grassS$.  Clearly, $\Gal$ acts on the classical affine module varieties
$\modlad$, $\toptd$ and $\laySS$ as well, and these actions are compatible
with the ones on their Grassmannian counterparts under the transfer maps for
orbits described in Proposition 2.2.   

Finally, we relate the above actions on the varieties $\grasstd$ to the
well-known fact that $\Gal$ also acts on the category
$\Lamod$, in that each $\phi \in \Gal$ gives rise to a ring automorphism of
$\la$ (actually, a $\Kcirc$-algebra automorphism), also denoted $\phi$, which
in turn induces a self-equivalence $F_\phi$ of $\lamod$. 

\proclaim{Observation 4.1 and Definition of $\Gal$-equivalence and
stability}  Let $\phi \in
\Gal$ and, for $M \in \Lamod$, denote by $F_\phi(M)$ the
$\la$-module with underlying abelian group $M$ and module multiplication
$\lambda *m = \phi^{-1}(\lambda) m$.  Then the assignment
$M\mapsto F_\phi(M)$ extends to an equivalence 
$$F_{\phi}: \Lamod \longrightarrow \Lamod$$ of categories.  For $C \in
\grasstd$, this equivalence sends $P/C$ to a module in the isomorphism
class of $P/C_{\phi}$, where $C_\phi$ is as above; in particular, $F_\phi$
sends $P = \bigoplus_{1 \le r \le \t} \la z_r$ to an isomorphic copy of
itself.  Moreover, $F_\phi(C) \cong C_\phi$ in
$\Lamod$, under the $\la$-isomorphism $P
\rightarrow P$ which sends $k z_r$ to $\phi(k) z_r$ for $k \in K$ and all
$r$.  

Clearly, $F_{\phi}$ preserves all properties of
$\la$-modules which are invariant under Morita-equivalence, fixes the
isomorphism classes of the indecomposable projectives and the simple
modules, and preserves dimension vectors of modules and their radical
layers; a fortiori,
$F_{\phi}$ preserves $K$-dimension.  Moreover, for $C, D \in \grasstd$, a
homomorphism $f: P/C
\rightarrow P/D$ is an isomorphism precisely when this is true for
$f_\phi = F_{\phi} (f): P/C_{\phi} \rightarrow P/D_\phi$.  Consequently,
the action of $\Gal$ permutes the $\autlap$-orbits within each irreducible
component of $\grasstd$, $\grassSS$, or $\grassS$. 
\smallskip

{\rm We call two $\la$-modules $M$ and $M'$}  equivalent {\rm \,or, more
precisely,} $\Gal$-equivalent{\rm , if there exists a map $\phi \in \Gal$
such that $M' \cong F_\phi(M)$.  Moreover, two connected components $\A$
and $\A'$ of the Auslander-Reiten quiver of $\lamod$ are said to be}
$\Gal$-equivalent {\rm \, if there exists $\phi \in \Gal$ with the property
that the functor
$F_\phi$ takes $\A$ to $\A'$.  

A property of modules is said to be} $\Gal$-stable {\rm if it is  passed on
from a module to any equivalent module.  Analogously, we define stability
for properties of} pairs {\rm of modules, calling such a  property}
$\Gal$-stable {\rm if, for all $\phi$, either both pairs  
$(M, M')$ and $(F_\phi(M), F_\phi(M'))$ have this property, or else neither
of them does.}
\qed \endproclaim

Note that the functors $F_\phi$ of Observation 4.1 are $\Kcirc$-linear
selfequivalences of $\Lamod$, but fail to be $K$-linear in all nontrivial
cases.  Next to categorically defined features (such as the standard
homological properties), the $\Gal$-stable properties of a module include
its set of skeleta, since $F_\phi(P/C) = P/C_\phi$ whenever $C \in
\grassS$.  The concept of $\Gal$-stability obviously translates into our
parametrizing varieties as follows:  A property (*) that pertains to
$d$-dimensional modules with top $T$, for instance, is
$\Gal$-stable precisely when it is preserved under module isomorphism and
the set of points $C$ in $\grasstd$ such that $P/C$ has property (*) is
stable under the $\Gal$-action on $\grasstd$.  

To illustrate the concept of $\Gal$-equivalence, we return to Example 3.3,
taking the base field to be $\CC$ and the subfield $\Kcirc$ to be the field
$\overline{\QQ}$ of algebraic numbers.  Then the following $\la$-modules
$M_1$, $M_2$ are $\Gal$-equivalent. They have projective cover $P = (\la
e_1)^2$, endowed with distinguished top elements $z_1, z_2$, both normed by
$e_1$.  We define
$M_j$ to be $P/C_j$ for $j = 1,2$, where
$$C_1 = \la (\beta_1 \alpha - \pi^3 \beta_2 \alpha) z_1\  +\
\la[\beta_1 \alpha z_1 +  (\pi \beta_1 \alpha - {\root 4 \of \pi}
\beta_2 \alpha) z_2 ],$$  and
$$C_2 = \la(\beta_1 \alpha - (1/ \pi)^3 \beta_2 \alpha) z_1\  +\
\la[\beta_1 \alpha z_1 +  \bigl( (1 / \pi) \beta_1 \alpha + i ( 1/{\root 4
\of \pi}) \beta_2 \alpha \bigr) z_2 ].$$ On the other hand, $M_1$ and $M_2$
fail to be isomorphic.

As pointed out, $\Gal$-equivalent modules have the same sets of skeleta;
clearly, they also have coinciding hypergraphs (see 3.9 for our
conventions).  In our example, the hypergraph of the equivalent modules
$M_1$ and $M_2$ with respect to the skeleton
$\S = \{z_1, \alpha z_1, \beta_1 \alpha z_1, z_2, \alpha z_2, \beta_1
\alpha z_2\}$ is the same as that of the module $P/C$ described in 3.9.

\definition{Definition 4.2 of $\SS$-, $(T,d)$-, and $d$-generic modules}  
Suppose that
$\C$ is an irreducible component of $\grassSS$.  A module $G \in
\lamod$ is called {\it $\SS$-generic for $\C$\/} (relative to $\Kcirc$) if 

$\bullet$  $G \cong P/C$ for some $C \in
\C$, and

$\bullet$ $G$ has all $\Gal$-stable generic properties of the modules in
$\C$.  By this we mean that $\autlap.C$ has nonempty intersection with
every $\Gal$-stable open subset of $\C$.
\smallskip 

If the closure of $\C$ in
$\grasstd$ is an irreducible component of the latter variety, then an
$\SS$-generic module $G$ for $\C$ is also called {\it
$(T,d)$-generic\/}, or, more precisely, {\it $(T,d)$-generic for the
closure of $\C$ in 
$\grasstd$\/}.  Finally, let $\C'$ be the irreducible component of
$\laySS$ which corresponds to $\C$ via the bijection of Proposition 2.2. 
If the closure of
$\C'$ in $\modlad$ is an irreducible component of $\modlad$, the module $G$
is also referred to as {\it $d$-generic\/}, or, by an abuse of language,
{\it $d$-generic for the closure of $\C$ in
$\modlad$\/}. 
\enddefinition

Note that the second requirement we imposed on a generic module for $\C$ has
the following equivalent formulation: $\autlap.C$ has nonempty intersection
with the $\Gal$-closure of any dense open subset of $\C$.

For several explicit samples of $\SS$-, $(T,d)$-, and $d$-generic modules,
see Example 4.8 and Section 5.  Clearly, the property of being generic on
one of the indicated levels is passed on from  a module $G$ to any module
isomorphic to $F_{\phi}(G)  = G_{\phi}$ for some $\phi \in \Gal$; in other
words, being $\SS$-generic (or 
$(T,d)$-generic, or $d$-generic) is $\Gal$-stable.  In
particular, this is a Zariski-dense property in light of the next theorem;
but it may fail to be an open property, as is easily seen for
$2$-dimensional modules over the Kronecker algebra.  With the next theorem,
we establish existence of generic modules, as well as their uniqueness up to
$\Gal$-equivalence, on all of the considered levels: $\SS$, $(T,d)$ and $d$.
A proof will be given after Lemma 4.4.

\proclaim{Theorem 4.3}  Let $\SS = (\SS_0, \dots, \SS_L)$ be any semisimple
sequence with $\grassSS \ne \varnothing$. 

Given any irreducible component $\C$ of $\grassSS$, there exists an
$\SS$-generic $\la$-module $G = G(\SS, \C)$ for $\C$.  It is unique up to
$\Gal$-equivalence {\rm {(}}i.e., unique up to a shift by some
auto-equivalence $F_\phi$ of $\Lamod${\rm {)}}.

{\bf Supplement {\rm 1:}\,}  Every skeleton $\S$ with $\grassS \cap \C
\ne
\varnothing$ is a skeleton of $G$, and any hypergraph of $G$ {\rm {(see
3.9)}} is shared by all other generic modules for $\C$. More precisely, for
each such skeleton $\S$, the generic module $G$ has a minimal projective
presentation as follows. Let
$P =
\bigoplus_{1 \le r \le \t} \la z_r$ be the distinguished projective cover of
the top $T = \SS_0$ of $G$.  Then $G \cong P/C$, where  
$$C \ = \  \sum_{\alpha p \zhat_r\ \S\text{-critical}} \la g_{\alpha
p,r} \qquad \text{with} \qquad g_{\alpha p,r} = \alpha p z_r - \sum_{q
\zhat_s \in \S(\alpha p \zhat_r)} c(\alpha p\zhat_r,q\zhat_s)\, q
z_s $$  
for a suitable family of scalars
$(c(\alpha p\zhat_r,q\zhat_s))$ in
$\grassS \cap \C$, such that the transcendence degree of the field
$\Kcirc\bigl(c(\alpha p\zhat_r,q\zhat_s) \mid \alpha p\zhat_r \text{\ is\ } \S
\text{-critical},\ q\zhat_s
\in \S(\alpha p\zhat_r)\,\bigr)$ over $\Kcirc$ equals the dimension of $\C$. 
Conversely, every module with a presentation of this type is generic for
$\C$.  
   
{\bf Supplement {\rm 2:}\,}  One can choose generic modules $G(\SS,\C)$ so
as to have the following properties relative to one another:  Given
two semisimple sequences $\SS$ and $\SS'$, together with
irreducible components $\C$ and $\C'$ of $\grassSS$ and $\grassSS'$,
respectively  --  in case $\SS = \SS'$, we require that $\C$ and $\C'$ be
distinct  --  any generic
$\Gal$-stable property of pairs of modules in $\C
\times \C'$ is shared by the pair $(G(\SS,\C), G(\SS',\C'))$. 
\endproclaim

In light of Section 2, every irreducible component $\D$ of $\modlad$ (or
$\toptd$) is the closure in $\modlad$ (resp. in
$\toptd)$) of a unique irreducible component $\C$ of $\laySS$.  In this
situation, the module $G(\SS, \C)$ is even $d$-generic for $\D$, or $(T,
d)$-generic for $\D$, respectively.  Consequently, Theorem 4.3 has two
spinoffs, namely exact analogues of its assertions for $d$- and
$(T,d)$-generic modules. 

The gist of Supplement 1 is this:  If one understands the geometry of a
nonempty intersection $\grassS \cap \C \ne \varnothing$ up to birational
equivalence, one obtains a concrete handle on the module $G$.  This
observation  has immediate applications to algebras for which one has
full grasp of the
$\grassS$, notably truncated path algebras (see Section 5).  In Supplement
2, the requirement that $\C \ne \C'$ is indispensable in general;
Corollary 4.7(d) indicates how to deal with the case $\C = \C'$.
\medskip

The lion's share of the theorem rests on a general fact concerning
``generic orbits" of algebraic group actions, which  allows for multiple
variations in different directions.  For our present
purpose, we let $\Gamma$ be a connected algebraic group over $K$ acting
morphically on an irreducible quasi-projective $K$-variety
$\C$, which is given in terms of an embedding into a projective $K$-space. 
We suppose that $\C$ is defined over
$\Kcirc$, so that the
$\Gal$-action on the encompassing projective space restricts to a
$\Gal$-action on
$\C$. This action is denoted in the form $(\phi,C) \mapsto C_\phi$, while
the $\Gamma$-action is given by $(g,C) \mapsto g.C$.  Moreover, we suppose
that $\C$ has a finite open affine cover
$(V_\sigma)$ such that each $V_\sigma$ is stable under this
$\Gal$-action.  In other words, we assume the
$V_\S$ to be isomorphic to closed subvarieties of some affine space
$\AA^N(K)$, cut out by suitable polynomials in
$\Kcirc[X_\nu \mid \nu \in N]$, such that the restriction to any $V_\sigma$ of
the $\Gal$-action on $\AA^N(K)$ coincides with the restricted $\Gal$-action
on $\C$.  Finally, we assume the following compatibility condition for the
two actions: If points $C$ and $D$ in $\C$ belong to the same
$\Gamma$-orbit, then so do $C_\phi$ and
$D_\phi$ for any $\phi$, that is, the $\Gal$-action on $\C$ induces an
action of $\Gal$ on the set of orbits of $\Gamma$ in
$\C$.  We note that, given any open subset $U$ of $\C$, the closure of $U$
under the $\Gal$-action is again Zariski-open. 

 In this situation, we call a
$\Gamma$-orbit $\Gamma.C$ of $\C$ {\it generic\/} (relative to $\Kcirc$) if
$\Gamma.C$ has nonempty intersection with the
$\Gal$-closure of any dense open subset of $\C$.   

\proclaim{Lemma 4.4\. Generic orbits}  Let $\C$ be a $\Gamma$-space as
above.  Then there exists a generic orbit $\Gamma.C$ in $\C$.  It is 
unique in the following sense: Any two generic orbits $\Gamma.C$ and
$\Gamma.C'$ are equivalent under the $\Gal$-action on the set of
$\Gamma$-orbits of $\C$, i.e., there exists
$\phi \in \Gal$ such that $(\Gamma.C')_\phi = \Gamma.C$. 

More strongly: If $(\C_k)_{k \in \NN}$ is a sequence of irreducible
$\Gamma$-spaces as specified above, then given any positive integer $n$,
there exists a generic
$\Gamma$-orbit in the $\Gamma$-space $\C_1 \times \cdots
\times \C_n$. It also has the described uniqueness property.  \endproclaim 

The existence proof is standard; yet, we include the argument for the
convenience of the reader.  It is the uniqueness argument which rests on
the specifics of the setup.

\demo{Proof of Lemma 4.4} We start by ensuring existence of a generic
orbit.  To provide a setting which works for the strengthened statement
concerning a sequence
$(\C_k)$ of irreducible $\Gamma$-spaces, let $B$ be an infinite subset of
$K$ which is algebraically independent over $\Kcirc$, and $B = \bigsqcup_{i \in
\NN} B_i$ a partition of $B$ into infinite disjoint subsets $B_i$.  For each
$i$, we denote by  $K_i$ the algebraic closure of $\Kcirc(b \mid b \in B_i)$ in
$K$.  We first focus on a fixed index $i \in \NN$, and suppose
$\C = \C_i$.   

Let $(V_\S)$ be a finite open affine cover of $\C$ as specified ahead of
the lemma; in particular, each
$V_\S$ is a closed affine subspace of some affine space $\AA(K)^N$, defined
by polynomials from $\Kcirc[X_\nu \mid \nu \in N]$.  Pick one of the
$V_\S$  -- call it $V$  --  and set $V_0 = V \cap \AA(\Kcirc)^N$, the latter
being a closed affine subvariety of $\AA(\Kcirc)^N$.  We abbreviate
$\dim V = \dim V_0 = \dim \C$ to $v$.   Further, we denote by
$\Q(V)$ the  rational function field of $V$ and write the coordinate
functions in $\Q(V)$ (relative to the standard coordinatization of
$V$ in $\AA^N$) as $X_\nu$, $\nu \in N$.  Suppose
$X_{\nu_1}, \dots, X_{\nu_v}$ form a transcendence basis for
$\Q(V)$ over $K$.  The corresponding coordinate functions of $V_0$, again
denoted by
$X_{\nu_1}, \dots, X_{\nu_v}$, then form a transcendence basis for
$\Q(V_0)$ over
$\Kcirc$.  Pick arbitrary distinct elements $c_{\nu_1},
\dots, c_{\nu_v}
\in B_i \subseteq K_i$  --  by construction they are algebraically
independent over $\Kcirc$  --  and embed the subfield
$\Kcirc(X_{\nu_1}, \dots, X_{\nu_v})$ of
$\Q(V_0)$ into $K$, by sending $X_{\nu_j}$ to $c_{\nu_j}$ and mapping $\Kcirc$
identically to itself.  Since $K_i$ is algebraically closed, this map
extends to a $\Kcirc$-embedding $\rho: \Q(V_0) \rightarrow K_i \subseteq K$. 
For
$\nu \in N \setminus \{\nu_1, \dots, \nu_v \}$, set $c_{\nu} =
\rho(X_{\nu})$.  Then
$C = \bigl( c_{\nu} \bigr)_{\nu \in N}$ is a point in $V$.    

To prove that the orbit $\Gamma.C$ is generic, we will, in fact, show that
$C$ belongs to any $\Gal$-stable dense open subset $U$ of $\C$.  It is
harmless to assume $U \subseteq V$, since $U
\cap V$ is again nonempty, open, and stable under the $\Gal$-action. If
$C \in U$, we are done.  Otherwise, we consider the nonempty proper closed
subset 
$\tilde{U} = V \setminus U$ of $V$, viewing the coordinate ring of
$\tilde{U}$ as a suitable factor ring of that of $V$, and hence of the
polynomial ring $K[X_{\nu} \mid \nu \in N]$; again, we do not make a
notational distinction between the coordinate functions of
$\tilde{U}$ and the variables $X_\nu$.  Since $\dim \tilde{U}$ is strictly
smaller than $v = \dim \C$, the coordinate functions
$X_{\nu_1}, \dots, X_{\nu_v}$ of $\tilde{U}$, where $\{\nu_1, \dots,
\nu_v\}$ are as in the construction of $C$, are algebraically dependent
over $K$.  Let $F$ be an intermediate field of the extension $\Kcirc \subseteq
K$ which has finite transcendence degree over
$\Kcirc$ such that $\tilde{U}$ is defined over $F$ and the coordinate
functions $X_{\nu_1}, \dots, X_{\nu_v}$ of $\tilde{U}$ are algebraically
dependent over $F$.  Next, let $(d_{\nu_1}, \dots, d_{\nu_v})$ be any
family of elements in $K$ which are algebraically independent over
$F$ and thus, a fortiori, over $\Kcirc$.  Then there exists a point $D$ in
$V$ whose $\nu_j$-th coordinate equals $d_{\nu_j}$ for $1 \le j \le v$,
together with a $\Kcirc$-automorphism $\psi$ of $K$ which sends $d_\nu$ to
$c_\nu$ for all $\nu \in N$; such a point is obtained as in the
construction of $C$.  This means $D_\psi = C$, showing that
$C$ is $\Gal$-equivalent to $D$.  On the other hand, in view of the
algebraic dependence over $F$ of the coordinate functions
$(X_{\nu_1}, \dots, X_{\nu_v})$ of $\tilde{U}$, the point $D$ does not
belong to $\tilde{U}$.  Thus $D \in U$, whence $C \in U$, due to
$\Gal$-stability of $U$.    

The extended existence assertion of the lemma concerning generic
$\Gamma$-orbits of finite direct products of the $\C_k$ is proved
analogously, given that $B_k$ is algebraically independent over the
composite of $K_1, \dots, K_m$,  whenever
$k > m$.

To verify uniqueness of $\Gamma.C$, suppose that $\Gamma.C'$ is another
generic orbit.  We let $V$, $v = \dim V$, and $X_{\nu_1},
\dots, X_{\nu_v}$ with $\nu_1, \dots, \nu_v \in N$ be as in the
construction of $C$.  In view of its generic status, the orbit $\Gamma.C'$
nontrivially intersects $V$, since $V$ is a $\Gal$-closed dense open subset
of
$\C$.  It is therefore harmless to assume that $C' \in V$.  In a first
step, we show that the orbit
$\Gamma.C'$ contains a point $D = (d_\nu)_{\nu
\in N} \in V$ with the property that $d_{\nu_1}, \dots, d_{\nu_v}$ are
algebraically independent over
$\Kcirc$.  

Consider the intersection $U = (\Gamma.C') \cap V$, a dense open subset of
the irreducible variety $\Gamma.C'$ (keep in mind that
$\Gamma$ is connected), let $u$ be its dimension, and $\overline{U}$ the
closure of $U$ in $V$.  Moreover, let
$\Q(U) = \Q(\overline{U})$ be the function field of $U$, the latter again
expressed as the field of fractions of a suitable factor ring of the
coordinate ring of the closed subvariety $V
\subseteq \AA^N = \AA(K)^N$.  Next, we choose an algebraically closed
intermediate field $F$ of the extension
$\Kcirc \subseteq K$ which has finite transcendence degree over $\Kcirc$, such
that $\overline{U}$ is defined over $F$, and denote by
$\Q(U_F)$ the function field of the restricted variety $U_F = U \cap
\AA(F)^N$.  The coordinate functions in each of these fields will be
labeled $X_\nu$, $\nu \in N$, but will be further identified by the
function field in which they live.  For any point
$A = (a_\nu) \in V$, we finally set $\Kcirc(A) = \Kcirc(a_{\nu_j} \mid j \le v)$
and $F(A) = F(a_{\nu_j} \mid j \le v)$, both subfields of $K$.  By the
choice of $\nu_1, \dots, \nu_v$ in the existence proof, the  subfield
$\Kcirc(X_\nu \mid \nu \in N)$ of
$\Q(V)$ is algebraic over the subfield $\Kcirc(X_{\nu_1}, \dots, X_{\nu_v})$ of
this function field, and a fortiori, the subfield $F(X_\nu \mid \nu \in N)$
of $\Q(U)$ is algebraic over the subfield $F(X_{\nu_1}, \dots, X_{\nu_v})$
of the latter function field.  Consequently, we may choose a transcendence
base of $\Q(U)$ over $F$ among the coordinate functions
$X_{\nu_1}, \dots, X_{\nu_v}$.  Without loss of generality, we may assume
that  $X_{\nu_1}, \dots, X_{\nu_u} \in \Q(U)$ constitute such a
transcendence base.  As we shift to the restricted setting of $U_F$, the
coordinate functions $X_{\nu_1}, \dots, X_{\nu_u}$ still form a
transcendence base for $\Q(U_F)= \Q(\overline{U}_F)$ over $F$.  Letting
$d_{\nu_1}, \dots, d_{\nu_u}$ in
$K$ be any sequence of scalars algebraically independent over
$F$, we pick a point $D \in U$, whose coordinates labeled by
$\nu_1, \dots, \nu_u$ coincide with the given scalars  $d_{\nu_1},
\dots, d_{\nu_u}$.  Say
$D = (d_{\nu})_{\nu \in N}$ in the coordinatization of $U$ and $V$.  Then
the subfield $F(D)$ of $K$ is $F$-isomorphic to $\Q(U_F)$, via an
isomorphism which maps the class of $X_\nu$ to
$d_{\nu}$ for all $\nu\in N$. Indeed, there is an $F$-algebra homomorphism
from the coordinate ring of $\overline{U}_F$ to $F[d_\nu \mid \nu\in N]$
sending the coordinate function $X_\nu$ of $\overline{U}_F$ to $d_\nu$,
since
$D$ satisfies the defining equations of $\overline{U}_F$. This map is, in
fact, an $F$-algebra isomorphism (compare transcendence degrees), and
therefore extends to the desired field isomorphism. In particular, it
restricts to the identity on $\Kcirc$.  

Let $w$ be the transcendence degree of
$\Kcirc(D)$ over $\Kcirc$.  Then $u \le w
\le v$, and it is clearly harmless to assume that $d_{\nu_1}, \dots,
d_{\nu_w}$ form a transcendence base for $\Kcirc(D)$ over $\Kcirc$.  If $w = v$,
our intermediate claim is proved.  So let us assume $w < v$.  By
construction, this implies that the coordinate function
$X_{\nu_v} \in \Q(U_F)$ is algebraically dependent over the subfield
$\Kcirc(X_{\nu_1}, \dots, X_{\nu_w})$ of $\Q(U_F)$, and any pertinent
algebraic dependence relation yields such a relation in the field $\Q(U)$. 
In other words, the coordinate function $X_{\nu_v}$ of
$U$ is algebraic over the subfield $\Kcirc(X_{\nu_1}, \dots, X_{\nu_w})$ of
$\Q(U)$.  We thus obtain a nontrivial polynomial
$g$ in the polynomial ring
$\Kcirc[X_{\nu_i}
\mid i \in \{1, \dots, w\} \cup\{v\}]$ which is satisfied by all points in
$U$.  On the other hand, the point $C \in V$ which we constructed in the
existence proof is not a root of $g$.  Consequently, the set of nonroots of
$g$ in $V$ is dense. It is clearly $\Gal$-stable, due to the invariance of
$g$ under $\Kcirc$-automorphisms of $K$.  But, by construction, $\Gamma.C'$ has
no point in common with this set, a contradiction to the generic status of
$\Gamma.C'$.  Thus $w = v$.  
        
In the following, let $D \in U$ be such that $d_{\nu_1}, \dots, d_{\nu_v}$
are algebraically independent over $\Kcirc$.  Thus, we obtain a field
isomorphism 
$$\tau: L= \Kcirc(C) \rightarrow R = \Kcirc(D)$$  fixing the elements of $\Kcirc$
and sending $c_{\nu_j}$ to
$d_{\nu_j}$ for $1 \le j \le v$.  We let $(Y_\nu)_{\nu \in N \setminus
\{\nu_1, \dots, \nu_v\}}$ be a family of independent variables over $K$,
abbreviated to $Y$, and also denote by $\tau$ the induced isomorphism 
$L[Y] \rightarrow R[Y]$ on the level of polynomial rings.  Let
$\mu$ be any index in
$N \setminus \{\nu_1, \dots, \nu_v\}$ and $f \in L[Y]$ the minimal
polynomial of $c_\mu$ over $L$.  By construction of $C$, replacement of
each $c_{\nu_j}$ by $X_{\nu_j}$ for $j \le v$ in the coefficients of $f$
yields the minimal polynomial of $X_\mu$ over $\Kcirc(X_j \mid 1 \le j \le v)$
in the function field of $V$.  Therefore, the coordinate
$d_\mu$ of the point $D \in V$ is a root of the irreducible polynomial
$\tau(f)$, which shows $\tau(f)$ to be the minimal polynomial of $d_\mu$
over $R$.  Thus $\tau$ extends to a field isomorphism $L(c_\mu)
\rightarrow R(d_\mu)$, which sends $c_\mu$ to
$d_\mu$.  An obvious induction provides us with an extension to an
isomorphism $\Kcirc(c_\nu \mid \nu \in N) \rightarrow \Kcirc(d_\nu \mid
\nu \in N)$ sending $c_\nu$ to $d_\nu$ for all $\nu \in N$.  Algebraic
closedness of $K$ now guarantees a further extension to a
$\Kcirc$-automorphism of $K$, which we still label $\tau$.  We thus obtain $D
= C_{\tau}$, whence $\Gamma.D = (\Gamma.C)_\tau$ as required. 

The uniqueness argument for the expanded claim requires no adjustments,
since it rests on a restriction of the focus to a finite subset of $\NN$.  
\qed \enddemo

Theorem 4.3 is readily deduced from Lemma 4.4 and its proof.

\demo{Proof of Theorem 4.3} We take $\Gamma =
\autlap$.  Given an irreducible component $\C$ of
$\grassSS$, together with the affine open cover $(V_\S)$ consisting of all
nonempty intersections $\C
\cap \grassS$, the hypotheses of Lemma 4.4 are satisfied.  Indeed, whenever
$\C\cap \grassS \ne
\varnothing$, this intersection is dense open in $\C$ and closed in $V_\S$;
so, in particular, it is again affine.  It is defined over $\Kcirc$ in the
standard coordinatization of $\grassS$, since both $\C$ and $\grassS$ are. 
If
$\autlap.D$ is a generic orbit as guaranteed by Lemma 4.4  --  say $D \in
V_\S$  --  the module $P/D$ is $\SS$-generic for $\C$ by definition, and the
proof of the lemma yields a point $C = (c_\nu)_{\nu \in N}
\in \autlap.D$ such that the transcendence degree of $\Kcirc(c_\nu \mid \nu \in
N)$ over $\Kcirc$ equals the dimension of $\C$.  This provides us with a
presentation $P/C$ of the generic module
$P/D$ as postulated in Theorem 4.3.  The remaining  claims are now
straightforward in light of Lemma 4.4.  \qed
\enddemo

The $\SS$-generic modules can be constructed from $Q$ and $I$, provided we
are handed a sufficiently large subset of $K$ which is algebraically
independent over $\Kcirc$.  As described in Observation 3.6, the irreducible
components of $\grassSS$ are determined by those of the covering
$\grassS$, and the latter can be computed algorithmically in terms of
generators (defined over $\Kcirc$) for their prime ideals; a program
implementing this algorithm is available at \cite{\codes}.   Suppose
$\C$ is an irreducible component of some
$\grassSS$ and $\S$ a skeleton with $\C \cap \grassS \ne
\varnothing$.  Then a projective presentation of the generic module $G$ for
$\C$ is available via the coordinate ring of the irreducible affine variety
$\C
\cap \grassS$.  Of course, a detailed analysis of representation-theoretic
features of $G$ hinges on the understanding of that coordinate ring.

Let $\C$ be an irreducible component of $\grassSS$, and $C$ a point in
$\C$, say $C = (c_\nu)_{\nu \in N} \in \C\cap \grassS$, where $N$ again
stands for the set of all pairs $(\alpha p\zhat_r, q\zhat_s)$ with a
$\S$-critical path
$\alpha p\zhat_r$ in the first slot and a path $q\zhat_s \in \S(\alpha
p\zhat_r)$ in the second.  If
$\Kcirc(c_\nu
\mid \nu \in N)$ has transcendence degree $\dim \C$ over $\Kcirc$, then any
module $G \cong P/C$ is
$\SS$-generic for $\C$ by Theorem 4.3.  In this situation, we call $P/C$ a
{\it generic presentation\/} of $G$.  Conversely, Theorem 4.3 gurantees
that every generic module for $\C$ has a generic presentation.  For a
simple example, consider $\la = KQ$, where $Q$ is
$1 \rightarrow 2$, and $\SS = (S_1^2, S_2)$; the only arrow of the quiver is
labeled $\gamma$.  Given a sequence of top elements for the projective cover
$P = (\la e_1)^2$ of the top $S_1^2$, say
$\{z_1 = e_1 z_1, z_2 = e_1 z_2\}$, the only module represented by
$\grassSS$, up to isomorphism, is $M = \la z_1 \oplus
\bigl( \la z_2/ J z_2 \bigr)$.  So $M$ is obviously $\SS$-generic. 
However, the given presentation is not generic, since the orbit dimension of
$M$ is 
$\dim \grassSS = 1$, while for either of the two possible skeleta $\S$ of
$M$ (one is $\{\zhat_1 = e_1 \zhat_1, \gamma \zhat_1, \zhat_2 = e_1
\zhat_2\}$, the other is symmetric), the $\S$-coordinates of the submodule
we factored out of
$\la z_1 \oplus \la z_2$ to obtain $M$ belong to the subset $\{0,1\}$ of 
$K$.  A generic presentation is, for instance,
$M \cong (\la z_1 \oplus \la z_2)/
\la (\gamma z_1 +  k \gamma z_2)$, where $k \in K$ is transcendental over
$\Kcirc$. 

We glean that the set of all points $C \in \C$ which give rise to generic
presentations (of generic modules for $\C$) is not $\autlap$-stable
in general, nor need it contain a nonempty open subset of $\C$ (in fact, not
even its closure under orbits contains a nonempty open subset of $\C$, in
general).  On the other hand, as is readily seen, this set {\it is\/} always
dense in $\C$.  We record the positive observations for future reference.   

\proclaim{Corollary 4.5\. Density of generic modules}  Retain the hypotheses
and notation of Theorem {\rm 4.3}, and suppose $G$ is $\SS$-generic for
$\C$.  Then $G$ has a generic presentation $G \cong P/C$, and the
$\Gal$-orbit of $C$ is dense in
$\C$.  In other words, any nonempty open subset $U$ of $\C$ contains a
point $D$ such that
$P/D$ is a generic presentation of a generic module for $\C$.  \qed 
\endproclaim  

In light of Sections 2 and 3, the list of all
$\SS$-generic modules (where
$\SS$ traces the $d$-dimensional semisimple sequences) includes the
$(T,d)$-generic modules, and those, in turn, include the
$d$-generic ones.  In the next corollary we therefore need not separate the
cases where  
$\C$ is an irreducible component of $\grassSS$, $\grasstd$, or
$\modlad$.
 
\proclaim{Corollary 4.6\. Syzygies and AR quivers}  Suppose that the module
$G$ is generic for some irreducible component $\C$ of $\grassSS$, or
$\grasstd$, or
$\modlad$. 

Then the syzygies of
$G$ are generic for the syzygies of the modules in
$\C$ in the following weakened sense:  If $\frak U \subseteq \C$ is the
subset consisting of the points $D$ such that $\Omega^k(P/D)$ is
$\Gal$-equivalent to $\Omega^k(G)$ for all $k \ge 0$, then $\frak U$
is dense in $\C$ {\rm {(}}however, $\frak U$ need not contain a nonempty
open subset of $\C${\rm {)}}.  

Suppose, in addition, that $G$ is indecomposable, and let $\A(G)$ be that
connected component of the Auslander-Reiten quiver of
$\lamod$ which contains $G$.  Then $\A(G)$ is generic for $\C$ in the same
sense:  Namely, there exists a dense subset $\V$ of $\C$ such that, for all
$D \in \V$, the module $P/D$ is indecomposable and the connected component
$\A(P/D)$ is $\Gal$-equivalent to $\A(G)$ in the sense of Definition
{\rm4.1}.    
\endproclaim

\demo{Proof}  It suffices to apply Corollary 4.5, to obtain a dense subset
$W \subseteq \C$ such that $P/D$ is $\Gal$-equivalent to
$G$ for all $D \in W$.  This means that there exists a map $\phi \in
\Gal$ such that $P/D \cong F_\phi (G)$ (see Observation 4.1 for notation). 
Since $F_\phi$ induces a category self-equivalence of
$\lamod$, both claims follow.  \qed
\enddemo    
        
All of the generic properties to which we will apply Theorem 4.3 are
actually stable under arbitrary automorphisms of the field $K$, so that the
choice of
$\Kcirc$ becomes irrelevant for practical purposes.  If the subfield $\Kcirc
\subseteq K$ is chosen larger than necessary, this simply excludes certain
elements in $K$ from being eligible as scalars in generic presentations. 
In the second installment of consequences to Theorem 4.3 and Lemma 4.4, the
direct sum of two semisimple sequences
$\SS$ and
$\SS'$ is defined in the obvious way, as resulting from componentwise
summation. 

 Part (a.3) of the next corollary is known; see \cite{\CBS}.  While most of
the assertions follow from the statement of Theorem 4.3, the final assertion
of part (d) is a consequence of its proof.  
  
\proclaim{Corollary 4.7\. Generic properties} Let $\C$ be an irreducible
component of
$\grassSS$ or $\grasstd$ or $\modlad$, and let $\C'$ be an irreducible
component of any other parametrizing variety ``on the same level" {\rm
{(}}referring to the radical layering, or top combined with dimension, or
else to dimension without further specification{\rm {)}}.  Moreover, let
$G$ be generic for
$\C$, and $G'$ generic for $\C'$.  Then: 
\medskip

\noindent{\rm{\bf{(a)}}}  The number of indecomposable summands of
$G$ is the generic number for $\C$, and the dimension vectors of these
direct summands are the generic ones.  In fact, the collection of radical
layerings of the indecomposable direct summands of $G$ is the generic one
for the indecomposable decompositions of the modules in $\C$.   

Moreover, given an indecomposable decomposition $\,G = M_1 \oplus \cdots
\oplus M_s$, with $\,T_i = M_i/J M_i$ and $\dim M_i = d_i$, the following
additional statements hold:
\smallskip

\noindent {\rm {(a.1)}}  If $G$ is $\SS$-generic, the $M_i$ are
$\SS(M_i)$-generic for suitable irreducible components of the
$\grass{\SS(M_i)}$. 
\smallskip

\noindent {\rm {(a.2)}} If $G$ is $(T,d)$-generic, the $M_i$ are
$(T_i, d_i)$-generic for suitable irreducible components of the
$\grass^{T_i}_{d_i}$.  
\smallskip

\noindent {\rm {(a.3)}} If $G$ is
$d$-generic, the $M_i$ are $d_i$-generic for suitable irreducible
components of the
$\mod _{d_i}(\la)$.  
\medskip

\noindent{\rm{\bf{(b)}}} The dimension of the socle of $G$ is the minimum
of the socle dimensions of the modules in $\C$.
\medskip

\noindent{\rm{\bf{(c)}}} $\pdim G = \min\{ \pdim M \mid M \in \C \}$ and
$\idim G = \min\{ \idim M \mid M \in \C \}$.  Moreover, 
$$\dim
\End_\la(G,G) = \min\{\dim \End_\la(M,M) \mid M \in \C\}$$
and
$$\dim \Ext^i_\la(G,G) = \min\{\dim \Ext^i_\la(M,M) \mid M \in \C\},$$ 
for $i \ge 1,$ 
with an analogous result holding for $\Tor$-spaces.   

In particular, the dimension of the $(\autlap$- or $\GL_d$-$)$
orbit corresponding to
$G$ is maximal among the orbit dimensions in $\C$.
\medskip

\noindent{\rm{\bf{(d)}}}
If we choose $G$ and
$G'$ with the relative generic properties spelled out in {\rm Supplement 2}
to  {\rm Theorem 4.3}, then
$$\dim \Ext^i_\la (G, G') = \min\{\dim
\Ext^i_\la(M,N) \mid M \in \C,\, N \in \C' \},$$ and
$$\dim \Tor_i^\la(G, G') = \min\{\dim
\Tor_i^\la (M,N) \mid M \in \C,\, N \in \C' \}$$ for $i \ge 0$,
whenever $\C \ne \C'$.

If $\C = \C'$, then the analogous result holds, provided that $G' =
G_{\phi}$, where $G = P/C$ with $C = (c_\nu)_{\nu \in N}$, and the join of
the fields $\Kcirc(c_\nu \mid \nu \in N)$ and $\Kcirc(\phi(c_\nu) \mid \nu \in N)$
has transcendence degree $2 \cdot \dim \C$ over $\Kcirc$.
\medskip

\noindent{\rm{\bf{(e)}}} The set of skeleta of $G$ equals the generic set
of skeleta of the modules represented by $\C$.  If $\C
\subseteq \grassSS$, this  set consists of all $\,\sigma$ with $\C \cap
\grassS \ne \varnothing$.
\endproclaim

\demo{Proof} (a)  In \cite{\Kac}, Kac showed that every irreducible
component 
$\C$ of $\modlad$ contains a nonempty open subset $V$ such that all modules
$\tilde{M}$ in $V$ have the same number of indecomposable summands
$\tilde{M}_i$ and the dimension vectors of the latter are invariant on
$V$.  Combined with the arguments of Section 2, this yields a nonempty open
$U
\subseteq V$ such that even the radical layerings $\SS(\tilde{M}_i)$ of the
indecomposable summands $\tilde{M}_i$ of the $\tilde{M}$ in $U$ are
constant.  Further, we deduce that the same statements hold for irreducible
components of
$\toptd$ (resp.
$\grasstd$) and $\laySS$ (resp. $\grassSS$).   Clearly, all these
quantities are
$\Gal$-stable, whence the first two assertions follow from Theorem 4.3.  

For (a.1), we set $\SS_i = \SS(M_i)$.  Moreover, we choose a point $C
\in
\C$ which represents $G$, next to an open subset $U \subseteq \C$ as
specified above.  In particular, each module $\tilde{M}$ in $U$ can be
written in the form $\tilde{M} = \tilde{M}_1 \oplus \cdots \oplus
\tilde{M}_s$ such that $\tilde{M}_i$ is indecomposable with top $T_i$ and
dimension $d_i$.  Without loss of generality, $C = \bigoplus_{1
\le i \le s} C_i$ with $C_i \subseteq JP_i$, where $P =
\bigoplus_{1 \le i \le s} P_i$ is a suitable decomposition of $P$.  We take
$P_i$ as the distinguished projective cover of $T_i$ in the realization of
$\grass^{T_i}_{d_i}$ and consider the morphism
$\Psi$ from the projective variety $\prod_{1 \le i \le s}
\grass^{T_i}_{d_i}$ to $\grasstd$, defined by 
$\Psi(D_1, \dots, D_s) = \bigoplus_i D_i$ whenever the $D_i$ are submodules
of the $JP_i$ of codimension $d_i$ in $P_i$, respectively.  Clearly, $\Psi$
induces an isomorphism between the direct product and a suitable closed
subvariety $W$ of $\grasstd$, which restricts to an isomorphism between
$\prod_{1 \le i \le s} \grass{\SS_i}$ and $W \cap
\grassSS$.  Since
$U$ is an irreducible open subset of $\grassSS$, the intersection
$U \cap W$ is irreducible and open in $W$.  Moreover, our choice of $U$
guarantess that every $\autlap$-orbit $Z$ in $\grassSS$ which has nonempty
intersection with $U$ intersects nontrivially with $U \cap W$, whence the
closure of $Z \cap W$ in $\grassSS$ coincides with the closure of $Z$.  We
infer the existence of irreducible components
$\C_i$ of $\grass{\SS_i}$ such that $U\cap W
\subseteq \Psi\bigl(\prod_{1 \le i \le s} \C_i \bigr) \subseteq \C$ and
conclude that the closure of $\Psi\bigl(\prod_{1 \le i \le s} \C_i
\bigr)$ in
$\grassSS$ equals $\overline{U \cap  W} = \overline{U} = \C$.  For each
$i$, we now pick an $\SS_i$-generic module $G(\SS_i,
\C_i)$ for $\C_i$, such that these generic objects even enjoy the relative
generic properties described in the final statement of Theorem 4.3.  By
construction, the direct sum 
$$\bigoplus_{1 \le i \le s} G(\SS_i, \C_i)$$ is $\SS$-generic for $\C$.  In
light of the uniqueness part of  Theorem 4.3, we deduce that this direct
sum is $\Gal$-equivalent to $G$, say isomorphic to
$G_\tau \cong \bigoplus_{1 \le i \le s} (M_i)_\tau$ for some 
$\tau \in \Gal$.  Using the Krull-Schmidt theorem, we conclude that
$G(\SS_i, \C_i)$ is isomorphic to $(M_{\pi(i)})_\tau$ for some permutation
$\pi$ of $\{1, \dots, s\}$, with the property that
$\SS\bigl( M_{\pi(j)} \bigr) = \SS(M_j)$ for all $j$.  This shows the
$(M_i)_\tau$, and therefore also the $M_i$, to be $\SS_i$-generic for the
components $\C_{\pi(i)}$.

Parts (a.2) and (a.3) are proved similarly. 
\smallskip

Part (b), as well as the statements regarding $\Hom$-dimensions under (c)
and (d), follow from upper semicontinuity of the function $\dim
\Hom(-, -)$ on $\C \times \C'$.  For part (b) note, moreover, that $\dim
\Hom_\la(S_i,M)$ is the dimension of the $S_i$-homogeneous component of
$M$, the automorphisms of $K$ leave the $S_i$ invariant, up to
isomorphism (cf\. Observation 4.1), and $\dim \Hom_\la(N,M) =  \dim
\Hom_\la\bigl(N_\phi,M_\phi\bigr)$ for all $N \in \lamod$ and $\phi \in
\Gal$.
\smallskip

Concerning the remaining statements under (c) and (d):  We use upper
semicontinuity of the functions $\pdim(-)$ and $\idim(-)$ to ascertain that
the minimal values on $\C$ are the generic ones (see \cite{\Sma, Corollary
5.4}), together with the fact that both $\pdim M$ and $\idim M$ are
$\Gal$-invariant attributes of a $\la$-module $M$.  Moreover, in light of
upper semicontinuity of
$\Hom_\la(-,-)$, Propositions 5.3 and 5.6 in \cite{\Sma} (attributed to
Bongartz), entail upper semicontinuity of $\Ext_\la ^i(-,-)$ and $Tor_i
^\la(-,-)$ (see also \cite{CBS, Lemma 4.3} for $\Ext$).  Given that
application of an automorphism
$\phi$ to both arguments leaves the dimensions of the resulting $\Ext$ or
$\Tor$-spaces unchanged, Theorem 4.3 yields all assertions under (c) and
(d), except for the final one under (d).  For the latter, we refer to the
construction in the proof of Lemma 4.4.
\smallskip

For (e), it suffices to observe that the collection of skeleta of the
modules in
$\C$ is generic and again invariant under application of automorphisms
$\phi \in \Gal$, by Observation 4.1.   
\qed  
\enddemo

The precaution ``$\C' \ne \C$" in part (d) of the preceding corollary is
not redundant:  Indeed, let $\la = \CC Q$ be the complex Kronecker algebra,
that is, $Q$ has the form $\xymatrix{1 \ar@/_/[r]_{\alpha}
\ar@/^/[r]^{\beta} & 2}$, and take
$\SS = \SS' =  (S_1, S_2)$.  Then $G = \la e_1 / \la(\beta - \pi \alpha)$
is an $\SS$-generic modules.  Clearly, $\End_\la(G) = K$, while
$\Hom_\la(M, M') = 0$ for all modules $M \not\cong M'$ with radical
layering $\SS$.

\definition{Example 4.8} Let $\la = KQ/I$, where $Q$ is the quiver

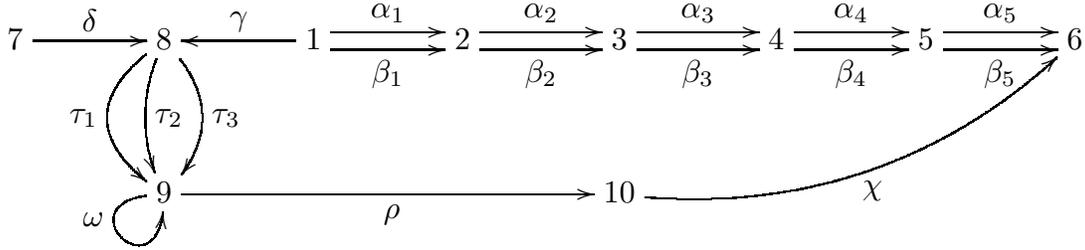
\ignore
$$\xymatrixrowsep{3pc}\xymatrixcolsep{3pc}
\xymatrix{   7 \ar[r]^{\delta}  &8 \ar@/_1.5pc/[d]_{\tau_1}
\ar@/_/[d]^{\tau_2}
\ar@/^1pc/[d]^{\tau_3}  &1 \ar[l]_{\gamma} \ar[r]<0.66ex>^{\alpha_1}
\ar[r]<-0.66ex>_{\beta_1}  &2 \ar[r]<0.66ex>^{\alpha_2}
\ar[r]<-0.66ex>_{\beta_2}  &3 \ar[r]<0.66ex>^{\alpha_3}
\ar[r]<-0.66ex>_{\beta_3}  &4 \ar[r]<0.66ex>^{\alpha_4}
\ar[r]<-0.66ex>_{\beta_4}  &5 \ar[r]<0.66ex>^{\alpha_5}
\ar[r]<-0.66ex>_{\beta_5}  &6  \\
 &9 \ar@(l,d)_(0.33){\omega} \ar[rrr]_{\rho}  &&&10
\ar@/_1.5pc/[urrr]_{\chi} }$$
\endignore

\noindent and $I \subseteq KQ$ is the ideal generated by the following
relations:
$\alpha_5 \alpha_4 \alpha_3 \beta_2 \beta_1 - \beta_5 \beta_4
\beta_3 \alpha_2 \alpha_1$, $\ \alpha_5 \alpha_4 \alpha_3 \alpha_2
\beta_1 - \alpha_5 \alpha_4 \alpha_3 \beta_2 \alpha_1$, $\ \alpha_5
\alpha_4 \beta_3 \alpha_2 \alpha_1 - \alpha_5 \beta_4 \alpha_3 \alpha_2
\alpha_1$, $\ \alpha_5 \beta_4 \alpha_3 \alpha_2 \alpha_1 - \beta_5
\alpha_4 \alpha_3 \alpha_2 \alpha_1$,  $\tau_i \gamma$ for $i = 1,2,3$, $\
\sum_{1 \le i \le 3} \tau_i \delta$,\  and $\ \omega^3$.  Then $L = 6$, and
we may take $\Kcirc$ to be the algebraic closure of the prime field of $K$. 
First choose 
$$\SS = \bigl( S_1 \oplus S_7, S_2 \oplus S_8, S_3 \oplus S_9^2, S_4
\oplus S_9, S_5 \oplus S_9, S_6 \bigr).$$   Then $\grassSS$ is irreducible
and coincides with $\grassS$, where $\S$ is the skeleton  $\S^{(1)} \sqcup
\ \S^{(2)}$ such that
$\S^{(1)}$ is the set of all initial subpaths of $\alpha_5 \alpha_4
\alpha_3 \alpha_2 \alpha_1 \zhat_1$ and
$\S^{(2)}$ the set of all initial subpaths of $\delta \tau_1 \omega_2
\zhat_2$ and $\delta \tau_2 \omega \zhat_2$.  One computes $\grassS
\cong V(Y^2 - X^3) \times \AA^2$.  We give a generic presentation of ``the"
generic module
$G$   --  we know it to be unique up to $\Gal$-equivalence  --  relative to
$\sigma$, and follow with the graph of $G$ relative to the same skeleton. 
Given scalars $c_1, c_2, c_3 \in K$ which are algebraically independent
over $\Kcirc$ and a projective cover $P = \la z_1
\oplus \la z_2$ of the top
$S_1 \oplus S_7$ of $\SS$ with $z_1 = e_1 z_1$ and $z_2 = e_7 z_2$, we
obtain: $G \cong P/(C_1 \oplus C_2)$, where $C_1$ is the submodule of
$\la z_1$ generated by 
$$\beta_1 z_1 - (\sqrt{c_1})^3 \alpha_1 z_1,\  
\beta_2 \alpha_1 z_1 - (\sqrt{c_1})^3 \alpha_2 \alpha_1 z_1, \
\beta_3 \alpha_2 \alpha_1 z_1 - c_1 \alpha_3 \alpha_2
\alpha_1 z_1,$$  
$$\beta_4 \alpha_3 \alpha_2 \alpha_1 z_1 - c_1
\alpha_4 \alpha_3 \alpha_2 \alpha_1 z_1,\ \ \beta_5 \alpha_4 \alpha_3
\alpha_2 \alpha_1 z_1 - c_1
\alpha_5 \alpha_4 \alpha_3 \alpha_2 \alpha_1 z_1;$$ here $\sqrt{c_1}$ is a
square root of $c_1$ in $K$; the submodule $C_2$ of $\la z_2$ is generated
by $\sum_{1 \le i \le 3} \tau_i \delta z_2$ and $\ \omega \tau_2 z_2 - c_2
\omega \tau_1 \delta z_2 - c_3 \omega^2
\tau_1 \delta z_2$.  The corresponding hypergraph is displayed at the end
of the discussion.

 Clearly, the generic module $G$ is decomposable into two local direct
summands
$\la z_1 / C_1$ and $\la z_2 / C_2$, these being $\SS(\la z_1 / C_1)$- and
$\SS(\la z_2 / C_2)$-generic, respectively, by Corollary 4.7.  In fact,
this example shows that not even for the generic module does every
$\S$-critical path
$\alpha p\zhat_r$  with $\S(\alpha p\zhat_r) \ne \varnothing$ lead to a
hyperedge in the corresponding hypergraph.  For instance, the $\S$-critical
path
$\gamma (e_1\zhat_1)$ leads to $\S(\gamma \zhat_1) = \{\delta
\zhat_7\}$, but the corresponding coefficient disappears in every module
with skeleton
$\S$, due to the relations
$\tau_i \gamma = 0$ in $\la$; hence, these latter relations are responsible
for the decomposability of $G$.  The socle of $G$ is
$S_6 \oplus S_9^2$, whence $S_6 \oplus S_9^2$ is contained in the socles of
all modules in $\grassSS$.  Moreover, generically, the modules with radical
layering $\SS$ have endomorphism rings of dimension at least $2$.  Further,
we observe that $G$ is also
$(T,d)$-generic, where $T = S_1
\oplus S_7$ and $d = 12$; in fact, $G$ is even $12$-generic.  While
$\grasstd$ is still irreducible, making $G$ the only
$(T,d)$-generic module, up to equivalence  --  it obviously has class
$0$ in the terminology of Observation 2.1  --  the full variety
$\modlad$ has many additional irreducible components; the one represented
by $G$ has class $1$.  Note that the
set of skeleta of
$\la z_1 / C_1$ has cardinality $2^5$, while $\la z_2/ C_2$ has $6$
different skeleta.  The singularities of $\grassS$, here identified with
$V(Y^2 - X^3) \times \AA^2$, are precisely the points that parametrize the
modules with fewer than $2^5 \cdot 6$ skeleta; in fact, each of the latter
modules has precisely $6$ skeleta.  Moreover, by Theorem 4.3, the syzygies
$\Omega^i(G)$ are the generic syzygies for the irreducible component
$\C$ of $\modlad$ obtained by closing
$\laySS$ in $\modlad$; they are obtainable combinatorially from the graph
of $G$, showing that $\pdim G = \infty$.  By Corollary 4.7, we conclude
further that all modules in $\C$ have infinite projective dimension. 

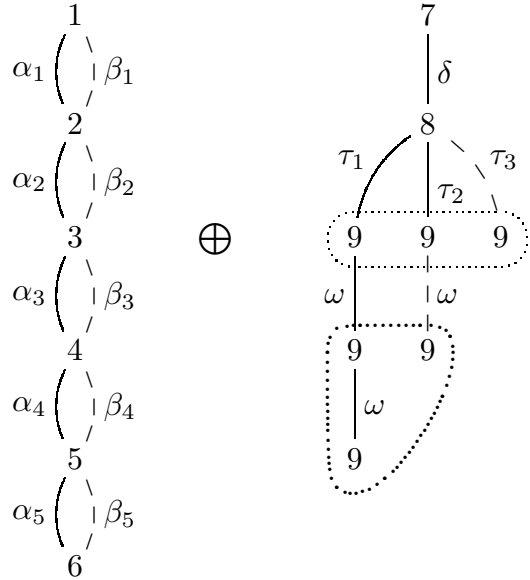
\ignore
$$\xymatrixrowsep{1.8pc}\xymatrixcolsep{1pc}
\xymatrix{  1 \edge@/_/[d]_{\alpha_1} \dashedge@/^/[d]^{\beta_1}  &&&&&7
\edge[d]^\delta  \\ 2 \edge@/_/[d]_{\alpha_2} \dashedge@/^/[d]^{\beta_2} 
&&&&&8
\edge@/_/[dl]_{\tau_1} \edge[d]^(0.6){\tau_2} \dashedge@/^/[dr]^{\tau_3} 
\\ 3 \edge@/_/[d]_{\alpha_3} \dashedge@/^/[d]^{\beta_3}  &&\bigoplus &&9
\levelpool2 \edge[d]_{\omega} &9 \dashedge[d]^{\omega} &9  \\ 4
\edge@/_/[d]_{\alpha_4} \dashedge@/^/[d]^{\beta_4}  &&&&9
\save[0,0]+(-3,0);[0,0]+(-3,0) **\crv{~*=<2.5pt>{.} [0,0]+(-3,3)
&[0,0]+(0,3) &[0,1]+(0,3) &[0,1]+(3,3) &[0,1]+(3,-4) &[1,0]+(3,-4)
&[1,0]+(-3,-4) &[1,0]+(-3,0) }\restore
\edge[d]^{\omega} &9  \\ 5 \edge@/_/[d]_{\alpha_5}
\dashedge@/^/[d]^{\beta_5}  &&&&9  \\ 6 }$$
\endignore            
 
Finally, we remark: If we replace $\SS$ by the sequence $\SS'$, which
differs from
$\SS$ only in the penultimate slot, taking $\SS'_4 = \SS_5 \oplus S_{10}$
versus $\SS_4 = S_5 \oplus S_9$, then the $\SS'$-generic module is
indecomposable and has socle $S_6 \oplus S_9$.  \qed
\enddefinition

\head 5.  Irreducible components and generic modules over truncated path
algebras
\endhead

{\it From now on, we assume that $\la$ is a truncated path algebra, i\.e\.,
$\la = KQ/I$, where $I = I(L) \subseteq KQ$ is the ideal generated by all
paths of length $L+1$ for some $L\ge1$. We will keep $L$ fixed in our
discussion. Note that in this situation, path lengths in $\la$ are well
defined, whence it is unnecessary to pass from projective
$\la$-modules $P$ to projective $KQ$-modules $\Phat$ in defining skeleta
{\rm(}cf\. Definition {\rm3.1)}. This means that in the following, all
$\zhat_r$ may be replaced by $z_r$, if desired.
\/}  
\medskip

As we will see in Theorem 5.3 below, in this situation, the $\grassSS$ are
irreducible smooth rational varieties.  To describe their structure in
greater detail, the fact that $\la$ is graded by path lengths will be
pivotal.  In fact, it will be advantageous to simultaneously explore, from
the start, varieties of graded modules with radical layering
$\SS$, next to the full class of modules parametrized by
$\SS$.  We recall from
\cite{\grassgraded} the tools required for studying the homogeneous points
in
$\grassSS$. 

As before, let $\SS$ be a $d$-dimensional semisimple sequence with top
$T = \SS_0$ and $P$ a projective cover of $T$, endowed with its
distinguished sequence $z_1, \dots, z_{\t}$ of top elements.  Moreover, let
$P = \bigoplus_{0 \le l \le L} P_l = \bigoplus_{1 \le r \le \t}
\bigl( \bigoplus_{0 \le l \le L} \la_l z_r \bigr)$ be the natural grading
of $P$, where $\la_l$ denotes the
$K$-subspace of $\la$ spanned by the paths of length $l$ and $P_l =
\bigoplus_{1 \le r \le \t} \la_l z_r$. The {\it homogeneous points\/} of
$\grasstd$ are those of the form
$C = \bigoplus_{1 \le l \le L} C_l$, where $C_l = C \cap P_l$.  For proofs
of the following observations, see \cite{\grassgraded}.  

\definition{Definitions and Observations 5.1}

{\rm (1)}  The set $\ggrasstd$ of all homogeneous points in
$\grasstd$ is a closed subvariety.  In particular, $\ggrasstd$ is in turn
projective.

{\rm (2)}  The set $\ggrassSS$ of all homogeneous points in
$\grassSS$ is a {\it projective\/} subvariety of $\grassSS$; in particular,
$\ggrassSS$ is closed in $\grassSS$.  (Note that, by contrast, $\grassSS$
fails to be projective, in general.)

{\rm (3)}  For each abstract skeleton $\S$ compatible with $\SS$, the set
$\ggrassS$ of all homogeneous points in $\grassS$ is an affine open
subvariety of $\ggrassSS$, and the family of these subvarieties, as $\S$
traces all eligible skeleta, covers $\ggrassSS$.

{\rm (4)}  The subgroup $\gautlap$ of $\autlap$, consisting of the
homogeneous automorphisms of $P$, acts morphically on $\ggrasstd$ and
$\ggrassSS$.  The orbits of $\ggrasstd$ under this action are in one-to-one
correspondence with the graded-isomorphism classes of those graded
$d$-dimensional modules with top $T$ which are generated in degree zero;
analogously, the $\gautlap$-orbits of
$\ggrassSS$ are in one-to-one correspondence with the graded-isomorphism
classes of those graded modules which have radical layering $\SS$ and are
generated in degree zero.
\enddefinition

For truncated path algebras $\la$, it is exceptionally easy to recognize
the semisimple sequences $\SS$ for which $\grassSS \ne
\varnothing$ or --  as it turns out, equivalently  --  $\ggrassSS \ne
\varnothing$.  Indeed, the final two of the equivalent conditions below are
immediately checkable via the (tree) graphs of the indecomposable
projective modules
$\la e_i$, and these graphs can in turn be simply read off the quiver.  
      
\proclaim{Observation 5.2} Suppose $\la = KQ/I$ is a truncated path
algebra, and $\SS = \bigl(\SS_0, \dots, \SS_L \bigr)$ a semisimple
sequence.  Then the following conditions are equivalent: 
\smallskip
\noindent {\rm (a)} $\grassSS$ is nonempty.

\noindent {\rm (b)} $\ggrassSS$ is nonempty.

\noindent {\rm (c)} There exists an abstract skeleton compatible with
$\SS$. 

\noindent {\rm (d)}  For each positive integer $l < L$, the following
holds:    
$$\sum_{1 \le i \le n} \bigl(\text{the number of arrows in}\ Q \ 
\text{from}\  e_i \ \text{to}\  e_j \bigr) \cdot \dim e_i \SS_l  \ \ \ \ge
\ \ \ \dim e_j \SS_{l+1},$$ 
\qquad  for all $j = 1, \dots, n$.  
\medskip

\noindent Moreover, any abstract skeleton $\S$ consisting of paths of
length $\le L$ in $Q$ is the skeleton of a {\rm {(}}graded{\rm {)}}
$\la$-module.  In other words, $\ggrassS \ne \varnothing$, and a fortiori,
$\grassS \ne \varnothing$, for every abstract skeleton
$\S$.   
\qed  \endproclaim

Condition (d) is always necessary for $\grassSS$ to be nonempty, but need
not be sufficient in the non-truncated case.  For example, consider the
algebra $KQ/I$, based on the quiver 

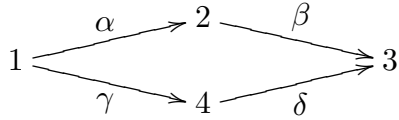
\ignore
$$\xymatrixrowsep{0.1pc}\xymatrixcolsep{4pc}
\xymatrix{  
 &2 \ar[dr]^{\beta} \\ 1 \ar[ur]^{\alpha} \ar[dr]_{\gamma} &&3 \\
 &4 \ar[ur]_{\delta} }$$
\endignore

\noindent where $I$ is generated by $\beta \alpha - \delta \gamma$.  Then
the sequence
$\bigl(S_1, S_2, S_3 \bigr)$ satisfies the two final conditions under 5.2,
but fails to arise as the sequence of radical layers of a
$\la$-module. 

As an obvious consequence of the following theorem, we obtain that the set
of skeleta of any $\SS$-generic module over a truncated path algebra is the
full set of skeleta compatible with $\SS$.  

For the definitions of $N(\SS)$ and $N_0(\SS)$, see Observation 3.10.
     
\proclaim{Theorem 5.3\. Irreducibility of  $\grassSS$ and 
$\ggrassSS$, and first structure results}  Suppose that
$\la$ is a truncated path algebra and
$\SS$ any semisimple sequence such that $\grassSS \ne \varnothing$.  

Then both $\grassSS$ and $\ggrassSS$ are irreducible varieties endowed with
finite open covers consisting of copies of an affine space, namely, of 
$\AA^{N(\SS)}$ in case of $\grassSS$, and $\AA^{N_0(\SS)}$ in case of
$\ggrassSS$. In particular, both
$\grassSS$ and $\ggrassSS$ are rational and smooth.   

More precisely: If
$\S$ is any skeleton compatible with
$\SS$, then $\grassS \cong \AA^{N(\SS)}$ is dense and open in
$\grassSS$, and $\ggrassS \cong \AA^{N_0(\SS)}$ is dense and open in
$\ggrassSS$.  
\endproclaim

\demo{Proof}  We only address the claims for ungraded modules.  The proofs
of the twin statements for $\ggrassSS$ are proved analogously.  However, we
will use some observations about the graded case to deal with the ungraded
one. 

As above, we let $\la = \bigoplus_{0 \le l \le L} \la_l$ and $P =
\bigoplus_{l \le L} P_l$ be the decompositions of $\la$ and $P$ into
homogeneous subspaces.  Note that, whenever $\S$ is an abstract skeleton
and $C \in \ggrasstd$, the requirement that $\S$ be a skeleton of $P/C$ is
equivalent to the formally weaker condition that $\S$ yields a basis for
$P/C$. 

Let $\S$ be compatible with $\SS$.  Whenever we do not explicitly  insist
on identification of $\S$ with a subset of $P$, we go back to the original
definition of skeleta in Definition 3.1, viewing the elements of
$\S$ as paths in the $KQ$-module $\Phat = \bigoplus_{1 \le r \le \t} (KQ)
\zhat_r$.  By $\S^{(r)}$ we denote the set of paths in $\S$ that start in
$\zhat_r$, and by $\S^{(r)}_{li}$ the subset of paths of length
$l$ in $\S^{(r)}$ which end in the vertex $e_i$. 

That $\grassS \cong \AA^{N(\SS)}$ is obvious.  Since the affine patches
corresponding to the skeleta form an open cover of $\grassSS$, it thus
suffices to show that, given any other skeleton compatible with
$\SS$  --  call it $\tilde{\S}$  --  the intersection $\grassS \cap
\grass{(\tilde{\S})}$ is non-empty;  indeed, $\grassS \cap
 \grass{(\tilde{\S})}$ is then a dense open subset of both $\grassS$ and
$\grass{(\tilde{\S})}$.  Clearly, $\S \cap\tilde{\S}$ is the disjoint union
of the intersections $\S^{(r)}_{li} \cap
\tilde{\S}^{(r)}_{li}$; the union $\S \cup \tilde{\S}$ breaks up similarly.

By induction on $L \ge 0$, we will construct a point $C \in
\ggrassSS$, that is, a homogeneous submodule $\bigoplus_{l \le L} C_l$ of
$P$, together with a bijection $f: \tilde{\S} \rightarrow \S$ whose
restriction to
$\S \cap
\tilde{\S}$ equals the identity on $\S \cap \tilde{\S}$, such that the
following additional conditions are satisfied:
\roster
\item"$\bullet$" $C \in \grassS \cap \grass{(\tilde{\S})}$; 
\item"$\bullet$"  whenever $\tilde{p} \in \tilde{\S}^{(r)}_{li}$, the path
$f(\tilde{p})$ belongs to
$\S^{(s)}_{li}$ for some $s$, and the difference  
$\tilde{p} - f(\tilde{p})$, now viewed as an element of $P$, belongs to the
homogenous subspace $C_{\len(\tilde{p})}$ of $C$.
\item"$\bullet$" $p z_r = 0$ in $P$, for any path $p \zhat_r \in \Phat$
which does not belong to $\S \cup \tilde{\S}$.   
\endroster

\noindent Note that, by the first condition, the families of residue
classes $\bigl(\tilde{p} + C \bigr)$ and $\bigl(f(\tilde{p})
 + C \bigr)$ in $P/C$ form bases for $P/C$; here we again identify $\S$ and
$\tilde{\S}$ with subsets of $P$ under the conventions of Definition 3.1. 
In particular, none of the elements
$\tilde{p}$ or $f(\tilde{p})$ of $P$ belongs to
$C$. Moreover, the second requirement forces
$f$ to induce a bijection $\tilde{\S}_{li} \rightarrow \S_{li}$, for each
$0 \le l \le L$ and $1 \le i \le n$.  

The proposed construction is trivial for
$L = 0$.  So, given $L$, we assume that our requirements can be met in any
setup of the described ilk whenever $L -1$ is an upper bound on the lengths
of nonvanishing paths.  Suppose $L' = L-1$.  Let
$\la'$ be the vector space $\bigoplus_{0 \le l \le L'} \la_l$ identified
with the algebra
$\la / J^L$, and let $P' = \bigoplus_{0 \le l \le L'} P_l$ be the
corresponding truncation of $P$, clearly a projective
$\la'$-module. (Up to isomorphism,
$P'$ is the same as $P/J^L P$ as a $\la'$-module; but as a vector space,
$P'$ is contained in $P$ by definition, so that, for a subset $U
\subseteq P'$, it makes sense to talk about the $\la$-submodule $\la U$ of
$P$ generated by $U$.)  Moreover, we set $\SS' = \bigl(\SS_0,
\dots, \SS_{L'} \bigr)$, and let $\S'$,  $\tilde{\S}'$ be the skeleta
resulting from $\S$ and $\tilde{\S}$ through deletion of the paths of
length $L$.  Observe that $\S'$ and  $\tilde{\S}'$ are both compatible with
$\SS'$.  Therefore, our induction hypothesis yields a homogeneous submodule
$C' \subseteq P'$ and a bijection $f':
\tilde{\S}' \rightarrow \S'$ satisfying the conditions listed above. 
Compatibility of $\S$ and $\tilde{\S}$ with $\SS$ allows us to extend
$f'$ to a bijection $f: \tilde{\S}
\rightarrow \S$ that coincides with the identity on $\S \cap \tilde{\S}$
and induces bijections $\tilde{\S}_{Li} \rightarrow \S_{Li}$ for $1 \le i
\le n$.  Our choice of $C'$ and $f$ then guarantees the following:  Suppose
$q$ is a path of length $L-1$, $\alpha$ an arrow and
$r \le \t$. In case $q \zhat_r \notin \S \cup \tilde{\S}$, we have
$\la \alpha q z_r = K  \alpha q z_r \subseteq \la C'$.  

To construct $C \in \grasstd$ with the required properties, we will add to
$\la C'$ suitable cyclic submodules of $P$ generated by linear combinations
of paths $q z_r$ in $P$, where $r \le \t$ and $q$ is a path of length $L$
in $KQ$.  In light of the preceding comment, our interest is focused on
paths of the form $\alpha p \zhat_r \in \Phat$, where $\alpha$ is an arrow
and $p \zhat_r$ belongs to
$\S_{L-1} \cup \tilde{\S}_{L-1}$.  We will now, for each path of this ilk,
construct a term $C(\alpha p \zhat_r) \subseteq P$ to be added to $\la C'$
so as to yield $C$ as desired.

First suppose $p \zhat_r \in \S_{L-1}$, and let $\alpha$ be an arrow such
that  $\alpha p \zhat_r$ is $\S$-critical (that is, $\alpha p
\zhat_r \notin \S$) and $\alpha f^{-1}(p \zhat_r) \notin
\tilde{\S}$.  Then also $\alpha p \zhat_r  \notin
\tilde{\S}$.  Indeed, suppose $\alpha p \zhat_r \in \tilde{\S}$.  Then
$p \zhat_r \in \S \cap \tilde{\S}$, whence $f^{-1}(p \zhat_r) = p
\zhat_r$, and consequently also $\alpha f^{-1}(p \zhat_r) = \alpha p
\zhat_r$ belongs to $\tilde{\S}$, contrary to our hypthesis.  In
particular, this means that we want the element $\alpha p z_r$ of $P$ to
belong to $C$ in order to meet the third of the conditions we are
targeting.  Accordingly, we let $C(\alpha p \zhat_r) \subseteq P$ be the
$\la$-submodule $K \alpha p z_r$ in that case. 

Now suppose $p \zhat_r \in \S_{L-1}$ such that $\alpha p \zhat_r$ is
$\S$-critical, while $\alpha f^{-1}(p \zhat_r) \in \tilde{\S}$. We define
the correction term as $C(\alpha p \zhat_r) = K\bigl(\alpha p z_r -
f(\alpha f^{-1}(p \zhat_r)) \bigr)$.  In the latter difference,
$f(\alpha f^{-1}(p \zhat_r))$ is identified with an element in $P$; that
is, 
$f(\alpha f^{-1}(p \zhat_r))$ stands for $\alpha q z_s$ if $f^{-1}(p
\zhat_r) = q \zhat_s$.  Thirdly, in case
$\alpha p \zhat_r \in \S$ but $\alpha f^{-1}(p \zhat_r) \notin
\tilde{\S}$, we set $C(\alpha p \zhat_r) = K \bigl(\alpha p z_r - f^{-1} (
\alpha p \zhat_r ) \bigr)$; in this difference, we again view 
$f^{-1} (\alpha p \zhat_r)$ as an element in $P$.  Finally, if $\alpha p
\zhat_r \in \S$ and $\alpha f^{-1}(p \zhat_r) \in
\tilde{\S}$, we take $C(\alpha p \zhat_r)$ to be zero; in that case, our
induction hypothesis guarantees that $p z_r -
 f^{-1}(p \zhat_r) \in C'$, whence  $\alpha p  z_r -
\alpha f^{-1}(p \zhat_r) \in \la C'$; once more, $f^{-1}(p \zhat_r)$ is
viewed as an element in $P$.   

Then 
$$C = \la C' \ + \ \sum_{\alpha \text{\ arrow},\ p\zhat_r\, \in \,
\S_{L-1} \cup \tilde{\S}_{L-1}} C(\alpha p \zhat_r)$$  is a $\la$-submodule
of $P$ satisfying all requirements imposed by our induction.  Indeed, 
since $C$ is homogeneous by construction, our initial comments show that,
to prove $C \in \grassSS \cap \grass{(\SS')}$, one only needs to observe
that each of the two families
$\bigl( p z_r \bigr)_{p \zhat_r \in \S}$ and $\bigl( p z_r \bigr)_{p
\zhat_r \in \tilde{\S}}$ of elements of
$P$ induces a basis for $P/C$, whence $C \in \grassSS \cap
\grass{(\SS')}$.  This completes the argument showing that $\grassS$ is
dense in $\grassSS$. 

Smoothness of $\grassSS$ now follows from the fact that the $\grassS$,
where $\S$ traces the skeleta compatible with $\SS$, form an open cover of
$\grassSS$.  Rationality of $\grassSS$ is obvious. \qed
\enddemo

To contrast Theorem 5.3 with the non-truncated situation:  In general  
{\it any\/} affine variety, not necessarily irreducible, arises as a variety
$\grassS$, up to isomorphism, for suitable $\la$ and $\S$; see
\cite{\GeomI, Theorem G}.  Moreover, arbitrary affine or projective
varieties arise in the form
$\grassSS$; see \cite{\grassgraded, Section 5} and \cite{\Hil, Example}. 

\proclaim{Corollary 5.4\. Rationality of the components of $\grasstd$}
Suppose 
$\la$ is a truncated path algebra, and consider the set
$$\{\overline{\grassSS} \mid \SS \text{\ is a\ } d\text{-dimensional
semisimple sequence with\ }\
\SS_0 = T \}$$
of closures of the varieties $\grassSS$ in $\grasstd$.  The maximal elements
in this set are precisely the irreducible components of $\grasstd$. 
They are determined by skeleta as follows:
$$\overline{\grassSS} \ = \ \overline{\grassS},$$ whenever $\S$ is a
skeleton compatible with $\SS$. 

Analogously, the irreducible components of $\ggrasstd$ are the maximal
ones among the closures of the varieties $\ggrassSS$ in
$\ggrasstd$.      

In particular, all irreducible components of $\grasstd$ and $\ggrasstd$ are
rational varieties. \qed \endproclaim

Theorem 5.3 and Corollary 5.4 can be transposed to the classical affine
setting via Proposition 2.2:  Namely, each variety $\laySS$ is irreducible
and smooth, and the irreducible components of the varieties $\toptd$ and
$\modlad$ are the maximal candidates among the closures of the
$\laySS$ in the larger varieties. 

For the next consequence of Theorem 5.3, we combine this theorem with
Observation 2.1 and Corollary 2.7. 

\proclaim{Corollary 5.5\. Number of irreducible components}
Let $\la$ be a truncated path algebra,
${\bold d}$ a dimension vector, and $\mu$ the number of
irreducible components of the variety $\mod_{{\bold d}}(\la)$
parametrizing the modules with dimension vector ${\bold d}$.  

Then $\mu$ is bounded from above by the number of semisimple sequences $\SS
= (\SS_0, \dots, \SS_L)$ with dimension vector
${\bold d}$ {\rm{(}}i.e., with  $\underline{\dim} \bigoplus_{0 \le l \le L}
\SS_l = {\bold d}${\rm {)}} such that
$\laySS \ne \varnothing$.
 \qed \endproclaim

On the other hand, the number $\mu$ of irreducible components of
$\mod_{{\bold d}}(\la)$ is bounded from below by the number of those
semisimple sequences
$\SS$ with dimension vector
${\bold d}$ which are minimal, with regard to $\laySS \ne \varnothing$, in
the domination order of Section 2; this lower bound is available over
arbitrary basic finite dimensional algebras (it is an immediate consequence
of  \cite{\grassIII, Corollary 2.12}).   Even in the case of a truncated
path algebra, $\mu$ lies strictly between the upper and lower bounds we thus
obtain, in general; see the example following Theorem 5.12. 

We remark that, for a given truncated path algebra $\la$ and semisimple
sequence $\SS$, it can be detected on sight from the quiver and Loewy length
of $\la$ whether $\laySS$ is nonempty.

While we expect the $\laySS$ to also be
rational, in tandem with the $\grassSS$, it is only in the following
slightly weakened form that we could carry over rationality to the
classical parametrizing varieties.   Fortunately, this weaker result
still suffices to yield Corollary 5.7.

For the representation-theoretic
audience, we recall that an irreducible variety $V$ over $K$ is
{\it unirational\/} if its function field embeds into a purely
transcendental extension of $K$.

\proclaim{Corollary 5.6\. Unirationality of the components of
$\toptd$ and $\modlad$}  For any positive integer
$d$, the irreducible components of the affine variety
$\modlad$ are unirational, as are the irreducible components of the
quasi-affine varieties $\toptd$.   

Moreover: Suppose $\C$ is any irreducible component of $\modlad$. If $\SS$
is the generic radical layering of the modules in $\C$, and $\S$ any
skeleton compatible with
$\SS$, then $\C$ contains a dense open subset consisting of modules with
skeleton $\S$. \endproclaim

\demo{Proof}  Let $\C$ be an irreducible component of $\modlad$ or some
$\toptd$, and let $\SS$ be the generic radical layering of the modules in
$\C$.  Since $\laySS$ is irreducible, $\C$ contains
$\laySS$ as a dense locally closed subvariety.  Let $\S$ be any skeleton
compatible with
$\SS$.

We construct a (Zariski-closed) subset, labeled $\modS$, of $\laySS$ as
follows:  First, we index a basis for $K^d$ by the elements of
$\S$, say $\bigl( b_{\phat} \bigr)_{\phat \in \S}$, and identify any map in
$\End_K(K^d)$ with its matrix relative to the sequence $(b_{\phat})$,
indexing rows and columns by the elements of $\sigma$.  Next, we consider
the
$K$-algebra generators of $\la = KQ/I$ given by (the $I$-residues) of the
elements in $Q^* = Q_0 \cup Q_1$, and define
$\modS$ as the set of those points $(x_u) \in \prod_{u \in Q^*}
\End_K(K^d)= \prod_{u \in Q^*} M_d(K)$ which satisfy the following
conditions:          For $u = e_i$, the matrix $x_u$ carries $1$ in position
$(\phat, \phat)$, whenever the path $\phat$ ends in the vertex
$e_i$, and carries $0$ in all other positions; for $u = \alpha \in Q_1$,
the $\phat$-th column of the matrix $x_u = x_\alpha$ has the following
form: in case $\alpha
\phat \in \S$, it carries
$1$ in position $(\alpha \phat, \phat)$ and zeros in the other positions,
while, in case $\alpha \phat \notin \S$ (meaning that
$\alpha \phat$ is $\S$-critical), the $\phat$-th column is required to have
zeros in all positions $(\qhat, \phat)$ having row index $\qhat$ outside
$\S(\alpha \phat)$, with no conditions imposed on the remaining entries. It
is readily seen that $\modS \subseteq \laySS$;  indeed, by construction,
the sequences of matrices in $\modS$ describe multiplication by the
elements $u \in Q^*$ on the factor modules $P/C$ for $C \in
\grassS$, relative to the basis $\bigl (p z_r + C \bigr)_{\phat = p \zhat_r
\in \S}$ of $P/C$, respectively.  Thus, Theorem 5.3 ensures that all those
positions in the matrix sequences above, which are not specified to be $0$
or $1$ in the definition of
$\modS$, will independently take on arbitrary values in $K$.  In other
words, $\modS \cong \AA^{N(\SS)}$.

Let $V$ be the closure of $\grassS$ in $\grassSS$ under the
$\autlap$-action.  Then 
$$V = \bigcup_{f \in \autlap,\, C \in \grassS} f.C$$ is open and dense in
$\grassSS$, since by Theorem 5.3 we know
$\grassS$ to be open and dense.  Proposition 2.2 therefore shows the
corresponding $\GL_d$-stable subvariety $W$ of $\laySS$ to be open and
dense in $\laySS$.  On the other hand, by construction, $W$ is the closure
of $\modS$ under the $\GL_d$-action.  This provides us with a dominant
morphism
$$\GL_d \times \modS \rightarrow \laySS,  \ \ \bigl(g, (x_u)_u \bigr)
\mapsto (g^{-1} x_u g)_u.$$  Therefore $\laySS$ is unirational, and so is
$\C$.  To conclude the argument, we observe that the $\GL_d$-orbits in $W$
are in one-to-one correspondence with the isomorphism classes of modules
having skeleton
$\S$.  \qed
\enddemo

Let us return to the hierarchy of irreducible components of
$\grasstd$ introduced in Corollary 2.7.  Over truncated path algebras, the
irreducible components of class $0$ are just the closures in $\grasstd$ of
the irreducible varieties $\grass{\SS^{(0,i)}}$,  where $\SS^{(0,1)},
\dots, \SS^{(0, n_0)}$ are the distinct semisimple sequences which are
minimal with respect to  
$\grass{\SS^{(0,i)}} \ne \varnothing$.  Hence, the irreducible components
of class $0$ of $\grasstd$ are in one-to-one correspondence with the
sequences $\SS$ that are minimal with respect to satisfying the equivalent
conditions of Observation 5.2.  They can be picked out on sight.  Next, let
$\SS^{(1,1)}, \dots, \SS^{(1, n_1)}$ be the distinct semisimple sequences
which are minimal among the sequences $\SS'$ with the property that
$\grass(\SS')$ is not contained in the closure of
$\bigcup_{1 \le i \le n_0} \grass{\SS^{(0,i)}}$ in $\grasstd$.  These
sequences are much harder to recognize in general, a task which will be
addressed in a sequel to this paper, with the aid of Corollary 5.7 below. 
The (necessarily distinct) closures of the $\grass{\SS^{(1,i)}}$ in
$\grasstd$ are precisely the irreducible components of class $1$ in the
bigger variety, and so forth.  Irreducible components of $\grasstd$ (alias
$\toptd$) of arbitrarily high class numbers are realizable over truncated
path algebras; see Examples 5.8(2) below. 

An analogous sifting process produces the irreducible components of class
$0$ of $\modlad$ as closures of certain varieties $\lay{\SS^{(0, n_0)}}$ in
$\modlad$, the $0$-th generation being again easy to detect. Recognizing
non-embeddedness in identifying the sequences that lead to irreducible
components of higher class numbers again requires additional theory
regarding closures.  The main tool in accessing such closures will result
from Corollaries 5.4 and 5.6 above, combined with the following fact, due
to J\. Koll\'ar (for a proof, see \cite{\grassII, Proposition 3.6}): 
Whenever $V$ is a unirational projective variety of positive dimension and
$x, y \in V$, there exists a curve $\psi: \PP^1 
\rightarrow V$, the image of which connects $x$ and $y$.  Applications as
indicated of the next corollary will follow in a sequel to
this article.

\proclaim{Corollary 5.7\. Curves in components} 

$\bullet$  Let $T$ be semisimple, $\C$ an irreducible component of the
projective variety  $\grasstd$, and $\SS$ the generic radical layering of
$\C$.  Moreover, suppose that
$\S$ is any skeleton compatible with $\SS$.  Then any point $C \in \C$
belongs to the image of a curve $\psi: \PP^1 \rightarrow \C$ mapping a
dense open subset of $\PP^1$ to $\grassS$.

$\bullet$  Let $\C$ be an irreducible component of the affine variety
$\modlad$, $\SS$ again its generic radical layering, and $\S$ a skeleton
compatible with $\SS$.  Then any point $x \in \C$ belongs to the image of a
rational map $\psi: \AA^1 \rightarrow \C$ which maps a dense open subset of
$\AA^1$ to the locally closed subvariety of
$\C$ consisting of the modules with skeleton $\S$. \qed 
\endproclaim   

Provided that the closure of some
$\lay{\SS^{(1,i)}}$ in $\modlad$ is maximal irreducible, its class among
the irreducible components of $\modlad$ may change in either direction when
compared with the class of its closure in the pertinent $\toptd$; see
Examples 5.8, (2) and (3).  The movement in the class numbers of
irreducible components, as one progresses to closures on the next level in
the hierarchy of parametrizing varieties, in fact, encodes a substantial
amount of information about
$\lamod$.   It warrants a separate study of these invariants, to link them
more directly to the quiver
$Q$ and the Loewy length of $\la$.   

\definition{Examples 5.8}
\smallskip
 
\noindent {\rm {(1)}}  Let $Q$ be the quiver 
\qquad\ignore
$\xymatrixrowsep{1pc}\xymatrixcolsep{2pc}
\xymatrix{   1 \ar@(ul,dl) \ar[r] &2 }$,
\endignore and take $L=2$, that is, $\la = KQ/I$, where $I$ is generated by
the paths of length $3$.  Moreover, let $d=3$, $T = S_1$, and $\SS = (S_1,
S_1, S_2)$.  Then $\grassSS \cong \AA^1$, and the closure of $\grassSS$ in
$\grasstd$ equals $\grasstd$, the latter variety being isomorphic to
$\PP^1$. Thus, the variety $\toptd$ is also irreducible.  Its closure in
$\modlad$ is an irreducible component of class $0$.  On the other hand, the
closure of the irreducible variety $\grass{(S_1, S_1 \oplus S_2, 0)}$ in
$\grasstd$ (a singleton) is not an irreducible component of $\grasstd$, as
it is properly contained in the closure of $\grassSS$. 
\smallskip

\noindent {\rm {(2)}} For every nonnegative integer $h$, there exist a
truncated path algebra $\la = KQ/I$ and a semisimple sequence
$\SS$, with top $T$ and dimension $d$ say, such that the closure of
$\grassSS$ in $\grasstd$ is an irreducible component of class $h$.  We give
examples for $h = 1, 2$, which make it clear how to move on to higher
values of $h$.

Let $Q$ be the quiver 

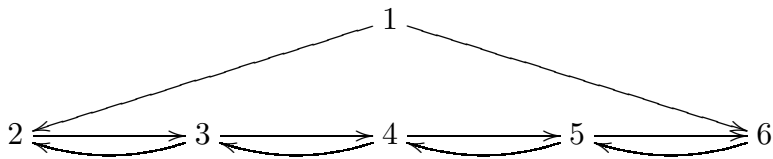
\ignore
$$\xymatrixrowsep{2pc}\xymatrixcolsep{4pc}
\xymatrix{  
 &&1 \ar[dll] \ar[drr] \\ 2 \ar[r]  &3 \ar[r] \ar@/^/[l]  &4 \ar[r]
\ar@/^/[l]  &5 \ar[r] \ar@/^/[l]  &6 \ar@/^/[l] }$$
\endignore

\noindent  and $\la$ the truncated path algebra with quiver
$Q$ and Loewy length $L+1 = 6$.  Moreover, take $T = S_1^2$ and $d = 7$. 
Then the closure of 
$$\grassSS : = \grass{(S_1^2, S_2, S_3, S_4, S_5, S_6)}$$ in $\grasstd$ is
an irreducible component of class $0$ of
$\grasstd$, clearly the only one.  The closure of 
$$\grass{(\SS')} : = \grass{(S_1^2, S_2 \oplus S_6, S_3, S_4, S_5, 0)}$$ 
in $\grasstd$ is an irreducible of component of class $1$: Indeed, given
that $\SS <
\SS'$, we only need to show that $\grass{(\SS')}$ is not contained in the
closure of $\grassSS$ in
$\grasstd$; but this follows from the fact that the arrow $6
\rightarrow 5$ annihilates every module in $\grassSS$, and hence
annihilates all modules in the closure, whereas $\grass{(\SS')}$ contains
an indecomposable module not annihilated by this arrow.  Similarly one
shows that the closure of
$$\grass{(\SS'')} : = \grass{(S_1^2, S_2 \oplus S_6, S_3 \oplus S_5, S_4,
0, 0)}$$ in $\grasstd$ is an irreducible component of class $2$ in
$\grasstd$.  Indeed, one notes that $\SS' < \SS''$ and that the closure of
$\grass{(\SS')}$ in $\grasstd$ is the only component of class $1$; then one
checks that $\grass{(\SS'')}$ is not contained in the closure of
$\grass{(\SS')}$.  By contrast: If we move on to the closures of the
corresponding irreducible subsets of $\toptd$ in $\modlad$, we obtain three
distinct irreducible components of class $0$ of $\modlad$, since
$(S_1)^2$ is minimal among the tops of the modules with dimension vector
$(2,1,1,1,1,1)$.
\smallskip  

\noindent {\rm {(3)}} For every nonnegative integer
$h$, there exists a truncated path algebra $\la$ of Loewy length $L+1 = 3$,
a dimension vector ${\bold d} = (d_1, d_2, \dots, d_n)$ of total
dimension $d = \sum_i d_i$, next to semisimple modules
$T{(0)} < T{(1)} < \cdots < T{(h)}$, such that the varieties
$\grass_{{\bold d}}^{T{(i)}}$ are irreducible, and hence, are
irreducible components of class zero of the varieties
$\grass_{d}^{T{(i)}}$  (for the connected components
$\grass_{{\bold d}}^{T{(i)}}$, see Observation 2.3).  On the other
hand, one can arrange for the closures of the corresponding
$\Top_{{\bold d}}^{T{(i)}}$ in
$\modlad$ to be irreducible components of class $i$ in $\modlad$,
respectively.  

To give a specific example, let $\la$ have quiver $Q$ as follows:

\ignore
$$\xymatrixrowsep{2pc}\xymatrixcolsep{4pc}
\xymatrix{   1 \ar[r]  &2 \ar[r]<0.75ex>  &3 \ar[l]<0.75ex> }$$
\endignore

\noindent and Loewy length $L+1 = 3$.  Set $T{(i)} = (S_1)^{h} \oplus
(S_3)^{i}$ for $0 \le i \le h$. Moreover, take $d = 3h$ and
${\bold d} = (h,h,h)$.  Observe that
$\grass_{{\bold d}}^{T(i)}$ equals
$\grass{\bigl( T(i), (S_2)^h, (S_3)^{h-i} \bigr)}$ and is therefore
irreducible by Theorem 5.3.  The closure of the corresponding irreducible
subset
$\Top_{{\bold d}}^{T{(i)}}$ of $\modlad$ is an irreducible component of
$\modlad$ of class $i$:  Clearly, the closure of
$\Top_{{\bold d}}^{T{(0)}}$ is an irreducible component of class
$0$ of $\modlad$, the only one in $\mod{_{{\bold d}}}$ in fact.  To see
that the closure of
$\Top_{{\bold d}}^{T{(1)}}$ in
$\modlad$ is an irreducible component of class $1$, it suffices to check
that $\Top_{{\bold d}}^{T{(1)}}$ is not contained in the closure of
$\Top_{{\bold d}}^{T{(0)}}$.  To see this, let $M$ be the unique
indecomposable module with radical layering $(S_1
\oplus S_3, S_2)$ and note that  $(\la e_1)^{h-1} \oplus M$ belongs to
$\Top_{{\bold d}}^{T{(1)}}$, but not to the closure of
$\Top_{{\bold d}}^{T{(0)}}$ in $\modlad$; indeed, the arrow 
$3 \rightarrow 2$ annihilates each module in
$\Top_{{\bold d}}^{T{(0)}}$, and consequently annihilates each module
in the closure.  Continue inductively. \qed
\enddefinition    

As we will show next, for a truncated path algebra, the  irreducible
varieties $\grassSS$ are bundles with affine fibres over a projective base
space.  The projective portion,
$\ggrassSS$, is recognized as a close kin to a flag variety.

Suppose that $V$ is an algebraic variety.  By a {\it Grassmann bundle over
$V$\/} we mean a fibre bundle over $V$ with fibre $F$, where
$F$ is a direct product of classical Grassmannians 
$\Grfrak(m_i, K^{n_i})$ and all of the pertinent maps are morphisms of
varieties; in particular, this means that the transition maps corresponding
to a suitable trivialization are automorphisms of $F$.  Moreover, we
call a bundle $\Delta$ over $V$ an {\it iterated Grassmann bundle\/} in
case there are bundles $\Delta_1, \dots, \Delta_r = \Delta$ such that
$\Delta_1$ is a Grassmann bundle over $V$ and each $\Delta_{i+1}$ is a
Grassmann bundle over $\Delta_i$.  In particular, the flag variety of any
finite dimensional vector space $W$ is an iterated Grassmann bundle over
$\Grfrak(\dim W - 1, W)$, and iterated Grassmann bundles may be viewed as
generalized flag varieties.  In a similar vein, a fibre bundle over
$V$ is referred to as an {\it affine bundle with fibre $F$\/} if $F$ is an
affine variety and, once again, all of the corresponding maps are
morphisms, resp., automorphisms of varieties.

\proclaim{Theorem 5.9\. Structure of the $\grassSS$} Suppose that $\la$ is
a truncated path algebra and $\SS = (\SS_0, \dots, \SS_L)$ a semisimple
sequence with $\grassSS \ne
\varnothing$.

Then $\grassSS$ is an affine bundle over the projective variety
$\ggrassSS$ with fibre $\AA^{N_1}$, where
$N_1 = N_1(\SS)$ is the invariant of $\SS$ introduced in  Observation {\rm
3.10}.  The base space of this bundle,
$\ggrassSS$, is an iterated Grassmann bundle over a direct product of
classical Grassmannian varieties.  More precisely, the base space of
$\ggrassSS$ is isomorphic to $\grass{(\SS_0, \SS_1, 0, \dots, 0)}$. 
\endproclaim

\demo{Proof} As before, we denote by $T$ the top $\SS_0$ of $\SS$, and  by
$P = \bigoplus_{l \ge 0} P_l$ a projective cover of $T$, endowed with the
natural grading and its distinguished sequence $z_1,
\dots, z_\t$ of top elements.   To any point $C \in \grassSS$ we assign the
following point in
$\ggrassSS$: 
$$C_{\h} = \bigoplus_{l \ge 1} C_l \ \ \ \text{with} \ \ \ C_l = P_l \cap
(C + \bigoplus_{m >l} P_m).$$   Suppose $\S$ is a skeleton compatible with
$\SS$. Throughout this argument, we will write the paths in our skeleton 
simply in the form $p\in \sigma$; in other words, $p$ stands for $p=
p'\zhat_r$ in $\Phat$ and also for $p'z_r$ in $P$ whenever we work in the
projective
$\la$-module
$P= \Phat/I\Phat$.

Let $C
\in \grassS$. Then $C_{\h}$ takes on the following form in the standard
affine coordinates for $\grassS$:   
$$C_{\h} = \bigl(c(\alpha p, q) \bigr)_{\alpha p  \ \S\text{-critical},\ q
\in \S_0(\alpha p)}$$  where $\S_0(\alpha p)$ is the set of all paths in
$\S(\alpha p)$ which have the same length as $\alpha p$. Analogously, we
denote by $\sigma_1(\alpha p)$ the set of all paths in
$\sigma(\alpha p)$ which are strictly longer than $\alpha p$. Recall that
$N_1$ is the cardinality of the disjoint union of the sets $\{\alpha
p\}\times
\S_1(\alpha p)$.

It is readily checked that the assignment
$$\pi: \grassSS \rightarrow \ggrassSS, \ \ \ C \mapsto C_{\h}$$ is a
morphism of varieties and that the inclusion $\ggrassSS
\hookrightarrow \grassSS$ is a section of $\pi$.  We will show that
$\pi$ makes $\grassSS$ an affine bundle with fibre $\AA^{N_1}$ over
$\ggrassSS$.  To do so, we specify trivializations over the open affine
subvarieties $\ggrassS$ covering $\ggrassSS$, where $\S$ traces the skeleta
compatible with
$\SS$.  Fixing such a skeleton, we define a morphism 
$$\ggrassS \times \AA^N_1 \rightarrow \pi^{-1}(\ggrassS)$$ as follows: 
Given a homogeneous point $D \in
\ggrassS$ with affine coordinates 
$$\bigl(d(\alpha p, q) \bigr)_{\alpha p\
\S\text{-critical},\ q \in \S_0(\alpha p)}$$  we send the pair 
$$\biggl(D, \bigl(c(\alpha p, q) \bigr)_{\alpha p\
\S\text{-critical},\ q \in \S_1(\alpha p)} \biggr)$$ to the following
submodule of $P$: 
$$\sum_{\alpha p\ \S\text{-critical}} \la \biggl(\alpha p
 - \sum_{q \in \S_1(\alpha p)} c(\alpha p, q)\, q
  - \sum_{q \in \S_0(\alpha p)} d(\alpha p, q)\, q \biggr).$$ By Theorem
5.3, arbitrary choices of coefficients $c(\alpha p, q)$ in $K^{N_1}$ lead
to such points in $\pi^{-1}(\ggrassS)$, independently of the given
homogeneous point $D$.  This makes the above assignment a well-defined
isomorphism of varieties.

To verify that the transition maps are  automorphisms of affine
$N_1$-space, suppose $D \in \ggrassS \cap
\ggrass{\tilde{\S}}$ with affine coordinates 
$$\bigl(d(\alpha p, q) \bigr)_{\alpha p\
\S\text{-critical},\ q \in \S_0(\alpha p)} \ \ \ \text{and}\ \ \   
\bigl(\tilde{d}(\alpha p, q) \bigr)_{\alpha p\
\S\text{-critical},\ q \in \S_0(\alpha p)}$$ relative to $\S$ and
$\tilde{\S}$, respectively.  Moreover, let $C \in \grassSS$ be a point in
the fibre above $D$.  Then $C \in
\grassS \cap \grass{(\tilde{\S})}$, since any homogeneous point $D$ has the
same set of skeleta as the points in $\pi^{-1}(D)$ (by an abuse of
language, we refer to skeleta of $P/D$ also as skeleta of $D$).  Let 
$$\bigl( c(\alpha p, q)\bigr)_{\alpha p \ \S\text{-critical},\ q \in
\S_1(\alpha p)}\ \  \text{and} \ \  \bigl( \tilde{c}(\alpha p,
q)\bigr)_{\alpha p \ \S\text{-critical},\ q \in \S_1(\alpha p)}$$  be the
supplementary nonhomogeneous coordinates relative to $\S$ and
$\tilde{\S}$, respectively.  For reasons of symmetry, it suffices to show
that the $\tilde{c}(\alpha p, q)$ are polynomials in the $c(\alpha p, q)$
with coefficients in $K \bigl( d(\alpha p, q) \bigr)$, if we treat the
$c(*)$'s, $\tilde{c(*)}$'s and
$d(*)$'s as independent variables over $K$.

Let $l$ be any integer between $1$ and $L$ and $\sigma_l =
\{p_1, \dots, p_s \}$; the set $\tilde{\S}$ has the same cardinality, say
$\tilde{\sigma}_l = \{ \tilde{p_1}, \dots, \tilde{p_s} \}$.  Then each
element $\tilde{p_i}$ in $P$ can be expanded, modulo $C$, in the format 
$$\tilde{p_i}
\equiv
\sum_{1
\le j \le s} a_{ij} p_j + A_i,  \tag{CONG}$$ where the $a_{ij}$ form an
invertible $s \times s$-matrix over
$K[d(*)]$, and each $A_i$ is a linear combination of paths in $\S$ of
lengths exceeding $l$, with coefficients in $K[c(*),d(*)]$; congruence
means congruence modulo $C$, and again we identify $\S$ and
$\tilde{\S}$ with subsets of $P$.  Solving for the elements 
$p_i \in P$ yields the latter as linear combinations
$$p_i \equiv
\sum_{1
\le j
\le s} b_{ij} \tilde{p}_j + B_i,  \tag{CONG-$l$}$$ modulo $C$, where the
$b_{ij}$ are coefficients in $K\bigl(d(*)\bigr)$ and $B_i$ is a linear
combination of paths in $\S$ which are longer than $l$, with coefficients
in $K\bigl(d(*)\bigr)[c(*)]$.  Now suppose that
$\tilde{\alpha} \tilde{p}$ is a $\tilde{\S}$-critical path with $\tilde{p}
\in \tilde{\S}_l$; say $\tilde{p} =
\tilde{p_1}$.  On multiplying the first of the $s$ congruences labeled
$(\text{CONG})$ from the left by
$\tilde{\alpha}$, we expand, modulo $C$, the element
$\tilde{\alpha} \tilde{p}$ of $P$ in terms of $\S$.  Then we successively
insert the congruences $\bigl(\text{CONG-}(l+1) \bigr)$,
$\text{CONG-}(l+2)$,
$\dots$ into the expansion.  This process terminates, because paths longer
than $L$ vanish.  It thus displays
$\tilde{\alpha} \tilde{p}$, modulo $C$, as a linear combination of terms
$\tilde{q}$ from $\tilde{\S}$ with coefficients in
$K\bigl(d(*)\bigr)[c(*)]$.  Comparison of coefficients shows the 
$\tilde{c}_{\tilde{\alpha} \tilde{p}}$ to have the required form.  This
completes the proof of the first assertion of the theorem.
\smallskip 

To prove the statement concerning $\ggrassSS$, we start by noting that
nonemptiness of $\grassSS$ implies that 
$$\ggrass{\SS_0,\SS_1, \dots, \SS_l,0, \dots, 0} \ne
\varnothing$$  for all $l \le L$. In a first step, we will show that
$\ggrass{\SS_0,\SS_1, 0, \dots, 0}$ is a direct product of classical
Grassmannians.  Indeed, this variety consists of all homogeneous submodules
$C = \bigoplus_{1 \le l \le L} C_l$ of $JP$ which are of the form $\,C_1
\oplus J^2 P = C_1 \oplus \bigoplus_{l \ge 2} P_l\,$ such that $\dim e_i
C_1 = \dim e_i P_1 - \dim e_i \SS_1$.  Moreover, we observe that, for any
$K$-subspace $U$ of $P_1$, the space
$U \oplus \bigoplus_{l \ge 2} P_l$ is a $\la$-submodule of $JP$.  Set
$m_i= \dim e_i \SS_1$ and $n_i= \dim e_i P_1$, that is, $\SS_1 =
\bigoplus_{1 \le i \le n} S_i^{m_i}$ and $JP/J^2P =
\bigoplus_{1 \le i \le n} S_i^{n_i}$.  Clearly, the map
$$\psi:  \prod_{1 \le i \le n} \Grfrak(n_i - m_i, P_1)
\rightarrow \ggrass{\SS_0, \SS_1, 0, \dots, 0}$$                  which
sends $(U_i)_{i \le n}$ to $\bigl( \bigoplus_{1
\le i \le n} U_i \bigr) \oplus \bigl( \bigoplus_{l \ge 2} P_l \bigr)$ is an
isomorphism of varieties, and the initial claim is established.

Finally, we prove that, for any integer $l \ge 2$ which is smaller than
$L$, the following map 
$\psi_l$ endows $\ggrass{\SS_0, \dots,
\SS_{l+1}, 0, \dots, 0}$ with the structure of a Grassmann bundle over 
$\ggrass{\SS_0, \dots, \SS_{l}, 0, \dots, 0}$: 
$$\psi_l: \ggrass{\SS_0, \dots, \SS_{l+1}, 0, \dots, 0} \rightarrow 
\ggrass{\SS_0, \dots, \SS_{l}, 0, \dots, 0}$$ sends $C' = \bigoplus_{1 \le
j \le L} C'_j$ = $\bigl( \bigoplus_{1 \le j
\le l+1} C'_j \bigr) \oplus \bigl( \bigoplus_{l+2 \le j \le L} P_j
\bigr)$ to 
$$C' + J^{l+1}P = \bigl( \bigoplus_{1 \le j \le l} C'_j \bigr)
\oplus \bigl( \bigoplus_{l+1 \le j \le L} P_j \bigr).$$   Next we provide a
trivialization over each patch of the open affine cover $(\ggrassS)_{\S}$
of $\ggrass{\SS_0, \dots, \SS_{l}, 0,
\dots, 0}$, where $\S$ traces the skeleta compatible with the semisimple
sequence $(\SS_0, \dots, \SS_{l}, 0, \dots, 0)$. Let
$\S = \bigcup_{1 \le r \le \t} \S^{(r)}$ be such a skeleton, suppose
$\SS_{l+1} = \bigoplus_{1 \le i \le n} S_i^{u_i}$, and consider the
subspaces
$$V_i = \sum_{1 \le r \le \t}\  \sum_{p \in \S_l^{(r)}} K \langle e_i Q_1 p
\rangle$$  of $P$, where $Q_1$ is the set of arrows in the quiver $Q$, and
$K\langle A \rangle$ the subspace generated by a subset $A$ of $P$ (again
we identify $\S$ with a subset of $P$ whenever called for).  If $v_i = \dim
V_i$, then $u_i \le v_i$ in view of the fact that 
$\ggrass{\SS_0, \dots, \SS_{l+1}, 0, \dots, 0} \ne \varnothing$. As we will
see, the fibre of $\psi_l$ over any point $C \in \ggrassS$ is isomorphic to
$$F := \prod_{1 \le i \le n} \Grfrak(v_i - u_i, V_i).$$ Note that $\sum_{1
\le i \le n} V_i = \bigoplus_{1 \le i \le n} V_i
\subseteq P_{l+1}$.  This setup permits us to describe a trivialization of
$\psi_l$ over $\ggrassS$ as follows:
$$\tau_{\S}: F \times \ggrassS \longrightarrow \psi_l^{-1}\bigl(
\ggrassS
\bigr)$$ with 
$$\bigl( (U_i)_{1 \le i \le n}\, , C \bigr) \mapsto C \oplus K\langle Q_1
C_l \rangle \oplus \bigl( \bigoplus_{1 \le i \le n} U_i \bigr)
\oplus\bigl( \bigoplus_{j \ge l+2} P_j\bigr),$$  where $K\langle Q_1 C_l
\rangle$ is the subspace of $P$ generated by all elements $\alpha c$, where
$\alpha$ is an arrow and $c \in  C_l$.  To ascertain well-definedness,
start by noting that the $K$-subspaces which are attained as images under
$\tau_{\S}$ are actually
$\la$-submodules of $JP$.  To see that each of these submodules  belongs to
$\grass{(\SS_0, \dots, \SS_{l+1}, 0, \dots, 0)}$, observe that
$$e_i P_{l+1} =  V_i \oplus  K\langle  e_i Q_1 C_l \rangle,$$ since the
layered graph of $P$ relative to $z_1, \dots, z_r$ is a forest with the
same tree sitting underneath each vertex labeled $e_i$ in the
$l$-th layer of this graph.  Therefore the codimension of $U_i \oplus
K\langle  e_i Q_1 C_l \rangle$ in $e_i P_{l+1}$ is $u_i$ as desired.  That
application of the map $\psi_l$ to any point in the image of $\tau_{\S}$
yields a point $D \in
\ggrass{\SS_0, \dots, \SS_{l}, 0, \dots, 0}$ with the property  that $P/D$
has skeleton $\S$, is obvious from our construction.  More strongly,
$\psi_l \circ \tau_{\S} \bigl( (U_i)_{i \le n}, C \bigr) = C$.  Thus
$\tau_{\S}$ is an isomorphism with inverse $C' \mapsto \bigl( (e_i
C_{l+1})_i\, , \psi_l(C') \bigr)$. 

Concerning compatibility of these trivializations:  Given two skeleta
$\sigma$ and $\tilde{\sigma}$ compatible with $(\SS_0, \dots, \SS_l, 0,
\dots, 0)$, it is routine to check that the corresponding transition map
$\pr_F \circ \tau_{\S}^{-1} \circ \tau_{\tilde{\S}} ( -, C)$ for $C
\in \ggrassS \cap \ggrass{\tilde{\S}}$ is an isomorphism of $F$, where
$\pr_F: F \times \ggrassS \rightarrow F$ denotes the projection onto the
fibre. \qed 
\enddemo

In general, the bundles of Theorem 5.9 are nontrivial.  A small example
exhibiting nontriviality is as follows.

\definition{Example 5.10} Let $\la = KQ/ I$, where $Q$ is the quiver

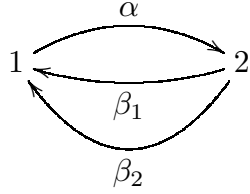
\ignore
$$\xymatrixrowsep{1pc}\xymatrixcolsep{5pc}
\xymatrix{   1 \ar@/^1pc/[r]^{\alpha}  &2 \ar@/^/[l]^{\beta_1}
\ar@/^2.25pc/[l]^{\beta_2} }$$
\endignore

\noindent and $I$ is generated by the paths of length $3$.  Let $T = S_1
\oplus S_2$, $d = 4$, and $\SS = (T, S_2, S_1)$.  Then it is readily
checked that
$\grassSS$ is a nontrivial $\AA^2$-bundle over $\ggrassSS \cong \PP^1$. On
the other hand, if $\SS' = (T, S_1, S_2)$, then $\grass{(\SS')}$ is
isomorphic to the trivial $\AA^1$-bundle over $\PP^1$.  This example will
be revisited after Theorem 5.12.  
\qed\enddefinition  
                                   
We follow with properties of the $\SS$-generic module for any sequence
$\SS$.  As in Section 4, we assume that $K$ has infinite transcendence
degree over its prime field. Moreover, we let $\Kcirc
\subseteq K$ be the algebraic closure of the prime field.  Clearly, any
truncated path algebra is defined over
$\Kcirc$ so that the prerequisites of Section 4 are in place.

The following lemma shows in particular that all syzygies of
$\la$-modules are direct sums of cyclic modules, each isomorphic to a left
ideal of $\la$ generated by a path. (That this holds for second syzygies
already follows from \cite{\domino, Theorem A}.)  

\proclaim{Lemma 5.11} Suppose $\la$ is a truncated path algebra and $C
\in \grasstd$. Then the syzygy $C$ of $P/C$ is a direct sum of cyclic
modules, each of which is isomorphic to a left ideal of $\la$. More
precisely, if $C\in \grassS$, the elements 
$$\omega_{\alpha p,r} :=\alpha pz_r - \sum_{q\zhat_s\in \sigma(\alpha
p\zhat_r)} c(\alpha p\zhat_r,q\zhat_s)\, qz_s\ \in \ P,$$ 
where $\alpha
p\zhat_r$ traces the
$\sigma\text{-critical}$ paths and the 
$c(\alpha p\zhat_r,q\zhat_s)$ are the affine coordinates of $C$ in
$\grassS$ {\rm {(}}see Section {\rm{3.C)}}, generate nonzero independent
cyclic submodules of
$P$ such that
$$C = \bigoplus_{\alpha p\zhat_r\ \S\text{-critical}} \la \omega_{\alpha
p,r}
\qquad
\text{and}\qquad  \la \omega_{\alpha p,r} \cong \la \alpha p\qquad 
\text{for each\ } \S\text{-critical path\ } \alpha p\zhat_r.$$
\endproclaim

\demo{Proof} Suppose $C\in
\grassS$.  We already know that $C$ is generated by the $\omega_{\alpha
p,r}$; see the remarks preceding Theorem 3.8.   By definition, all of the
$\omega_{\alpha p,r}$ are nonzero in
$P$, as the paths involved in the pertinent expansions have lengths $\le
L$.  Assume, to the contrary of our claim, that $0$ is  a sum of certain
nontrivial $\la$-multiples of
$\omega_{\alpha p,r}$'s in $P$.  Let
$$\sum_{1\le i\le \nu} \lambda_i\, \omega_{\alpha_i p_i,r_i} =0 
\tag{$\dagger$}$$  
be such a nontrivial dependence relation in $P$, such
that $\nu$ is minimal and the $\lambda_i$ are $K$-linear combinations of
paths in $\la$, the lengths of which are bounded by the differences $L -
\len(\alpha_i p_i)$, respectively; the latter assumption is legitimate as,
by construction, each path in $\S(\alpha_i p_i \zhat_{r_i})$ is at least as
long as $\alpha_i p_i$.  To reach a contradiction, we use the
$K$-linear independence of the formally distinct elements of the form $p
z_r$ in $P$, where $p$ is a path of length at most $L$ in $\la$ which
starts in the vertex $e(r)$ that norms the top element $z_r$.   Let $l_0$
be the minimum of the lengths of the
$\alpha_i p_i$, and assume that $\alpha_1 p_1$ has length $l_0$.  Next,
observe the following: If $\lambda_1= \sum_j k_j q_j$, where the $k_j$ are
nonzero elements of $K$ and the $q_j$ are distinct paths with $0\le
\len(q_j) \le L-l_0$, then the occurrence of the path
$q_1 \alpha_1 p_1 z_{r_1}$ in the term $\lambda_1 \alpha_1 p_1 z_{r_1}$ is
the only one in $(\dagger)$.   Indeed, as we already remarked, the paths
$p\zhat_r$ in the union of the
$\sigma(\alpha_i p_i \zhat_{r_i})$ are at least as long as $\alpha_1 p_1$,
and they all belong to $\sigma$; since
$\alpha_1 p_1 \zhat_{r_1}$ does not belong to $\sigma$, we infer that
$\alpha_1 p_1 z_{r_1}$ does not arise in the $K$-linear path expansion in
$P$ of any
$\la$-multiple of a path $qz_s$ with $q\zhat_s \in
\bigcup_{1\le i\le \nu} \sigma(\alpha_i p_i \zhat_{r_i})$.  That the path
expansions of the
$\la$-multiples  of the $\alpha_i p_i z_{r_i}$ with $2 \le i\le \nu$ in $P$
do not contain nontrivial $K$-multiples of $q_1 \alpha_1 p_1 z_{r_1}$
either is clear, because the minimal initial subpath not belonging to $\S$
of any
$q \alpha_i p_i \zhat_{r_i}$ is $\alpha_i p_i \zhat_{r_i}$ and thus
different from
$\alpha_1 p_1 \zhat_{r_1}$ for $i > 1$ (due to minimality of
$\nu$). This yields the required contradiction and thus proves independence
of the cyclic modules
$\la
\omega_{\alpha p,r}$.

To see that each $\la \omega_{\alpha p,r} \subseteq P$ is isomorphic to
some left ideal of $\la$, let $e_s$ be the endpoint of $\alpha$. Then one
readily checks that $\la \omega_{\alpha p,r} \cong \la \alpha p$, in view
of the fact that all paths $q\zhat_s\in \sigma(\alpha p\zhat_r)$ are at
least as long as $\alpha p\zhat_r$ and end in the vertex $e_s$ as well.
\qed\enddemo

\definition{Remark} Lemma 5.11 clearly extends to syzygies of nonfinitely
generated modules, the same argument being applicable.  For that purpose,
one generalizes the concept of a skeleton to the infinite dimensional case
and adapts the definitions of
$\S$-critical paths and the scalars $c(\alpha pz_r,qz_s)$ in the obvious
way.  This strengthening is immaterial here, but is used in \cite{\DHL}.
\enddefinition    

Most of the statements of the following theorem are immediate consequences
of Theorems 4.3 and 5.3 and their proofs; for the sharper result concerning
generic syzygies under (1), which is key to easy applicability, we use
Lemma 5.11.  Note in particular that, for  our present choice of $\la$, the
$\SS$-generic modules can be read off the quiver, without requiring any
computation, since the auxiliary varieties
$\grassS$ are known to be affine spaces in this case.  The uniqueness
assertion under (1) also improves on the corresponding one in Theorem 4.3;
the considerably stronger version holding in the truncated situation is
apparent from the high level of symmetry of the hypergraph of $G(\SS)$,
again a consequence of the projective presentation of
$G(\SS)$ given below.  Concerning part (2) of the upcoming theorem: The
results of Theorem 4.3 and Corollary  4.7 for $\SS$-generic modules carry
over to
$\SS$-generic graded modules, if we replace ``isomorphism" by
``graded-isomorphism".  Mutatis mutandis, the final assertion of part (1)
is true for $\SS$-generic graded modules as well.    

We keep the notation of the general Theorem 4.3, but recall that it can be
simplified for a truncated path algebra, due to the well-definedness of
path length in $\la$. Indeed, it is now harmless to view skeleta as subsets
of the distinguished projective cover $P$, whence the $\zhat_r$ may be
replaced by $z_r$ throughout.

\proclaim{Theorem 5.12\. The $\SS$-generic modules and the
$\SS$-generic graded modules over a truncated path algebra}  Suppose
$\la$ is a truncated path algebra, and let $\SS$ be any semisimple sequence
with $\grassSS \ne \varnothing$.  As before, $P$ is the distinguished
projective cover of the top $T = \SS_0$.
\medskip

\noindent {\rm (1)} The $\SS$-generic module $G(\SS)$ for the irreducible
variety $\grassSS$  --  equivalently, for the classical variety $\laySS$ 
--  has a projective presentation as follows.  Let
$\S$ be any skeleton compatible with
$\SS$ and $N$ the {\rm(}disjoint{\rm)} union of the sets $\{\alpha
p\zhat_r\} \times
\S(\alpha p\zhat_r)$. Moreover, let $\bigl( x(\alpha p\zhat_r,q\zhat_s)
\bigr)$ be any family of scalars, indexed by $N$, which is algebraically
independent over
$\Kcirc$.  Then  $G(\SS) = P/C(\SS)$, up to equivalence, where
$C(\SS)$ is the  submodule of $P$ generated by the differences
$\alpha pz_r - \sum_{q\zhat_s \in \S(\alpha p\zhat_r)}  x(\alpha
p\zhat_r,q\zhat_s)\, qz_s$, with $\alpha p\zhat_r$ tracing the $\S$-critical
paths. {\rm (By definition,
$\S(\alpha p\zhat_r)$ consists of those paths in $\S$ which have lengths
between $\len(\alpha p)$ and $L$ and end in the same vertex as $\alpha p$.)}
In particular, the syzygy
$$C(\SS)= \bigoplus_{\alpha p\zhat_r\, \sigma\text{-critical}} \la
\biggl( \alpha pz_r - \sum_{q\zhat_s \in \S(\alpha p\zhat_r)}  
x(\alpha p\zhat_r,q\zhat_s)\ qz_s \biggr)\ \ \ \ \cong \ \
\bigoplus_{\alpha p\zhat_r\
\S\text{-critical}} \la \alpha p$$  of $G(\SS)$ is a direct sum of cyclic
modules, each of which is isomorphic to a left ideal of $\la$ generated by
a path.  It is completely determined by $\SS$, up to isomorphism {\rm
{(}}not only up to $\Gal$-equivalence{\rm {)}}.

Moreover, any submodule $H$ of $G(\SS)$, which  --  modulo $C(\SS)$  -- is
generated by homogeneous elements of a fixed degree in $P$, is the
$\SS(H)$-generic module, up to equivalence.
\medskip

\noindent {\rm (2)} The $\SS$-generic graded module $\GrG(\SS)$ generated
in degree $0$ has a projective presentation as follows. Denote by $N_0$ the
{\rm(}disjoint{\rm)} union of the sets $\{\alpha p\zhat_r\} \times
\S_0(\alpha p\zhat_r)$, where $\S_0(\alpha p \zhat_r)$ consists of those
paths in
$\S(\alpha p \zhat_r)$ which have the same length as $\alpha p$.
 Let
$\S$ be any skeleton compatible with
$\SS$ and $\bigl( x(\alpha p\zhat_r,q\zhat_s) \bigr)$ a
$\Kcirc$-algebraically independent family of scalars indexed by $N_0(\SS)$. 
Then  $\GrG(\SS)$ is equivalent to
$P/(\GrC(\SS))$, where
$$\GrC(\SS)= \bigoplus_{\alpha p\zhat_r\, \sigma\text{-cricital}} \la
\biggl( \alpha p\zhat_r - \sum_{q\zhat_s \in \S_0(\alpha p\zhat_r)} x(\alpha
p\zhat_r,q\zhat_s)\ qz_s
\biggr)\  \ \cong \ \ C(\SS). \qquad \square$$ 
\endproclaim

In light of Corollary 4.7, Theorem 5.12 thus reduces structural problems
regarding the $\SS$-generic modules to combinatorial tasks.  In particular,
this is true for the questions regarding decomposability of the
$\SS$-generic module, the structure of its indecomposable summands, generic
socles, higher syzygies, etc.  For example, the module
$G(\SS)$ is indecomposable if and only if all of its (finitely many)
hypergraphs  relative to full sequences of top elements are connected. 

\definition{Example 5.10 revisited}  Let $d= 4$.  For the dimension vector
${\bold d} = (2,2)$, we list all irreducible components of $\modlad$
which are contained in the connected component $\mod_{{\bold d}}$ and
display their generic modules.  There are precisely four, and the generic
modules have top dimension $\le 2$.  We will start by displaying the
$\SS$-generic modules for each of the six nonempty varieties
$\grassSS$, where $\SS$ is a semisimple sequence of dimension vector
${\bold d}$ with $\dim \SS_0 \le 2$; these are: $\SS^{(1)}$ $=$
$(S_1^2, S_2^2, 0)$, $\SS^{(2)}$ $=$ $(S_2^2, S_1^2, 0)$, $\SS^{(3)}$ $=$
$(S_1 \oplus S_2, S_1 \oplus S_2, 0)$,
$\SS^{(4)}$ $=$ $(S_2, S_1^2, S_2)$, $\SS^{(5)}$ $=$ $(S_1 \oplus S_2, S_2,
S_1)$, and $\SS^{(6)}$ $=$ $(S_1 \oplus S_2, S_1, S_2)$.  Instead of
formally presenting the generic modules $G_i = G(\SS^{(i)})$ (which is
straightforward in view of Theorem 5.12), we provide hypergraphs relative
to suitable top elements and skeleta, which is more informative at a
glance.     

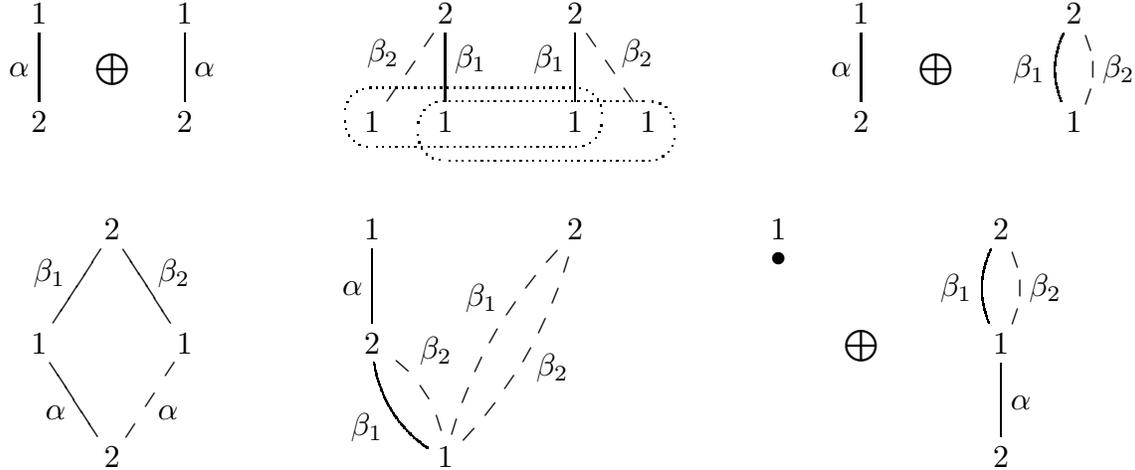
\ignore
$$\xymatrixrowsep{1.8pc}\xymatrixcolsep{1pc}
\xymatrix{  1 \edge[d]_{\alpha} & \ar@{}[d]|{\tsize\bigoplus} &1
\edge[d]^{\alpha}  &&&&2
\dashedge[dl]_{\beta_2} \edge[d]^(0.4){\beta_1} &&2
\edge[d]_(0.4){\beta_1} \dashedge[dr]^{\beta_2}  &&&&1 \edge[d]_{\alpha} &
\ar@{}[d]|{\tsize\bigoplus} &&2 \edge@/_/[d]_{\beta_1}
\dashedge@/^/[d]^{\beta_2}  \\ 2 &&2  &&&1 \save
[0,0]+(-3,4);[0,3]+(3,-2.5) **\frm<10pt>{.}\restore
 &1 \save [0,0]+(-3,2.5);[0,3]+(3,-4) **\frm<10pt>{.}\restore  &&1 &1  &&&2
&&&1  \\
 &2 \edge[dl]_{\beta_1} \edge[dr]^{\beta_2}  &&&&1 \edge[d]_{\alpha} &&&2
\dashedge@/_/[ddll]_(0.4){\beta_1} \dashedge@/^/[ddll]^{\beta_2}  &&&1
\drbl &&&2 \edge@/_/[d]_{\beta_1} \dashedge@/^/[d]^{\beta_2}    \\  1
\edge[dr]_{\alpha} &&1 \dashedge[dl]^{\alpha}  &&&2
\dashedge@/^/[dr]^(0.33){\beta_2} \edge@/_/[dr]_{\beta_1}   &&&&
&&&\bigoplus &&1 \edge[d]^{\alpha}    \\
 &2  &&&&&1 &&&& &&&&2  &&
 }$$
\endignore

\noindent  Concerning $G_6$: 
The displayed hypergraph does not correspond to a generic presentation of
$G_6$ in the sense of the paragraph preceding Corollary 4.5; the given
non-generic presentation more clearly exhibits decomposability.  

Let
$\C_i$ be the closure of $\lay{\SS^{(i)}}$ in $\modlad$ for $i = 1,
\dots, 6$.  Aided by Corollary 4.7, we will sift out the $\C_i$ which are
maximal irreducible in $\modlad$.  Clearly, $\C_3$ is contained in each of
$\C_4, \C_5, \C_6$.  Next, it is straightforward to construct a curve
$\psi: \AA^1 \rightarrow \modlad$ with $\psi(t) \in \lay{\SS^{(4)}}$ for
$t \ne 0$, such that $\psi(0)$ represents $G_6$; hence $\C_6 \subseteq
\C_4$.  Since the dimension vector of top$\,G_i$ for $i = 1, 4$ is minimal
among the dimension vectors of the tops of the modules in
$\C_1, \dots, \C_6$ (see Observation 2.1), $\C_1$ and $\C_4$ are
irreducible components of $\modlad$.  Comparing tops, we further note that
the only $\C_j$ which potentially contains $\C_2$ is $\C_4$; but, in light
of Corollary 4.7, the containment $\C_2 \subseteq \C_4$ is ruled out by the
fact that $S_2$ is evidently a summand of $\Soc G_4$, but not a summand of
$\Soc G_2$.  Again comparing socles of generic modules and invoking
Corollary 4.7, we conclude that $\C_5$ is not contained in any of
$\C_1, \C_2, \C_4$.  Thus the $\C_i$ for $i = 1, 2, 4, 5$ constitute a full
irredundant list of those irreducible components of
$\modlad$ which are contained in the connected component representing the
modules with dimension vector
${\bold d}$.  Consequently, the $G_i$ for $i = 1, 2, 4, 5$ are even
$4$-generic. In the top-order, $\C_1$ and $\C_4$ are irreducible components
of $\modlad$ of class $0$, while $\C_2$ and $\C_5$ are of class $1$.

The descriptions of the generic modules $G_i$ as in Theorem 5.11 allow
for their represen\-tation-theoretic evaluation.  For instance,  
generically, the modules in $\C_5$ have socle $S_1$ (this is not visible
from the given hypergraph of $G_5$, but is immediate from the
projective presentation).  In particular, the modules in
$\C_5$ are generically indecomposable.  Moreover, generically, they satisfy
$\dim \End_\la(M) = 2 = \dim \Ext^1_\la(M,M)$, have generic syzygy
isomorphic to $S_1 \oplus (\la e_1/ J^2 e_1)^2$ and generic projective
dimension $\infty$.  \qed
\enddefinition

We illustrate Theorem 5.11 with another, somewhat more complex, example. We
also display the generic graded module $\GrG(\SS)$ for the considered
semisimple sequence $\SS$.    

\definition{Example 5.13}  Let $Q$ be the quiver

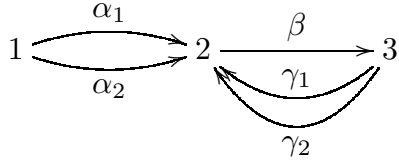
\ignore
$$\xymatrixrowsep{1pc}\xymatrixcolsep{4pc}
\xymatrix{   1 \ar@/^/[r]^{\alpha_1} \ar@/_/[r]_{\alpha_2} &2
\ar[r]^{\beta} &3
\ar@/^1.25pc/[l]_{\gamma_1} \ar@/^2pc/[l]^{\gamma_2} }$$
\endignore

\noindent and $\la$ the truncated path algebra $KQ/I$, where $I$ is
generated by the paths of length $4$; i.e., $L = 3$.  We will consider the
$14$-dimensional semisimple sequence 
$$\SS= (S_1^2 \oplus S_2 \oplus S_3, S_2^5 \oplus S_3, S_3^3, S_2).$$ We
give the hypergraph of the corresponding $\SS$-generic module
$G(\SS)$, relative to the skeleton $\S$ shown below; the broken edges again
indicate the $\S$-critical paths.

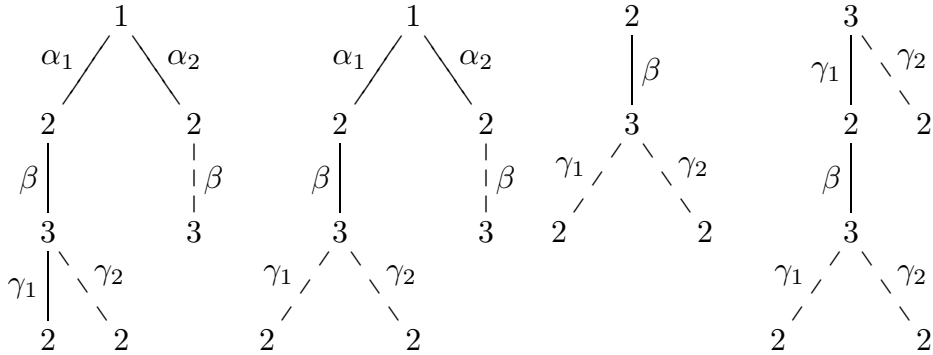
\ignore
$$\xymatrixrowsep{1.8pc}\xymatrixcolsep{1pc}
\xymatrix{ 
 &1 \edge[dl]_{\alpha_1} \edge[dr]^{\alpha_2}  &&&&1
\edge[dl]_{\alpha_1} \edge[dr]^{\alpha_2}  &&&2
\edge[d]^{\beta}  &&&3 \edge[d]_{\gamma_1}
\dashedge[dr]^{\gamma_2}  \\ 2 \edge[d]_{\beta} &&2
\dashedge[d]^{\beta}  &&2
\edge[d]_{\beta} &&2 \dashedge[d]^{\beta}  &&3
\dashedge[dl]_{\gamma_1} \dashedge[dr]^{\gamma_2}  &&&2 
\edge[d]_{\beta} &2  \\ 3 \edge[d]_{\gamma_1} \dashedge[dr]^{\gamma_2} &&3 
&&3
\dashedge[dl]_{\gamma_1} \dashedge[dr]^{\gamma_2} &&3  &2 &&2  &&3
\dashedge[dl]_{\gamma_1} \dashedge[dr]^{\gamma_2}  \\ 2 &2  &&2 &&2 && &&
&2 &&2 }$$
\endignore

\noindent The hypergraph of $G(\SS)$ results from superposition of the
following diagrams (1)--(4):

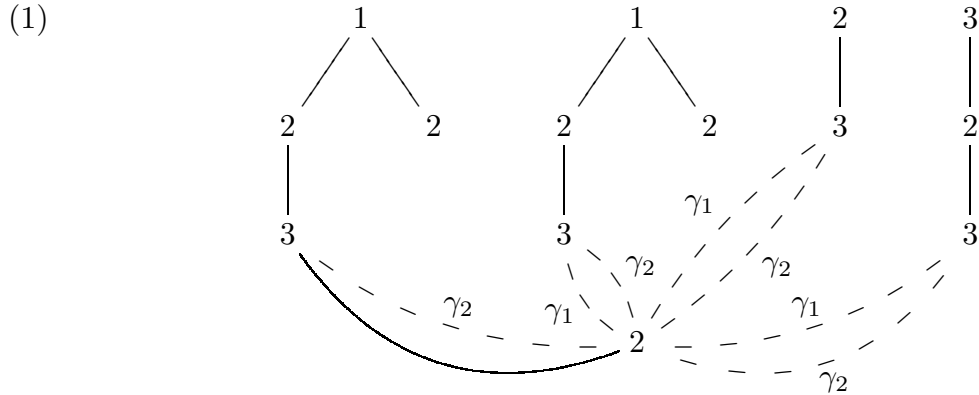
\ignore
$$\xymatrixrowsep{1.8pc}\xymatrixcolsep{1pc}
\xymatrix{ 
 &1 \edge[dl] \edge[dr]  &&&&1 \edge[dl] \edge[dr]  &&&2 \edge[d]  &&3
\edge[d]  \\ 2 \edge[d] &&2  &&2 \edge[d] &&2  &&3
\dashedge@/_/[ddlll]_{\gamma_1}
\dashedge@/^/[ddlll]^{\gamma_2}  &&2 \edge[d]  \\ 3
\dashedge@/_1pc/[1,5]^{\gamma_2} \edge@/_2pc/[1,5]  && &&3
\dashedge@/^/[dr]^{\gamma_2} \dashedge@/_/[dr]_{\gamma_1}  && && &&3
\dashedge@/^1pc/[1,-5]_{\gamma_1} \dashedge@/^2pc/[1,-5]^{\gamma_2} 
\\
 && && &2 }  \tag1	$$
\endignore

\ignore
$$\xymatrixrowsep{1.8pc}\xymatrixcolsep{1pc}
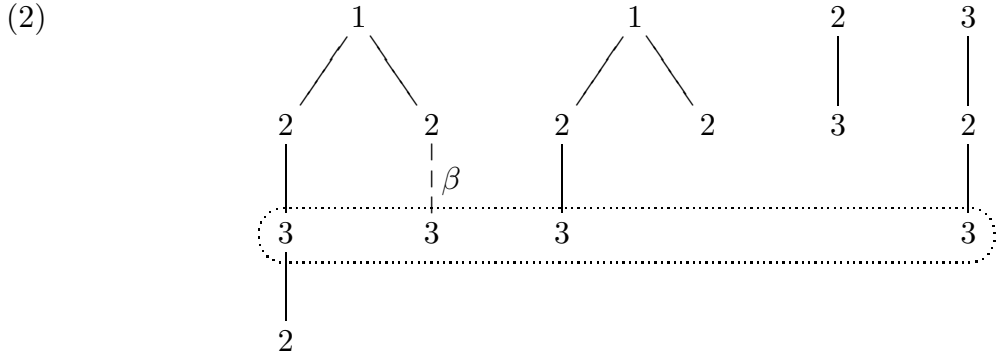
\xymatrix{ 
 &1 \edge[dl] \edge[dr]  &&&&1 \edge[dl] \edge[dr]  &&&2 \edge[d]  &&3
\edge[d]  \\ 2 \edge[d] &&2 \dashedge[d]^{\beta}  &&2 \edge[d] &&2  &&3 
&&2 \edge[d] 
\\  3 \edge[d] \levelpool{10}  &&3  &&3 && && &&3  \\ 2 }  \tag2 $$
\endignore

\ignore
$$\xymatrixrowsep{1.8pc}\xymatrixcolsep{1pc}
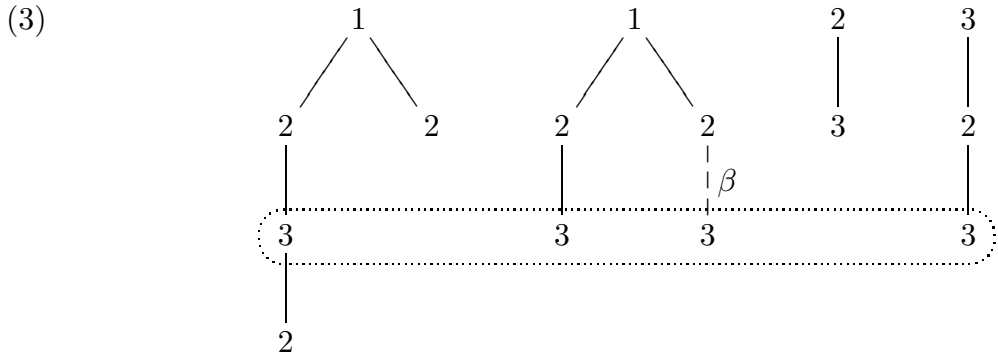
\xymatrix{ 
 &1 \edge[dl] \edge[dr]  &&&&1 \edge[dl] \edge[dr]  &&&2 \edge[d]  &&3
\edge[d]  \\ 2 \edge[d] &&2  &&2 \edge[d] &&2 \dashedge[d]^{\beta}  &&3 
&&2 \edge[d] 
\\  3 \edge[d] \levelpool{10}  &&  &&3 &&3 && &&3  \\ 2 }  \tag3 $$
\endignore

\ignore
$$\xymatrixrowsep{1.8pc}\xymatrixcolsep{1pc}
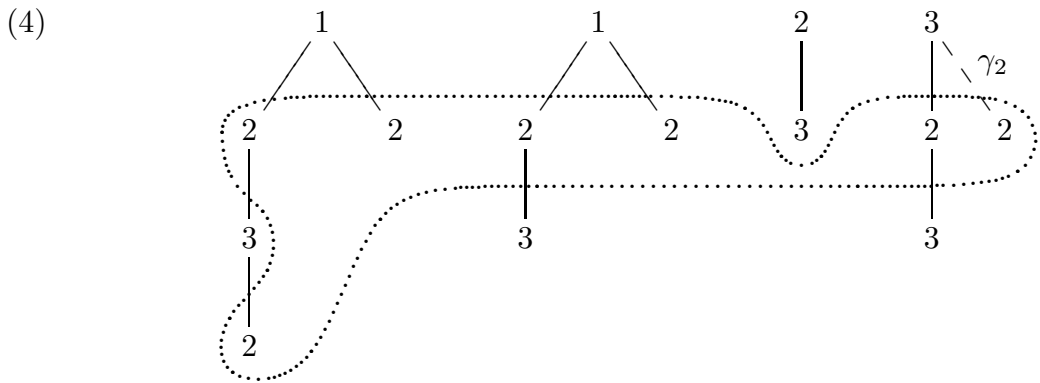
\xymatrix{ 
 &1 \edge[dl] \edge[dr]  &&&&1 \edge[dl] \edge[dr]  &&&2 \edge[d]  &&3
\edge[d] \dashedge[dr]^{\gamma_2}  \\ 2 \edge[d]
\save[0,0]+(-3,0);[0,0]+(-3,0) **\crv{~*=<2.5pt>{.} [0,0]+(-3,4)
&[0,1]+(3,4) &[0,6]+(-3,4) &[0,7]+(3,4) &[0,8]+(-3,-4) &[0,8]+(3,-4)
&[0,9]+(-3,4) &[0,10]+(0,4) &[0,11]+(4,4) &[0,11]+(4,-6) &[0,10]+(-3,-6)
&[0,3]+(3,-6) &[0,2]+(-3,-6) &[2,1]+(0,-4) &[2,0]+(-4,-4) &[2,0]+(-3,4)
&[1,0]+(3,-4) &[1,0]+(3,4) &[0,0]+(-3,-6) }\restore &&2  &&2
\edge[d] &&2  &&3  &&2 \edge[d] &2  \\   3 \edge[d]  &&  &&3 && && &&3  \\
2 && }  \tag4 $$
\endignore

\noindent The module $G = G(\SS)$ is indecomposable, as can already be
gleaned from Diagram (1); indeed, the full collection of scalars involved
in the relations, which are indicated by dashed edges, is algebraically
independent over $\Kcirc$.  The generic socle is
$S_2^2$; one copy of
$S_2$ in the socle is obvious, the other can be read off Diagram (3).  The
closure of $\grassSS$ is not an irreducible component of $\grasstd$, where
$T = S_1^2 \oplus S_2 \oplus S_3$ and $d = 14$; in fact, $G$ arises as a
top-stable degeneration of a generic module for the semisimple sequence
$\SS' = (S_1^2 \oplus S_2 \oplus S_3, S_2^5, S_3^4, S_2)$.  Moreover, a
purely combinatorial process yields the first syzygy of $G$ to be
$\Omega^1(G) \cong S_1^5 \oplus (\la e_2/ J^2 e_2)^3
\oplus (\la e_2/ J^3 e_2) \oplus (\la e_3/ J^2 e_3)^2$, which makes the
higher generic syzygies for the modules with radical layering $\SS$ readily
available.  In particular, $\pdim G = \infty$, which, in view of Corollary
4.7, shows that all modules with the given radical layering have infinite
projective dimension.

The $\SS$-generic graded module $\GrG(\SS)$ with top generated in a fixed
degree has the following modified hypergraph relative to $\S$:  It is the
superposition of the diagram below with Diagrams (2) and (3) above.

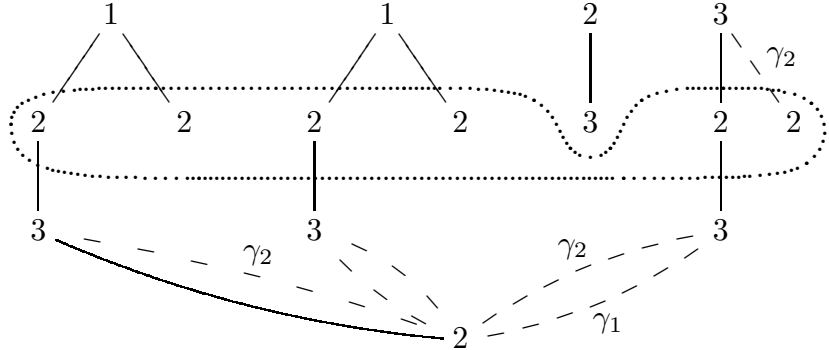
\ignore
$$\xymatrixrowsep{1.8pc}\xymatrixcolsep{1pc}
\xymatrix{ 
 &1 \edge[dl] \edge[dr]  &&&&1 \edge[dl] \edge[dr]  &&&2 \edge[d]  &&3
\edge[d] \dashedge[dr]^{\gamma_2}  \\ 2 \edge[d]
\save[0,0]+(-3,0);[0,0]+(-3,0) **\crv{~*=<2.5pt>{.} [0,0]+(-3,4)
&[0,1]+(3,4) &[0,6]+(-3,4) &[0,7]+(3,4) &[0,8]+(-3,-4) &[0,8]+(3,-4)
&[0,9]+(-3,4) &[0,10]+(0,4) &[0,11]+(4,4) &[0,11]+(4,-6) &[0,10]+(-3,-6)
&[0,1]+(3,-6) &[0,0]+(-3,-6) }\restore &&2  &&2 \edge[d] &&2  &&3  &&2
\edge[d] &2  \\  3 \dashedge@/^0.3pc/[1,6]^{\gamma_2} \edge@/_/[1,6]  && &&3
\dashedge@/_0.3pc/[drr] \dashedge@/^/[drr]   && && &&3 
\\
 && && &&2 \dashedge@/^/[urrrr]^{\gamma_2} \dashedge@/_/[urrrr]_{\gamma_1}
}$$
\endignore

\noindent Generically, the graded modules with radical layering $\SS$  have
two indecomposable summands, with dimension vectors $(0,1,1)$ and
$(2,6,4)$. The generic socle may shrink as one passes from the graded to
the ungraded situation; indeed, $\Soc \GrG(\SS) \cong \Soc G(\SS) \oplus
S_3$ in our example.  \qed
\enddefinition

In general, it is cumbersome to explicitly state combinatorial equivalent
conditions for indecomposability of $G(\SS)$ or $\GrG(\SS)$ in terms of
$\SS$, $Q$, and $L$.  We content ourselves with presenting a
straightforward necessary condition for the graded situation.  The
following auxiliary graph depends on a choice of skeleton, but the vertex
sets of its connected components do not.  The vertex set is the set $Z =
\{z_1,
\dots, z_{\t}\}$ of distinguished top elements in the projective cover
$P$ of $\SS_0$.  Given a skeleton
$\S$ compatible with $\SS$, there is an edge connecting $z_r$ and
$z_s$ if and only if either there exists a $\S$-critical path
$\alpha p\zhat_r$ and a path $q\zhat_s \in \sigma$ with $\len(\alpha p)=
\len(q)$ and
$\text{endpt}(\alpha p)= \text{endpt}(q)$, or else this condition holds
with the roles of
$r$ and
$s$ reversed.

\proclaim{Corollary 5.14\. Generic indecomposability}  Let $\SS$ be a
semisimple sequence over a truncated path algebra $\la$ such that $\grassSS
\ne
\varnothing$.  Moreover, let $Z_1, \dots, Z_{\mu} \subseteq Z$ be the
vertex sets of the connected components of any of the auxiliary graphs
introduced above, and let 
$U_i$, for $1 \le i \le \mu$, be the submodule of $\GrG(\SS) = P/\GrC(\SS)$
generated by the  residue classes $z + \GrC(\SS)$, $z \in Z_i$.

Then $\GrG(\SS) = \bigoplus_{1 \le i \le \mu} U_i$.  In particular,
indecomposability of $\GrG(\SS)$ implies connectedness of the auxiliary
graphs.

If the top of $\SS$ is squarefree,  the $U_i$ are the indecomposable direct
summands of $\GrG(\SS)$, and
$\GrG(\SS)$ is indecomposable if and only if any of the auxiliary graphs is
connected. \qed    
\endproclaim

We leave the easy proof to the reader.  There are obvious analogues of the
two statements of Corollary 5.14 for the ungraded situation.  The first is
always true, while the second is not.  The following general connection
between the graded and ungraded situations is helpful.

\proclaim{Corollary 5.15}  Let $\SS$ be a semisimple sequence over a
truncated path algebra $\la$ such that $\grassSS \ne \varnothing$. If the
generic graded module $\GrG(\SS)$ with radical layering $\SS$ is
indecomposable, then so is the generic ungraded module $G(\SS)$.
\endproclaim

\demo{Proof}  Note that the module $\GrG(\SS)$ (resp., $G(\SS)$) is
indecomposable precisely when every hypergraph of $\GrG(\SS)$ (resp.,
$G(\SS)$) is connected. Since the set of skeleta of $\GrG(\SS)$ coincides
with that of
$G(\SS)$, our claim is easily deduced from Theorem 5.12. \qed
\enddemo

The converse of Corollary 5.15 fails in general.  Indeed, let $\la = KQ$,
where $Q$ is the quiver
$$1 \longrightarrow 2 \longrightarrow 3 \longleftarrow 4$$ and $\SS$ the
sequence $(S_1 \oplus S_4, S_2, S_3)$.  Then the
$\SS$-generic module is indecomposable, while the $\SS$-generic graded
module decomposes. 

If ``structural symmetry" of a module is measured by the dimension of its
endomorphism ring, the generic module $G(\SS)$ (or $\GrG(\SS)$) has minimal
structural symmetry among the modules represented by $\SS$; see Corollary
4.7.  Yet, note that the endomorphism ring of an indecomposable
$\SS$-generic module $G(\SS)$ need not be trivial:  Let
$Q$ be the quiver of Examples 5.8(1), and take $L = 2$ and
$\SS = (S_1, S_1, S_2)$. Then
$\End_{\la}G(\SS)$ has dimension $2$.  By contrast, if indecomposable, the
generic graded module $\GrG(\SS)$ always has trivial homogeneous
endomorphism ring.

On the other hand, in terms of its submodule structure, the module $G(\SS)$
(resp\. $\GrG(\SS)$) displays maximal symmetry among the (graded) modules
with radical layering $\SS$, in the following sense.  Let $\S$ be any
skeleton compatible with
$\SS$, and $G(\SS) = P/ C(\SS)$ the corresponding generic presentation
described in Theorem 5.12.  Then the theorem shows in particular that,
given any two paths $p\zhat_r, q\zhat_s \in \S$ with coinciding length and
endpoint, the canonical images of $\la pz_r$ and $\la qz_s$ in
$G(\SS)$ are $\Gal$-equivalent and thus have the same hypergraphs.

\enddocument